\documentclass{amsart}
\usepackage{amssymb,amsmath,amsthm}
\usepackage{enumerate,color}

\newtheorem{theorem}{Theorem}
\newtheorem{lemma}[theorem]{Lemma}
\newtheorem{proposition}[theorem]{Proposition}
\newtheorem{corollary}[theorem]{Corollary}

\theoremstyle{definition}

\newtheorem{remark}[theorem]{Remark}
\newtheorem{question}{Question}

\numberwithin{equation}{section}

\newcommand{\R}{\mathbf{R}}
\newcommand{\Z}{\mathbf{Z}}
\newcommand{\N}{\mathbf{N}}
\newcommand{\bS}{\mathbf{S}}
\newcommand{\cA}{\mathcal{A}}
\newcommand{\cE}{\mathcal{E}}
\newcommand{\eps}{\varepsilon}
\newcommand{\loc}{{\rm loc}}
\newcommand{\mi}{{\rm min}}
\newcommand{\ma}{{\rm max}}
\newcommand{\Carl}{{\rm Carl}}
\newcommand{\cov}{{\rm cov}}

\DeclareMathOperator{\dist}{\operatorname{dist}}

\begin{document}

\title[Critical norm blow-up rates]{Critical norm blow-up rates 
for the energy supercritical nonlinear heat equation}

\author[T. Barker]{Tobias Barker}
\address{Department of Mathematical Sciences, University of Bath, Bath BA2 7AY. UK}
\email[Corresponding author]{tb2130@bath.ac.uk}

\author[H. Miura]{Hideyuki Miura}
\address{Department of Mathematics, 
Institute of Science Tokyo, Tokyo 152-8551, Japan}
\email{hideyuki@math.titech.ac.jp}

\author[J. Takahashi]{Jin Takahashi}
\address{Department of Mathematical and Computing Science, 
Institute of Science Tokyo, Tokyo 152-8552, Japan}
\email{takahashi@c.titech.ac.jp}



\begin{abstract}
We prove the first classification of blow-up rates of the critical norm 
for solutions of the energy supercritical nonlinear heat equation, 
without any assumptions such as radial symmetry or sign conditions. 
Moreover, the blow-up rates we obtain are optimal, 
for solutions that blow-up with 
bounded $L^{n(p-1)/2,\infty}(\R^n)$-norm up to the blow-up time.

We establish these results by proving quantitative estimates 
for the energy supercritical nonlinear heat equation 
with a robust new strategy based on 
quantitative $\varepsilon$-regularity criterion 
averaged over certain comparable time scales. 
With this in hand, we then produce the quantitative estimates 
using arguments inspired by Palasek \cite{Pa22} and Tao \cite{Ta21} 
involving quantitative Carleman inequalities applied to the Navier-Stokes equations.

Our work shows that energy structure is not essential 
for establishing blow-up rates of the critical norm 
for parabolic problems with a scaling symmetry. 
This paves the way for establishing such critical norm blow-up rates 
for other nonlinear parabolic equations. 
\end{abstract}

\maketitle

\tableofcontents

\section{Introduction}

\subsection{Background}
We consider the following nonlinear heat equation: 
\begin{equation}\label{eq:fujitaeq1}
\left\{ 
\begin{aligned}
	&u_t=\Delta u+|u|^{p-1}u
	&&\mbox{ in }\R^n\times(0,T), \\
	&u(\cdot,0)=u_0
	&&\mbox{ on }\R^n, 
\end{aligned}
\right. 
\end{equation}
where $n\geq1$ and $p>1$. 
The equation enjoys the invariance under the scaling 
$u(x,t) \mapsto \lambda^{2/(p-1)} u(\lambda x,\lambda^2 t)$ for $\lambda>0$. 
The invariance defines the critical Lebesgue space $L^{q_c}(\R^n)$, where 
\[
	q_c:= \frac{n(p-1)}{2}. 
\]
By local-existence theory, 
\eqref{eq:fujitaeq1} has a unique local-in-time classical solution 
(see Subsection \ref{subsec:nos} for a definition)
if $u_0\in L^{q_c}(\R^n)$ with $q_c>1$. 
The solution can be continued as a classical solution 
up to the maximal time of existence $T\leq \infty$. 
If $T<\infty$, the blow-up occurs in the sense of $L^\infty$. 
In the pioneering work of 
Giga and Kohn \cite{GK87}, they studied the blow-up rate of the $L^\infty$-norm, 
and then they proposed the following question in \cite[Section 7]{GK87}.

\begin{question}\label{que:GK}
But it is natural to ask about other norms as well, 
for example $\|u\|_{L^q}$: do they blow up as $t\to T$, 
and if so, at what rate?
\end{question}

For $q>q_{c}$, local-existence theory provides an answer. 
Indeed,  
\[
	\|u(\cdot,t)\|_{L^q(\R^n)}
	\geq C(n,p,q) (T-t)^{-\frac{n}{2}(\frac{1}{q_c}-\frac{1}{q})} 
\]
holds for finite time blow-up solutions $u$, 
see \cite[Section 6]{We81} and \cite[Remark 16.2 (iii)]{QSbook2}. 
However, in the critical case $q=q_c$, 
even the first part of Question \ref{que:GK} 
does not follow from local-existence theory and 
has been asked in Brezis-Cazenave \cite[Open problem 7]{BC96} 
and Quittner-Souplet \cite[OP2.1, Section 55]{QSbook2}.

For $q=q_c$, the subtlety of Question \ref{que:GK} 
is shown by the fact \cite{dPMW19,dPMWZZpre,Sc12} that 
there exist type II blow-up solutions satisfying 
$\sup_{0<t<T} \|u(\cdot,t)\|_{L^{q_c}(\R^n)}<\infty$ 
for the energy critical nonlinearity $p=p_S$ with $3\leq n\leq 5$, where 
\[
	p_S:= \frac{n+2}{n-2} \quad \mbox{ for } n\geq 3. 
\]
Nevertheless, 
under the type I assumption, 
Mizoguchi and Souplet \cite{MS19} showed that 
$\lim_{t\to T} \|u(\cdot,t)\|_{L^{q_c}(\R^n)}=\infty$ 
holds for each $p>1$. 
Here we note that 
the blow-up is of type I if 
$\limsup_{t\to T} (T-t)^{1/(p-1)}\|u(\cdot,t)\|_{L^\infty(\R^n)}<\infty$ 
and type II if it is not of type I. 
In addition, the proof in \cite{MS19} is inspired 
by Escauriaza-Seregin-\v{S}ver\'{a}k's seminal result 
showing the regularity of solutions $v$ of the 
$3$-dimensional Navier-Stokes equations under the assumption 
that $\sup_{0<t<T} \|v(\cdot,t)\|_{L^3(\R^n)}<\infty$.

For the energy supercritical nonlinearity $p>p_S$, 
the second and third author \cite{MTap,MTpre} 
removed the type I assumption to show 
$\limsup_{t\to T} \|u(\cdot,t)\|_{L^{q_c}(\R^n)}=\infty$, 
and subsequently 
\[
	\lim_{t\to T} \|u(\cdot,t)\|_{L^{q_c}(\R^n)}=\infty. 
\]
Thus, for $p>p_S$, the first part of Question \ref{que:GK} 
is completely resolved. 
In this paper, we are interested in the second part, 
which we label Question \ref{que:GKq}.

\begin{question}\label{que:GKq}
For $p>p_S$, at what rate does $\|u(\cdot,t)\|_{L^{q_c}(\R^n)}$ 
blow up as $t\to T$? 
\end{question}

Such an open problem has also been stated in \cite[Remark 1.6]{MTpre}. 
The main challenge for resolving Question \ref{que:GKq} is that 
the qualitative results \cite{MTap,MTpre} (see also \cite{MTno}) crucially use 
the Giga-Kohn monotonicity formula to classify 
the blow-up limit as backward self-similar solutions. 
Such qualitative features are challenging to use in a quantitative way. 
We are able to overcome this obstacle to obtain our main results below 
by using quantitative Carleman inequalities established 
by Palasek \cite{Pa22} and Tao \cite{Ta21}, which in turn use Carleman inequalities 
established by Escauriaza-Seregin-\v{S}ver\'{a}k \cite{ESS03a}. 
As a by-product, our method also can be used to estimate 
the blow-up rate of infinite-time blow-up solutions. 
In the next subsection, we list our main results.

\subsection{Statement of results}

Our first main result addresses Question \ref{que:GKq} 
by means of a quadruple logarithmic rate in the general setting.

\begin{theorem}\label{thm:critrateheat}
Let $n\geq 3$,  $p>p_S$ and 
$u$ be a classical solution 
of \eqref{eq:fujitaeq1} with $u_0\in L^{q_c}(\R^n)$. 
If the maximal time of existence $T>0$ is finite, then 
\[
	\limsup_{t\to T} \frac{ \|u(\cdot,t)\|_{L^{q_{c}}(\R^n)} }
	{\left(\log\log\log\log\left(
	\frac{1}{(T-t)^{\frac{1}{2(p-1)}}}
	\right)\right)^c}=\infty, 
\]
where $c>0$ is a constant depending only on $n$ and $p$. 
\end{theorem}

While Theorem \ref{thm:critrateheat} provides a rate 
for Question \ref{que:GKq}, we do not know if the rate is optimal. 
To the best of our knowledge, for constructed blow-up solutions 
to the nonlinear heat equation, the $L^{q_c}$-norm typically diverges 
at a logarithmic rate even when the blow-up is type II, 
see for example \cite[Corollary 1.5]{Se18} for 
$p=p_{JL}:=(n-2\sqrt{n-1})/(n-4-2\sqrt{n-1}) (>p_S)$ with $n\geq11$, 
\cite[Corollary 3]{MS21} for $p>p_{JL}$ and 
\cite[Subsection 8.2.1]{Ha20s} for $p=p_S$ with $n=6$. 
In particular, the solutions constructed by Seki \cite{Se18} and 
Mukai-Seki \cite{MS21} satisfy
\begin{align}
	&\label{eq:Lobdd} 
	\sup_{0<t<T} \|u(\cdot,t)\|_{L^{q_c,\infty}(\R^n)} <\infty, \\
	&\label{eq:lograte} 
	C' \log\left( \frac{1}{T-t} \right) \leq \int_{\R^n} |u(x,t)|^{q_{c}} dx
	\leq C'' \log\left( \frac{1}{T-t} \right), 
\end{align}
where we also refer to \cite[Subsection 4.1]{MM09} and \cite[Proposition C.1]{MM11} 
for deriving \eqref{eq:Lobdd}.

Our second main result shows that for solutions satisfying \eqref{eq:Lobdd}, 
the lower bound of the blow-up rate in \eqref{eq:lograte} is generic.

\begin{theorem}\label{thm:critratetypeIheat}
Let $n\geq 3$, $p>p_S$ and $u$ be a classical solution 
of \eqref{eq:fujitaeq1} with $u_0\in L^{q_c}(\R^n)$. 
There exist constants $M_0, C>0$ depending only on $n$ and $p$ such that 
if maximal time of existence $T>0$ is finite and 
$u$ satisfies 
\[
	\sup_{0<t<T}\|u(\cdot,t)\|_{L^{q_{c},\infty}(\R^n)}\leq M
\]
with a constant $M\geq M_0$, then 
\[
	\int_{\R^n} |u(x,t)|^{q_{c}} dx 
	\geq  e^{-e^{e^{M^C}}} \log\left( \frac{t}{T-t} \right) 
\]
for all $( 1+ M^{-2(2p^2+3p+2)})^{-1} T<t<T$. 
\end{theorem}

Our third result gives the first classification of infinite-time blow-up solutions 
in terms of the critical norm, which answers a question raised in 
\cite[Remark 1.9]{MTap}. 
In what follows, we write $B(x,r):= \{y\in\R^n; |x-y|<r\}$ 
for $x\in\R^n$ and $r>0$.

\begin{theorem}\label{thm:critratetypeIheatinfinite}
Let $n\geq 3$, $p>p_S$ and $u$ 
be a global-in-time classical solution of \eqref{eq:fujitaeq1}. 
There exist  constants $C,M_0>0$ depending only on $n$ and $p$ such that 
if $\limsup_{t\to \infty}|u(0,t)|=\infty$ and 
\[
	\sup_{0<t<\infty}\|u(\cdot,t)\|_{L^{q_{c},\infty}(\R^n)}\leq M
\]
with a constant $M\geq M_0$, then we conclude that there exists a sequence $t_{n}\to \infty$ 
such that
\begin{equation}\label{eq:BPrateheatinfinite}
	\int_{B(0, t_{n}^{\frac{1}{2}}e^{e^{M^C}}) } |u(x,t_n)|^{q_{c}} dx 
	\geq e^{-e^{e^{M^C}}}\log t_n.  
\end{equation}
\end{theorem}

Theorem \ref{thm:critratetypeIheatinfinite} can be applied 
to the infinite-time blow-up solutions 
constructed by Pol\'{a}\v{c}ik-Yanagida's works 
\cite[Theorem 1.3]{PY03} for $p\geq p_{JL}$ 
and \cite[Theorem 1]{PY14} for $p_S<p< p_{JL}$. 
We note that the blow-up rate of $L^\infty$-norm is studied in 
\cite{FKWY06,FWY04,Mi06} for $p\geq p_{JL}$ 
and that Theorem~\ref{thm:critratetypeIheatinfinite} 
is also applicable to the solutions in these papers.

\begin{remark}
In order for the assumptions in Theorem \ref{thm:critratetypeIheatinfinite} 
to apply to the known solutions, 
it seems necessary that 
the integral in \eqref{eq:BPrateheatinfinite} is over a finite ball. 
In particular, we expect that there is no infinite-time blow-up solutions 
with $u_0\in L^{q_c}(\R^n)$. 
As evidence for this, under $p>p_S$ and $\nabla u_0 \in L^{q_*}(\R^n)$ 
with $q_*:=n(p-1)/(p+1)$ (this implies $u_0 \in L^{q_c}(\R^n)$ 
by the Sobolev inequality), we can prove that 
any associated global-in-time classical solution $u$ satisfies 
\[
	\sup_{1<t<\infty} t^\frac{1}{p-1}\|u(\cdot,t)\|_{L^\infty(\R^n)}<\infty, 
\]
by applying the argument of Souplet \cite[Proof of Theorem 2, (7.1)]{So17} 
with the aid of \cite[Proposition 5.1(ii)]{So17}, the convolution inequalities 
and the approximation 
$v_0\in C^\infty_0(\R^n)$ 
satisfying $\|\nabla (u_0-v_0)\|_{L^{q_*}(\R^n)}<\delta$ 
for a small constant $\delta>0$.
\end{remark}

The above theorems hinge on quantitative estimates 
for solutions of 
\begin{equation}\label{eq:fujitaeq}
	u_t=\Delta u+|u|^{p-1}u
	\quad 
	\mbox{ in }\R^n\times(-1,0], 
\end{equation}
which are classical up to (and including) $t=0$. 
Our main quantitative estimate is as follows. 
It will be applied to the rescaled function 
\[
	u_{t}(x,s):=(t^\frac{1}{2})^\frac{2}{p-1}u(t^\frac{1}{2} x,ts+t)
\]
with $t\in (0,T)$ fixed.

\begin{proposition}[main quantitative estimate]\label{pro:quantregheat}
Let $n\geq3$, $p>p_S$ and $u$ be a classical solution 
of \eqref{eq:fujitaeq}. 
Then there exist constants $C_0, C_1, M_0>0$ 
depending only on $n$ and $p$ such that 
the following statement holds true. 
Suppose that $u$ satisfies 
\begin{equation}\label{eq:Type1bound2.0}
	\sup_{-1<t<0}\|u(\cdot,t)\|_{L^{q_{c},\infty}(\R^n)}\leq M
\end{equation}
with a constant $M\geq M_{0}$. Let $x_0\in \R^n$ and set 
\begin{align}\label{eq:Ndef}
	&N:= 
	\int_{ B(x_0, e^{e^{M^{C_0}}}) }
	|u(x,0)|^{q_c} dx, \\
	&\label{tstardef}
	t_{*}:= - M^{-p(2p+3)-1} \exp\left(-Ne^{e^{e^{M^{C_{0}}}}}\right). 
\end{align}
Then, 
\begin{equation}\label{eq:quantest}
	\|u\|_{L^{\infty}(B(x_0, \frac{1}{4}(-t_*)^{\frac{1}{2}})
	\times (\frac{1}{16}t_{*},0))}
	\leq C_1 (-t_*)^{-\frac{1}{p-1}}. 
\end{equation}
\end{proposition}

In the next subsection, 
we compare our method with the previous literature.

\subsection{Comparison with previous literature}
In the cases $q_c=1,2$, 
from \cite{We86} and \cite[Propositions 16.3, 16.3a]{QSbook2}, 
it is possible to obtain a single logarithmic blow-up rate 
of the critical norm by elementary direct arguments. 
Otherwise, there are very few results in the literature 
regarding quantitative blow-up rates of critical norms 
for evolution equations with a scaling symmetry. 
To the best of our knowledge, the first 
quantitative blow-up rate for a critical norm 
for more general nonlinearities was obtained by 
Merle and Raph\"ael \cite{MR08}
for radial solutions of the $L^2$ supercritical nonlinear 
Schr\"odinger equation. 

In the parabolic setting, a breakthrough work of Tao \cite{Ta21} 
established that if a solution $v$ of the $3$-dimensional 
Navier-Stokes equations in $\R^n\times(-1,0)$ 
first loses smoothness at $t=0$, 
then the critical $L^3$-norm becomes unbounded 
with a quantitative estimate 
\begin{equation}\label{eq:taorate1}
	\limsup_{t\to 0} 
	\frac{ \| v(\cdot,t)\|_{L^3(\R^3)} }{
	\left(\log \log \log \left( \frac{1}{(-t)^{c_0}} \right) \right)^{c_1} }
	=\infty. 
\end{equation}
To show this, Tao's aim  
is the following: \\

\noindent \textit{Tao's objective.}
Assume 
\begin{equation}\label{eq:NSEcritbounded}
    A:=\sup_{-1<t<0} \|v(\cdot,t)\|_{L^3(\R^n)}<\infty.
    \end{equation}
If the following statement 
\begin{equation}\label{eq:TaoobC}
	N^{-1} \sup_{-1/2<t<0} \| P_N v\|_{L^3(\R^n)} <A^{-c} 
	\quad \mbox{ for all } N\geq N_* 
\end{equation}
fails, find a quantitative upper bound $N_*(A)$ for $N$, 
where $P_N$ is a Littlewood-Paley projection on the frequency $N$.\\

Once \eqref{eq:TaoobC} is established for all $N\geq N_*(A)$, 
this produces quantitative estimates of 
$\|v\|_{L^\infty(\R^3\times (-1/4,0))}$ in terms of $A$, 
which can then be used to prove \eqref{eq:taorate1}. 
See the survey \cite{BP24} for more details. 

In \cite{Ta21}, Tao achieves the above objective by 
\begin{enumerate}
\item
Backward propagation mechanism 
based on the contraposition of \eqref{eq:TaoobC}.
\item
Establishing regions of quantitative regularity.
\item
Applying quantitative Carleman inequalities 
based on \cite{ESS03a} to the quantity 
$\nabla \times v$. 
\item
Summing of scales of the $L^3$-norm of $v$ at the final time over disjoint annuli 
to produce the upper bound $N_*(A)$. 
\end{enumerate}
For obtaining our main quantitative estimate 
(Proposition \ref{pro:quantregheat}), 
which in turn implies Theorems \ref{thm:critrateheat}, 
\ref{thm:critratetypeIheat} and 
\ref{thm:critratetypeIheatinfinite}, 
we encounter several challenges in the above (1), (2) and (3) 
in the context of the energy supercritical heat equation. 
Let us describe these now in more details. 
The novelty of the results will also be described.

First, we discuss the difficulties regarding the quantitative backward propagation. 
For $3$-dimensional Navier-Stokes equations, Tao's objective 
(the contraposition of \eqref{eq:TaoobC}) produces a sequence of 
the following `frequency bubbles of concentration': 
For all $n\in \N$, there exists a frequency 
$N_n\in (0,\infty)$ and $(x_n,t_n)\in \R^3\times(-1,t_{n-1})
\subset \R^3\times (-1,0)$ such that 
\[
	N_n^{-1} |P_{N_n}u(x_n,t_n)| > A^{-c}
\]
with 
\[
	x_n=x_0 + A^C (-t_n)^\frac{1}{2}, \quad 
	N_n \sim A^C (-t_n)^{-\frac{1}{2}}, 
\]
see \cite{BP24} for more details. 
The frequency bubbles of concentration utilizes the Duhamel formula, 
localized estimates of Fourier multipliers and paraproduct decompositions. 
It seems non-trivial to obtain such 
a frequency based backward propagation mechanism 
for non-smooth nonlinearities, such as those 
we are faced with for the energy supercritical nonlinear heat equations 
in the general case.

The first author and Prange \cite{BP21} extended Tao's results 
and proved an analogue of Theorem \ref{thm:critratetypeIheat} 
for the $3$-dimensional Navier-Stokes equations 
by using a different objective. 
For brevity, we state this objective in a less general setting than in \cite{BP21} 
(same setting as Tao \cite{Ta21}). \\

\textit{Barker and Prange's objective.}
Assume $A:=\sup_{-1<t<0} \|v(\cdot,t)\|_{L^3(\R^n)}<\infty$. 
If the following statement 
\begin{equation}\label{eq:BPob}
	(-t)^\frac{1}{2} \int_{B(x_0, A^c(-t)^\frac{1}{2})} 
	|\nabla \times v(x,t)|^2 dx \leq A^{-c}
\end{equation}
fails, find a quantitative upper bound $t_*(A)\in (-1,0)$ for $t$ 
such that necessarily $t\leq t_*(A)$.\\

Once \eqref{eq:BPob} is established for all $t_{*}(A)<t<0$, 
this produces quantitative estimates of 
$\|v\|_{L^\infty(\R^3\times (-1/4,0))}$ in terms of $A$. 
A key part in achieving Barker and Prange's objective in \cite{BP21} 
is showing `backward propagation of vorticity concentration'. 
Namely, the initial concentration 
\begin{equation}\label{eq:initialvortconc}
	(-t)^\frac{1}{2} \int_{B(x_0, (-t)^\frac{1}{2})} 
	|\nabla\times v(x,t)|^2 dx > A^{-c} 
\end{equation}
propagates backwards in time and holds true for all times that are sufficiently 
in the past of $t$. 
Both Barker and Prange's objective and 
backward propagation of vorticity concentration crucially hinge on 
`local-in-space smoothing' for the $3$-dimensional Navier-Stokes equations. 
Namely, if the initial data $v_0$ satisfies 
$\|v_0\|_{L^3(\R^n)}\leq M$ and 
$\|v_0\|_{L^6(B(0,1))}\leq N$, then the solution $v$ is quantitatively 
bounded on $B(0,1/2)\times (0,T(M,N)]$. 
The proof of local-in-space smoothing crucially uses the local energy 
structure of the $3$-dimensional Navier-Stokes equations. 
In particular, it is not obviously clear that it holds 
for other nonlinear parabolic equations such as 
the energy supercritical nonlinear heat equation. 

To overcome these difficulties, 
we pursue a different strategy for producing quantitative 
estimates for the energy supercritical nonlinear heat equation 
(Proposition \ref{pro:quantregheat}). 
Our strategy is based on quantitative partial regularity 
over comparable time scales 
under the assumption \eqref{eq:Type1bound2.0}.\\

\textit{New objective.}
Assume that $u$ is a solution of the energy supercritical nonlinear heat equation 
on $\R^n\times(-1,0)$ satisfying the Lorentz norm bound \eqref{eq:Type1bound2.0}. 
Let $N$ be defined by \eqref{eq:Ndef}.
If the following statement 
\begin{equation}\label{eq:newob}
	(-t)^{\frac{2}{p-1} - \frac{n}{2}} 
	\int_t^{t/2} 
	\int_{B(x_0, A(M,p) (-t)^\frac{1}{2})} 
	|u(x,t)|^{p+1} dx dt \leq M^{-c(p)}
\end{equation}
fails, find a quantitative upper 
bound $t_*(M,N)\in (-1,0)$ for $t$ such that necessarily $t\leq t_*(M,N)$.\\

Once \eqref{eq:newob} is established for $t_*(M,N)<t<0$, 
this will then imply the quantitative bounds in Proposition \ref{pro:quantregheat}.
To pursue this objective, we must first quantify 
the $\eps$-regularity criterion (Proposition \ref{pro:epsreg}) 
concerning \eqref{eq:newob} with \eqref{eq:Type1bound2.0}, 
which is accomplished using the Giga-Kohn monotonicity formula 
\cite[Proposition 2.1]{GK87} and 
Blatt-Struwe's $\eps$-regularity criterion \cite[Proposition 4.1]{BS15}.
With this in hand, we can show that the initial concentration 
\[
	(-t)^{\frac{2}{p-1} - \frac{n}{2}} 
	\int_t^{t/2} 
	\int_{B(x_0, A(M,p) (-t)^\frac{1}{2})} 
	|u(x,t)|^{p+1} dx dt > M^{-c(p)}
\]
propagates backwards in time and holds for times that 
are (quantifiably) sufficiently in the past of $t$. 
This forms a crucial part of our proof of Proposition \ref{pro:quantregheat}.

Next, we discuss the obstacle on the quantitative regions of regularity. 
In the works \cite{Ba23,BP21, Hu24, Pa22,Ta21} on the Navier-Stokes equations, 
a key part in the production of the quantitative estimates 
is the use of quantitative Carleman inequalities in Tao \cite[Section 4]{Ta21}. 
Applying such Carleman inequalities requires the vorticity 
to satisfy certain differential inequalities, which requires 
determining regions of quantitative regularity of the velocity. 
In \cite{Ta21}, to apply quantitative unique continuation Carleman inequality, 
it is used that a solution to $3$-dimensional Navier-Stokes equations 
satisfying 
\[
	\sup_{t\in I} \|v(\cdot,t)\|_{L^3(\R^n)} \leq M
\]
for some interval $I$, there exists a quantifiable 
sub-interval $I'\subset I$ (epoch of regularity) 
such that $v$ is quantitatively bounded on $\R^3 \times I'$. 
This crucially uses the energy structure of the Navier-Stokes equations 
and the Sobolev embedding theorem in $n=3$, 
which imply that solutions in $\R^3$ with bounded energy are 
in subcritical spaces on many time slices. 
Here we refer to a Lebesgue space as `subcritical' for 
a parabolic evolution equation with scaling symmetry, 
if it has higher integrability than the (scale-invariant) critical Lebesgue spaces. 
For the higher-dimensional Navier-Stokes equations ($n\geq4$) with 
\[
	\sup_{-1 <t<0} \|v(\cdot,t)\|_{L^n(\R^n)} \leq M, 
\]
it is not known if solutions possess quantitative epoch of regularity. 
This is also the case for solutions of the energy supercritical 
nonlinear heat equation with \eqref{eq:Type1bound2.0}.

To overcome this obstacle, we utilize ideas from Palasek \cite{Pa22} 
for the higher-dimensional Navier-Stokes equations. 
In particular, we use the quantitative space-time partial regularity 
\eqref{eq:newob} with \eqref{eq:Type1bound2.0} 
to find quantitative space-time `slices of regularity' for the solution. 
In doing this, we face an additional obstacle compared to \cite{Pa22} 
in that the norm $L^{q_c,\infty}(\R^n)$ in \eqref{eq:Type1bound2.0} 
can have an equal presence at many disjoint scales. 
We overcome this by performing a Calder\'on splitting 
of the solution \cite{Ca90}, 
together with the fact that the integrability exponent 
$p+1$ in the quantitative partial regularity 
\eqref{eq:newob} with \eqref{eq:Type1bound2.0} 
is lower than the critical exponent $q_c$ 
due to the supercriticality $p>p_S$.

Once we have established quantitative backward propagation 
and regions of regularity, we then utilize quantitative Carleman 
inequalities in \cite{Pa22} to obtain 
our main quantitative estimate (Proposition \ref{pro:quantregheat}). 
Due to cases involving a non-smooth nonlinearity, 
we are required to implement such Carleman inequalities 
in a lower regularity setting compared 
to previous works on the Navier-Stokes equations.

\begin{remark}
The quadruple logarithmic blow-up rate 
in Theorem \ref{thm:critrateheat} is due to the quadruple exponential 
quantitative estimate in Proposition \ref{pro:quantregheat}. 
Tao's triple logarithmic blow-up rate \eqref{eq:taorate1} 
for the $3$-dimensional Navier-Stokes equations 
is due to a triple exponential quantitative estimate, 
with the reason for the triple exponential estimate 
outlined in \cite[Remark 1.5]{Ta21}. 
By comparison, in Palasek's work \cite{Pa22} on the higher-dimensional 
Navier-Stokes equations and this paper, 
the quantitative estimates involve a further exponential loss 
for the following reason. 
For the $3$-dimensional Navier-Stokes equations 
satisfying \eqref{eq:NSEcritbounded} and \eqref{eq:initialvortconc}, 
quantitative unique continuation implies that 
for certain temporal scales $T_{1}$ and sufficiently large $R$, 
\[
	 \int_{-T_1}^{-T_1/2} \int_{R/2\leq |x|\leq 2R} 
	 |\nabla\times v (x,t)|^2 dxdt\geq T_{1}^{\frac{1}{2}}e^{-\frac{A^c R^2}{T_1}}
\]
(see \cite[(99)]{BP24}, for example). 
Yet in \cite{Pa22} and our setting, 
the use of iterated quantitative unique continuation 
produces a much smaller lower bound  
than the above Gaussian lower bound for the $L^2$-norm of the solution. 
In both settings, this necessitates 
the subsequent use of far larger spatial scales 
for applying quantitative backward uniqueness 
than those used for the $3$-dimensional Navier-Stokes equations. 
This requirement, together with pigeonhole arguments, 
produces an extra (quadruple) exponential 
in the quantitative estimates compared to the $3$-dimensional Navier-Stokes equations.
\end{remark}

\subsection{Final remarks}
Our work shows that energy structure is not an essential feature 
when proving quantitative blow-up rates of critical norms 
of nonlinear parabolic equations with a scaling symmetry. 
To prove our results, we only required the quantitative 
partial regularity at comparable time scales 
\eqref{eq:newob} for solutions satisfying \eqref{eq:Type1bound2.0}, 
along with the fact that regular solutions satisfy the correct differential inequality 
to apply quantitative Carleman inequalities. 
This opens the door for obtaining quantitative blow-up rates of critical norms 
for other nonlinear parabolic equations. 

\subsection{Organization of the paper}
The rest of this paper is organized as follows. 
In Section \ref{sec:epsreg}, we prove the quantitative 
$\eps$-regularity theorem. 
As applications, in Section \ref{sec:estiveps}, 
we prepare regularity estimates concerning propagation of concentration 
and regions of regularity. 
In Section \ref{sec:Carl}, 
we give Carleman inequalities regarding 
quantitative backward uniqueness and unique continuation. 
By combining these ingredients, we prove 
our main quantitative estimate (Proposition \ref{pro:quantregheat}) 
in Section \ref{sec:prmain}. 
In Section \ref{sec:proof}, we show Theorems \ref{thm:critrateheat}, \ref{thm:critratetypeIheat} and \ref{thm:critratetypeIheatinfinite}.

\subsection{Notation}
For $(x,t)\in\R^n \times \R$ and $r>0$, 
we write $B(x,r):= \{y\in\R^n; |x-y|<r\}$ and $B(r) := B(0,r)$. 
We denote by $Q((x,t), r):= B(x,r)\times (t-r^2,t)$ 
and $Q(r):=Q((0,0), r)$ the backward parabolic cylinders. 
Throughout this paper, 
$C$ denotes positive constants 
depending only on $n$ and $p$ unless otherwise stated. 
To stress the dependence, we also write $C(a,b,\ldots)$ 
when it depends only on $a$, $b$, \ldots.
Each of the constants may change from line to line during the proofs. 
We set $q_c:=n(p-1)/2$ and $q_*:=n(p-1)/(p+1)$, where 
$L^{q_*}$ is the scaling invariant critical space for the gradient of solutions. We denote $C^{2,1}$ to be the space of functions 
which are twice continuously differentiable in the space variable  
and once in the time variable.

\subsection{Notions of solutions}\label{subsec:nos}
We say that `$u$ is a classical solution of \eqref{eq:fujitaeq1} 
with $u_0\in L^q(\R^n)$' ($q\geq1$) 
if $u$ belongs to 
$C([0,T); L^q(\R^n))\cap L^\infty_\loc((0,T);L^\infty(\R^n)) 
\cap C^{2,1}(\R^n\times (0,T))$ 
and satisfies \eqref{eq:fujitaeq1}. 
Such a solution exists uniquely for any $u_0\in L^{q_c}(\R^n)$ 
with $q_c>1$, see \cite{BC96,We80} and \cite[Remark 15.4 (i)]{QSbook2}. 
We note that Theorems \ref{thm:critrateheat} and \ref{thm:critratetypeIheat} 
can be applied to classical solutions of \eqref{eq:fujitaeq1} 
with $u_0\in L^{q_c}(\R^n)$.

By a classical solution of \eqref{eq:fujitaeq1} 
(resp. \eqref{eq:fujitaeq})
without mentioning the initial data, 
we mean a function in $C^{2,1}(\R^n\times (0,T))$ 
(resp. $C^{2,1}(\R^n\times (-1,0])$)
without specifying the initial data. 
In particular, we do not impose 
$C([0,T); L^q(\R^n))$ ($q\geq1$) or 
$L^\infty_\loc((0,T);L^\infty(\R^n))$. 
We can apply Theorem \ref{thm:critratetypeIheatinfinite} 
and Proposition \ref{pro:quantregheat} 
to classical solutions without mentioning the initial data.

\section{Quantitative $\eps$-regularity}\label{sec:epsreg}
Let $u$ be a classical solution of \eqref{eq:fujitaeq} 
satisfying the Lorentz norm bound \eqref{eq:Type1bound2.0}. 
The goal of this section is to show the following 
quantitative $\eps$-regularity result.

\begin{proposition}[Quantitative $\eps$-regularity]\label{pro:epsreg}
Assume $p>p_S$ and \eqref{eq:Type1bound2.0} with a constant $M>1$. 
Then there exist constants $C>0$ and $0<\eps_1<1$ 
depending only on $n$ and $p$ such that  
the following holds for 
any $0<\eps<\eps_0:= M^{-2(p+1)^2} \eps_1^{2(p+1)}$: 
If 
\begin{equation}\label{eq:1scale}
\left\{ 
\begin{aligned}
	&\delta^{\frac{4}{p-1}-n} 
	\int_{t_0-\delta^2}^{t_0-\delta^2/2} \int_{B(x_0, A\delta)} 
	 |u(x,t)|^{p+1} dxdt \leq \eps \quad  \mbox{ for some }x_0\in\R^n, \\
	&-1/16<t_0\leq0, 
	0<\delta< 1/4 \mbox{ and } 
	A>(24\log (M^{p+1}/\eps) )^{1/2}, 
\end{aligned} \right. 
\end{equation}
then 
\begin{equation*}
\begin{aligned}
	&\|u\|_{L^\infty( Q((x_0,t_0), \delta/4))} \leq 
	C M \eps^\frac{1}{2(p+1)^2} \delta^{-\frac{2}{p-1}}, \\
	&\|\nabla u\|_{L^\infty( Q((x_0,t_0), \delta/4))} \leq 
	C M \eps^\frac{1}{2(p+1)^2} \delta^{-\frac{p+1}{p-1}}. 
\end{aligned}
\end{equation*}
\end{proposition}

We prove Proposition \ref{pro:epsreg} based on the analysis of 
the Giga--Kohn weighted energy. 
For $\tilde x\in \R^n$, $-1<t<\tilde t\leq 0$, 
we define the energy $E_{(\tilde x,\tilde t)}$ by 
\[
\begin{aligned}
	E_{(\tilde x,\tilde t)}(t) &:= (\tilde t-t)^\frac{p+1}{p-1}
	\int_{\R^n} \left( \frac{1}{2} |\nabla u(x,t)|^2
	- \frac{1}{p+1}|u(x,t)|^{p+1} \right. \\
	&\quad \left. + \frac{1}{2(p-1)(\tilde t-t)} |u(x,t)|^2
	\right)
	K_{(\tilde x,\tilde {t})}(x,t)  dx, 
\end{aligned}
\]
where $K$ is the backward heat kernel given by 
\[
	K_{(\tilde x,\tilde {t})}(x,t):= 
	(\tilde t-t)^{-\frac{n}{2}} 
	e^{-\frac{|x-\tilde x|^2}{4(\tilde t-t)}}. 
\]
We recall the backward similarity variables $(\eta,\tau)$ by 
\[
	\eta:=\frac{x-\tilde x}{(\tilde t- t)^{1/2}}, \quad \tau:= -\log(\tilde t -t). 
\]
Then the backward rescaled solution $w_{(\tilde x, \tilde t)}$ is defined by 
\[
	w_{(\tilde x, \tilde t)}(\eta,\tau)
	:= e^{-\frac{1}{p-1}\tau} 
	u(\tilde x+e^{-\frac{1}{2}\tau}\eta,  \tilde t - e^{-\tau}) 
	= (\tilde t-t)^\frac{1}{p-1}u(x,t). 
\]
The corresponding Giga--Kohn energy $\cE_{(\tilde x, \tilde t)}$ is given by 
\begin{equation}\label{eq:cEdef}
	\begin{aligned}
	\cE_{(\tilde x,\tilde t)}(\tau) 
	&:= \int_{\R^n} \left( 
	\frac{1}{2} |\nabla w_{(\tilde x, \tilde t)}(\eta,\tau)|^2  
	- \frac{1}{p+1} |w_{(\tilde x, \tilde t)}(\eta,\tau)|^{p+1}  \right. \\
	&\quad \left. + \frac{1}{2(p-1)} |w_{(\tilde x, \tilde t)}(\eta,\tau)|^2 \right) 
	\rho(\eta) d\eta, 
	\end{aligned}
\end{equation}
where $\rho(\eta):=e^{-|\eta|^2/4}$. 
We note that 
\begin{equation}\label{eq:EcEeqi}
	E_{(\tilde x, \tilde t)}(t)
	=\cE_{(\tilde x, \tilde t)}(\tau) \quad \mbox{ with }\tau=-\log(\tilde t-t). 
\end{equation}
Moreover, by setting $\tau_0:=-\log(\tilde t+1)$, 
we see that $w$ satisfies 
\begin{equation}\label{eq:weq}
	\rho w_\tau 
	= \nabla \cdot (\rho \nabla  w) - \frac{1}{p-1} w\rho +|w|^{p-1}w\rho
	\quad \mbox{ in } \R^n\times(\tau_0, \infty).  
\end{equation}
 Here and below, 
we often suppress the subscript $(\tilde x, \tilde t)$ 
if it is clear from the context.

We recall the 
Giga--Kohn monotonicity formula \cite{GK85,GK87}. 

\begin{lemma}\label{lem:Emono}
For any $\tilde x\in \R^n$, $-1<t'<t<\tilde t\leq 0$ 
and $\tau=-\log(\tilde t-t)$, 
\begin{align}
	&\label{eq:GKmn}
	\frac{d \cE_{(\tilde x,\tilde t)}}{d\tau}(\tau)
	= - \int_{\R^n} 
	|\partial_\tau w_{(\tilde x,\tilde t)}(\eta,\tau)|^2 
	\rho(\eta) d\eta,  \\
	&\label{eq:GKmnori}
	E_{(\tilde x,\tilde t)}(t) 
	= E_{(\tilde x,\tilde t)}(t') 
	- \int_{t'}^t \int_{\R^n} (\tilde t-s)^{\frac{2}{p-1}-1} 
	|S_{(\tilde x,\tilde t)}(x,s)|^2 
	K_{(\tilde x,\tilde t)}(x,s) dx ds, \\
	&\notag
	S_{(\tilde x,\tilde t)}(x,t) := 
	\frac{1}{p-1} u(x,t) + \frac{1}{2}(x-\tilde x) \cdot \nabla u(x,t) 
	- (\tilde t-t ) u_t(x,t). 
\end{align}
In particular, $\cE_{(\tilde x,\tilde t)}(\tau)$ 
and $E_{(\tilde x,\tilde t)}(t)$  
are nonincreasing in $\tau$ and $t$, respectively. 
\end{lemma}

\begin{proof}
See \cite[Proposition 2.1]{GK87}. 
\end{proof}

The monotonicity guarantees the nonnegativity of $E_{(\tilde x,\tilde t)}$. 
Moreover, the critical norm bound \eqref{eq:Type1bound2.0} 
implies the uniform boundedness of $E_{(\tilde x,\tilde t)}$.

\begin{lemma}\label{lem:Ebdd}
There exists $C>0$ depending only on $n$ and $p$ such that 
for any $\tilde x\in \R^n$ and  $-1/4\leq t<\tilde t\leq 0$, 
\[
	0\leq E_{(\tilde x,\tilde t)}(t) \leq C M^{2p}. 
\]
\end{lemma}

\begin{proof}
From \eqref{eq:weq} and the integration by parts, it follows that 
\[
	\frac{1}{2} \frac{d}{d\tau} \int_{\R^n} |w|^2 \rho d\eta
	= \int_{\R^n} \left( 
	- |\nabla  w|^2 - \frac{1}{p-1}|w|^2 +|w|^{p+1} 
	\right) \rho d\eta. 
\]
Then, \eqref{eq:cEdef} shows that 
\begin{equation}\label{eq:w22} 
	\frac{1}{2} \frac{d}{d\tau} 
	\int_{\R^n} |w|^2 
	\rho d\eta = -2\cE(\tau) 
	+ \frac{p-1}{p+1} \int_{\R^n} |w|^{p+1} \rho d\eta. 
\end{equation}
This together with 
Lemma \ref{lem:Emono} and the H\"older inequality yields 
\[
	\frac{1}{2} \frac{d}{d\tau} \int_{\R^n} |w(\eta,\tau)|^2 \rho(\eta) d\eta
	\geq  -2\cE(\tau') 
	+ C \left( \int_{\R^n} |w(\eta,\tau)|^2 \rho(\eta) d\eta \right)^\frac{p+1}{2} 
\]
for all $\tau_0<\tau'<\tau<\infty$ with $\tau_0:=-\log(\tilde t+1)$. 
This inequality and a contradiction argument show that 
$E_{(\tilde x,\tilde t)}(t)=\cE_{(\tilde x,\tilde t)}(\tau)\geq 0$ 
for any $\tilde x\in \R^n$ and $-1\leq t<\tilde t\leq 0$. 
For details, 
see \cite[Propositions 2.1, 2.2]{GK87} 
and \cite[Proposition 23.8]{QSbook2}.

We show the upper bound, which hinges on the strategy 
in \cite[Lemma 2.4]{MTpre}. 
From \eqref{eq:GKmnori} and \eqref{eq:Type1bound2.0}, 
it follows that 
\begin{equation}\label{eq:EM2nab}
\begin{aligned}
	&E_{(\tilde x,\tilde t)}(t) \leq E_{(\tilde x,\tilde t)}(-1/4) \\
	&\leq 
	C (\tilde t+1/4)^\frac{p+1}{p-1} 
	\int_{\R^n} \left( 
	|\nabla u(x,-1/4)|^2 + \frac{|u(x,-1/4)|^2}{\tilde t+1/4} \right)
	K_{(\tilde x,\tilde t)}(x,-1/4) dx \\
	&\leq 
	CM^2 + 
	C (\tilde t+1/4)^\frac{p+1}{p-1} 
	\int_{\R^n} |\nabla u(x,-1/4)|^2 
	K_{(\tilde x,\tilde t)}(x,-1/4) dx
\end{aligned}
\end{equation}
for $-1/4\leq t<\tilde t\leq 0$.  
By using the Duhamel formula for solutions of the nonlinear heat equation, 
we have an integral equation for $u(x,-1/4)$. 
Differentiating the integral equation and 
using $|\nabla G| \leq C K_1$ with 
the heat kernel $G(x,t):=(4\pi t)^{-n/2}e^{-|x|^2/(4t)}$ 
and $K_1(x,t):=t^{-(n+1)/2}e^{-|x|^2/(8t)}$,
we see that  
\[
\begin{aligned}
	|\nabla u(x,-1/4)| &\leq 
	C \int_{\R^n} K_1(x-y,1/4) |u(y,-1/2)| dy  \\
	&\quad 
	+ C \int_{-1/2}^{-1/4} \int_{\R^n} 
	K_1(x-y,-1/4-s) |u(y,s)|^p dyds \\
	&=: C U_1(x) + C U_2(x). 
\end{aligned}
\]

The H\"older inequality in the Lorentz spaces 
(see \cite[Section IV]{ON63}) gives 
\[
\begin{aligned}
	U_1(x) &\leq C 
	\left\| |u(\cdot,-1/2)| e^{-\frac{|x-\cdot|^2}{2}} 
	\right\|_{L^1(\R^n)} \\
	& \leq 
	C \| u(\cdot,-1/2) \|_{L^{q_c,\infty}(\R^n)}
	\left\| e^{-\frac{|x-\cdot|^2}{2}} 
	\right\|_{L^{\frac{q_c}{q_c-1},1}(\R^n)}
	\leq CM. \\
\end{aligned}
\]
Again by the H\"older inequality 
with $q_*:=n(p-1)/(p+1)>2$ for $p>p_S$, we have 
\[
\begin{aligned}
	&\int_{\R^n} |\nabla u(x,-1/4)|^2 K_{(\tilde x,\tilde t)}(x,-1/4) dx 
	\leq 
	C\int_{\R^n} (|U_1|^2+|U_2|^2) K_{(\tilde x,\tilde t)}(x,-1/4) dx  \\
	&\leq 
	CM^2
	+ C (\tilde t+1/4)^{-\frac{n}{2}} 
	\left\| |U_2|^2 e^{-\frac{|\cdot-\tilde x|^2}{4(\tilde t+1/4)}}
	\right\|_{L^1(\R^n)} \\
	&\leq 
	CM^2
	+ C (\tilde t+1/4)^{-\frac{n}{2}} \| |U_2|^2 \|_{L^{\frac{q_*}{2},\infty}(\R^n)}  
	\left\| e^{-\frac{|\cdot-\tilde x|^2}{4(\tilde t+1/4)}} 
	\right\|_{L^{\frac{q_*}{q_*-2},1} (\R^n)} \\
	&\leq 
	CM^2
	+ C (\tilde t+1/4)^{-\frac{p+1}{p-1}}  \| U_2 \|_{L^{q_*,\infty}(\R^n)}^2
	\leq 
	C(\tilde t+1/4)^{-\frac{p+1}{p-1}} 
	(M^2 +  \| U_2 \|_{L^{q_*,\infty}(\R^n)}^2). 
\end{aligned}
\]
Recall \eqref{eq:EM2nab}. Then, for $-1/4\leq t<\tilde t\leq 0$, we also have 
\begin{equation}\label{eq:EM2U2mid}
\begin{aligned}
	E_{(\tilde x,\tilde t)}(t) 
	&\leq 
	CM^2 + 
	C (\tilde t+1/4)^\frac{p+1}{p-1} 
	\left( 
	M^2
	+ (\tilde t+1/4)^{-\frac{p+1}{p-1}}  \| U_2 \|_{L^{q_*,\infty}(\R^n)}^2
	\right) \\
	&\leq 
	CM^2 + 
	C \| U_2 \|_{L^{q_*,\infty}(\R^n)}^2. 
\end{aligned}
\end{equation}

We estimate $U_2$ by a modification of \cite[Theorem 18.1]{Me97}. 
By the change of variables in the definition of $U_2$, 
we set 
\[
\begin{aligned}
	&U_2(x)= \int_0^\infty S(x,s) ds, \\
	&S (x,s):= 
	\chi_{(0,1/4)}(s) \int_{\R^n} K_1(x-y,s) |u(y,-1/4-s)|^p dy. 
\end{aligned}
\]
For $\lambda>0$ and $\tau>0$, 
define $D_\lambda:=\{x\in\R^n; U_2(x)>\lambda\}$ and 
\[
	U_2(x)=\left( \int_0^\tau + \int_\tau^\infty \right) S(x,s) ds 
	=: V_\tau(x) + W_\tau(x). 
\]
Let us estimate the Lebesgue measure $|D_\lambda|$. 
By the H\"older inequality, we have 
\[
\begin{aligned}
	S(x,s) 
	&\leq C s^{-\frac{n}{2}-\frac{1}{2}} \chi_{(0,1/4)}(s)
	\left\||u(\cdot,-1/4-s)|^p e^{-\frac{|x-\cdot|^2}{8s}}
	\right\|_{L^1(\R^n)} \\
	&\leq C s^{-\frac{n}{2}-\frac{1}{2}} \chi_{(0,1/4)}(s)
	\||u(\cdot,-1/4-s)|^p \|_{L^{\frac{q_c}{p},\infty}(\R^n)} 
	\left\| e^{-\frac{|x-\cdot|^2}{8s}} 
	\right\|_{L^{\frac{q_c}{q_c-p},1}(\R^n)}  
	\\
	&\leq C s^{-\frac{np}{2q_c}-\frac{1}{2}} 
	\chi_{(0,1/4)}(s)
	\left\| u(\cdot,-1/4-s)\right\|_{L^{q_c,\infty}(\R^n)}^p 
	\leq C M^p s^{-\frac{p}{p-1}-\frac{1}{2}}
\end{aligned}
\]
for any $s>0$. Then, 
\[
	W_\tau(x) 
	\leq C M^p \int_\tau^\infty s^{-\frac{p}{p-1}-\frac{1}{2}} ds 
	=C' M^p \tau^{-\frac{p+1}{2(p-1)}}, 
\]
where $C'>0$ is independent of $t$.  
For $\lambda>0$, we choose $\tau$ such that 
\begin{equation}\label{eq:tausig}
	C' M^p \tau^{-\frac{p+1}{2(p-1)}} = \frac{\lambda}{2}. 
\end{equation}
Then $W_\tau\leq \lambda/2$, and so 
$\tilde D_\lambda:=\{x\in\R^n; V_\tau(x) >\lambda/2\}$ 
satisfies $D_\lambda\subset \tilde D_\lambda$.

From the change of variables, it follows that 
\[
\begin{aligned}
	S(x,s) &= \int_{\R^n} K_1(z,s) \chi_{(0,1/4)}(s) |u(z+x,-1/4-s)|^p dz \\
	&\leq 
	C s^{-\frac{n}{2}-\frac{1}{2}} 
	\int_{\R^n}  \chi_{(0,1/4)}(s) |u(z+x,-1/4-s)|^p e^{-\frac{|z|^2}{8s}} dz 
\end{aligned}
\]
and so by O'Neil's convolution inequality \cite[Theorem 2.6]{ON63} 
\[
\begin{aligned}
	\|S(\cdot,s)\|_{L^{\frac{q_c}{p},\infty}(\R^n)} 
	&\leq 
	C s^{-\frac{n}{2}-\frac{1}{2}} 
	\int_{\R^n}  \chi_{(0,1/4)}(s) 
	\| u(\cdot,-1/4-s)|^p\|_{L^{\frac{q_c}{p},\infty}(\R^n)}
	e^{-\frac{|z|^2}{8s}} dz \\
	&\leq C M^p s^{-\frac{n}{2} -\frac{1}{2}} 
	\int_{\R^n} e^{-\frac{|z|^2}{8s}} dz
	\leq C M^p s^{-\frac{1}{2}}
\end{aligned}
\]
for any $s>0$. Thus, 
by Minkowski's inequality for Lorentz spaces \cite[Lemma 1]{mandel23}, 
\[
	\|V_\tau\|_{L^{\frac{q_c}{p},\infty}(\R^n)} 
	\leq \int_0^\tau \|S(\cdot,s)\|_{L^{\frac{q_c}{p},\infty}(\R^n)} ds 
	\leq C M^p \tau^\frac{1}{2}. 
\]
This together with the H\"older inequality 
for the Lorentz spaces shows that
\[
	\int_{\tilde D_\lambda} V_\tau(x) dx 
	\leq C \| \chi_{\tilde D_\lambda} \|_{L^{\frac{q_c}{q_c-p},1}(\R^n)} 
	\|V_\tau\|_{L^{\frac{q_c}{p},\infty}(\R^n)} 
	\leq C M^p |\tilde D_\lambda|^{1-\frac{p}{q_c}} \tau^\frac{1}{2}. 
\]
On the other hand, 
$\int_{\tilde D_\lambda} V_\tau(x) dx\geq (\lambda/2) |\tilde D_\lambda|$. 
By $D_\lambda\subset \tilde D_\lambda$ and \eqref{eq:tausig}, 
we obtain
\[
	|D_\lambda|\leq |\tilde D_\lambda|\leq C \lambda^{-\frac{q_c}{p}} 
	M^{q_c} \tau^\frac{q_c}{2p} 
	=C M^{p q_*} \lambda^{-q_*}, 
\]
and so $\lambda  |\{ x\in\R^n;  U_2(x)>\lambda\}|^{1/q_*} \leq C M^p$
for $\lambda>0$. 
Hence we obtain 
\[
	\|U_2\|_{L^{q_*,\infty}(\R^n)} 
	\leq 
	\|U\|_{L^{q_*,\infty}(\R^n)}
	\leq CM^p. 
\]
This together with \eqref{eq:EM2U2mid} and $M>1$ 
gives the desired inequality. 
\end{proof}

We prove the following lemma based on Blatt-Struwe's 
$\eps$-regularity criterion \cite[Proposition 4.1]{BS15}.

\begin{lemma}\label{lem:epspre}
There exist $0<\tilde \eps_1<1$ and $C>0$ 
depending only on $n$ and $p$ such that 
the following holds for any $0<\tilde \eps\leq \tilde \eps_1$: 
If there exists $0<\delta<1/2$ such that 
\begin{equation}\label{eq:Irdef}
	I_r(\tilde x,\tilde t):= 
	(r/2)^{\frac{4}{p-1}-n}
	\int_{\tilde t-r^2/4}^{\tilde t-r^2/16} \int_{B(\tilde x, r/2)}
	|u(x,t)|^{p+1} dxdt \leq \tilde \eps
\end{equation}
for any $r>0$, $\tilde x\in\R^n$ and $-1<\tilde t\leq 0$ 
satisfying $Q( (\tilde x,\tilde t), r)\subset Q(\delta/2)$, 
then 
\[
	\|u\|_{L^\infty(Q(\delta/4))}
	\leq C\tilde \eps^\frac{1}{p+1} \delta^{-\frac{2}{p-1}}, 
	\quad 
	\|\nabla u\|_{L^\infty(Q(\delta/4))} 
	\leq C\tilde \eps^\frac{1}{p+1} \delta^{-\frac{p+1}{p-1}}. 
\]
\end{lemma}

\begin{proof}
Let $\tilde \eps_1>0$ be a constant chosen later and let $0<\tilde \eps<\tilde \eps_1$.  
Set $v(y,s):= \delta^{2/(p-1)} u(\delta y, \delta^2 s)$. 
Note that $v_t=\Delta v+|v|^{p-1}v$ in $\R^n\times(-4,0)$, 
since $-4\delta^2>-1$. 
If $Q( (\tilde x,\tilde t), r)\subset Q(1/2)$ and $(r/4)^2\leq -\tilde t$, then 
$Q( (\tilde x,\tilde t+(r/4)^2), r)\subset Q(1/2)$, 
and so the change of variables and \eqref{eq:Irdef} show that 
\begin{equation}\label{eq:r24p1}
\begin{aligned}
	&(r/2)^{\frac{4}{p-1}-n} \int_{\tilde t-r^2/16}^{\tilde t} 
	\int_{B(\tilde x, r/2)} |v(y,s)|^{p+1} dyds \\
	&\leq 
	(r/2)^{\frac{4}{p-1}-n} 
	\int_{(\tilde t+r^2/16)-r^2/4}^{(\tilde t+r^2/16)-r^2/16} 
	\int_{B(\tilde x, r/2)} |v(y,s)|^{p+1} dyds \leq \tilde \eps. 
\end{aligned}
\end{equation}

Let $(x',t')\in Q(1/3)$ and 
$\lambda:=(-t')^{1/2}$. 
Then $\tilde v(x,t)
:=\lambda^{2/(p-1)} v(\lambda x+x', \lambda^2 t+t')$ satisfies 
$\tilde v_t=\Delta \tilde v+|\tilde v|^{p-1} \tilde v$ 
in $\R^n\times ( -4/\lambda^2, 0)$. 
If $B( (\tilde x-x')/\lambda, r/\lambda)\subset B(1/2)$ 
and $((\tilde t-r^2)/\lambda^2, \tilde t/\lambda^2)\subset ( -1/4, 0)$, then 
$Q( (\tilde x, \tilde t + t'), r)\subset Q(1/2)$ and 
$(r/4)^2\leq -(\tilde t + t')$. 
Therefore, \eqref{eq:r24p1} shows that 
\[
\begin{aligned}
	&(r/4\lambda)^{\frac{4}{p-1}-n}
	\int_{(\tilde t-(r/4)^2)/\lambda^2}^{\tilde t/\lambda^2} 
	\int_{ B((\tilde x-x')/\lambda, r/4\lambda) }
	|\tilde{v}(x,t)|^{p+1} dxdt \\
	&=(r/4)^{\frac{4}{p-1}-n} 
	\int_{\tilde t+t'-(r/4)^2}^{\tilde t+t'}
	\int_{B(\tilde x, r/4)}
	|v(y, s)|^{p+1} dyds 
	\leq C \tilde \eps
\end{aligned}
\]
if $Q( ((\tilde x-x')/\lambda, \tilde t/\lambda^2), r/\lambda) 
\subset Q(1/2)$. 
Replacing $(\tilde x,\tilde t)$ and $r$ with 
$(\lambda \tilde x+x',\lambda^2 \tilde t)$ and $\lambda r$, 
respectively, we see that 
\[
	(r/4)^{\frac{4}{p-1}-n}
	\iint_{Q((\tilde x,\tilde t), r/4)} |\tilde{v}(x,t)|^{p+1} dxdt 
	\leq C \tilde \eps
\]
if $Q( (\tilde x,\tilde t), r) \subset Q(1/2)$. 
Hence $\|\tilde{v}\|_{M^{p+1,2(p+1)/(p-1)}(Q(1/3))} 
\leq C \tilde \eps^{1/(p+1)}$, 
where $\|\cdot\|_{M^{p+1,2(p+1)/(p-1)}(Q(1/3))}$ 
is the parabolic Morrey norm on the backward parabolic cylinder $Q(1/3)$, 
see \cite[Section 2]{BS15} for the definition. 
Taking $\tilde \eps_1$ ($>\tilde \eps$) 
sufficiently small depending only on $n$ and $p$, 
we apply \cite[Proposition 4.1]{BS15} to see that 
\[
	\|\tilde v\|_{L^\infty(Q(1/4))} 
	+ \|\nabla \tilde v\|_{L^\infty(Q(1/5))} 
	\le C \|\tilde v \|_{M^{p+1,\frac{2(p+1)}{p-1}} (Q(1/3))} 
	\leq C\tilde \eps^\frac{1}{p+1}, 
\]
and so 
\[
	\lambda^\frac{2}{p-1} |v(\lambda x+x', \lambda^2 t+t')| 
	+ \lambda^\frac{p+1}{p-1} |\nabla v(\lambda x+x', \lambda^2 t+t')| 
	\leq C\tilde \eps^\frac{1}{p+1}
\]
for $(x,t)\in Q(1/5)$. 
Letting $(x,t)\to (0,0)$ and using $\lambda=(-t')^{1/2}$ yield  
$|v(x', t' )| \leq C\tilde \eps^{1/(p+1)} (-t')^{-1/(p-1)}$ 
and 
$|\nabla v(x', t' )| \leq 
C\tilde \eps^{1/(p+1)} (-t')^{-(p+1)/(2(p-1))}$. 
Recall that $(x',t')\in Q(1/3)$.  Hence we obtain 
\begin{align}
	&\label{eq:GKas}
	|v(y,s)| \leq C\tilde \eps^\frac{1}{p+1} (-s)^{-\frac{1}{p-1}} 
	&& \mbox{ for }y\in B(1/3), \; -1/9 \leq s <0,  \\ 
	&\label{eq:GKnabas}
	|\nabla v(y,s)| \leq C\tilde \eps^\frac{1}{p+1} (-s)^{-\frac{p+1}{2(p-1)}} 
	&& \mbox{ for }y\in B(1/3), \; -1/9 \leq s <0.  
\end{align}

We give an $L^\infty$ bound of $v$ 
based on the argument of \cite[Theorem 2.1]{GK89}. 
Let $\phi\in C_0^\infty(\R^n)$ satisfy 
$\phi=1$ in $B(7/24)$, 
$\phi=0$ in $\R^n\setminus B(1/3)$ and 
$0\leq \phi\leq 1$. 
Set $w:= \phi v$. Then $w$ satisfies 
\[
	w_t -\Delta w = |v|^{p-1} w -2\nabla \cdot (v\nabla \phi) + v\Delta \phi. 
\]
From the Duhamel formula and 
convolution inequalities on the heat kernel $G$ and $\nabla G$,
together with \eqref{eq:GKas} and $(t-s)^{-1/2}>1$ for $-1/9<s<t<0$, 
it follows that 
\begin{equation}\label{eq:wcom19}
\begin{aligned}
	&\| w(\cdot,t)\|_{L^\infty(\R^n)} \\
	&\leq \|w(\cdot,-1/9)\|_{L^\infty(\R^n)} 
	+\int_{-1/9}^t \|v(\cdot,s)\|_{L^\infty(B(1/3))}^{p-1} 
	\|w(\cdot,s)\|_{L^\infty(\R^n)} ds \\
	&\quad 
	 +C \int_{-1/9}^t (t-s)^{-\frac{1}{2}} 
	\|v(\cdot,s)\|_{L^\infty(B(1/3))} ds
	 \\
	&\leq 
	C\tilde \eps^\frac{1}{p+1} 
	\left( 1+ \int_{-1/9}^t (t-s)^{-\frac{1}{2}} (-s)^{-\frac{1}{p-1}} ds\right)\\
	&\quad+
	C\tilde \eps^\frac{p-1}{p+1} \int_{-1/9}^t (-s)^{-1} \|w(\cdot,s)\|_{L^\infty(\R^n)} ds. 
\end{aligned}
\end{equation}

Similar computations in \cite[Lemma 2.2]{GK89} show that 
\begin{equation}\label{eq:GK89l22}
	\int_{-1/9}^t (t-s)^{-\frac{1}{2}} (-s)^{-\frac{1}{p-1}} ds 
	\leq \left\{ \begin{aligned}
	&C\left(  1+ \log \frac{1}{-t} \right) 
	&&\mbox{ if }\frac{1}{2}\geq \frac{1}{p-1}, \\
	&C(-t)^{-(\frac{1}{p-1}-\frac{1}{2})} &&\mbox{ if }\frac{1}{2}< \frac{1}{p-1}. 
	\end{aligned}
	\right.  
\end{equation}
By setting 
\[
	\alpha:=\left\{ 
	\begin{aligned}
	&\frac{1}{4(p-1)} &&\mbox{ if }\frac{1}{2}\geq \frac{1}{p-1}, \\
	&\frac{1}{p-1}-\frac{1}{2} &&\mbox{ if }\frac{1}{2} < \frac{1}{p-1}, 
	\end{aligned}
	\right.
\]
we see that 
\[
	\int_{-1/9}^t (t-s)^{-\frac{1}{2}} (-s)^{-\frac{1}{p-1}} ds 
	\leq 
	C(-t)^{-\alpha}. 
\]
Then there exist $C_1,C_2>0$ depending only on $n$ and $p$ such that 
\[
\begin{aligned}
	\|w(\cdot,t)\|_{L^\infty(\R^n)} 
	\leq 
	C_1\tilde \eps^\frac{1}{p+1} (-t)^{-\alpha}
	+
	C_2\tilde \eps^\frac{p-1}{p+1} 
	\int_{-1/9}^t (-s)^{-1} \|w(\cdot,s)\|_{L^\infty(\R^n)} ds. 
\end{aligned}
\]
We fix $\tilde \eps_1$ such that 
$C_2\tilde \eps_1^{(p-1)/(p+1)}\leq \min((4(p-1))^{-1}, 4^{-1})$.

Since $t\mapsto (-t)^{-\alpha}$ is nondecreasing, 
Gronwall's inequality yields 
\[
\begin{aligned}
	&\|w(\cdot,t)\|_{L^\infty(\R^n)}  
	\leq C_1\tilde \eps^\frac{1}{p+1} (-t)^{-\alpha} 
	\exp\left( C_2\tilde \eps^\frac{p-1}{p+1} \int_{-1/9}^t (-s)^{-1} ds \right) \\
	&\leq 
	C_1 \tilde \eps^\frac{1}{p+1} (-t)^{-\alpha} 
	\exp\left( C_2 \tilde \eps_1^\frac{p-1}{p+1} \log\frac{1}{-t} 
	\right) \\
	&= 
	C_1\tilde \eps^\frac{1}{p+1} (-t)^{-\alpha - C_2\tilde \eps_1^\frac{p-1}{p+1}} 
	\leq 
	\left\{ \begin{aligned}
	&C_1 \tilde \eps^\frac{1}{p+1} (-t)^{-\frac{1}{2(p-1)} } 
	&&\mbox{ if } \frac{1}{2}\geq \frac{1}{p-1},  \\
	&C_1 \tilde \eps^\frac{1}{p+1} (-t)^{-\frac{1}{p-1}+\frac{1}{4}} 
	&&\mbox{ if } \frac{1}{2}< \frac{1}{p-1}.  \\
	\end{aligned} \right.
\end{aligned}
\]
We observe that the above estimate is an improvement of \eqref{eq:GKas}. 
In the first case ($1/2\geq 1/(p-1)$), 
we can bootstrap the above argument one more time 
to get the required boundedness of $v$ on $Q(1/4)$. 
In the second case ($1/2< 1/(p-1)$), 
we can apply finitely many bootstraps to enter 
a similar scenario to the first case. 
In all scenarios, we get
\[
	\|v\|_{L^\infty(Q(1/4))} \leq 
	\|w\|_{L^\infty(\R^n\times (-1/9,0))}
	\leq C\tilde \eps^\frac{1}{p+1}, 
\]
and so 
$\| u\|_{L^\infty(Q(\delta/4))} 
\leq C\tilde \eps^{1/(p+1)} \delta^{-2/(p-1)}$.

Let us next consider an $L^\infty$ bound of $\nabla v$. 
We note that $i$-th derivative $v_{y_i}$ satisfies 
$\partial_t v_{y_i}=\Delta v_{y_i} + p|v|^{p-1} v_{y_i}$. 
Setting $\tilde w:=\phi v_{y_i}$ gives 
\[
	\tilde w_t -\Delta \tilde w 
	= p|v|^{p-1} \tilde w -2\nabla \cdot (v_{y_i} \nabla \phi) 
	+ v_{y_i} \Delta \phi. 
\]
Similar computations to \eqref{eq:wcom19} with \eqref{eq:GKnabas} show that 
\[
\begin{aligned}
	&\| \tilde w(\cdot,t) \|_{L^\infty(\R^n)} \\
	&\leq 
	\|\tilde w(\cdot,-1/9)\|_{L^\infty(\R^n)} 
	+ p\int_{-1/9}^t \|v(\cdot,s)\|_{L^\infty(B(1/3))}^{p-1} 
	\| \tilde w(\cdot,s)\|_{L^\infty(\R^n)} ds  \\
	&\quad
	+ C\int_{-1/9}^t (t-s)^{-\frac{1}{2}} 
	\| v_{y_i}(\cdot,s) \|_{L^\infty(B(1/3))} ds \\
	&\leq 
	C \tilde \eps^\frac{1}{p+1} 
	\left( 1+ \int_{-1/9}^t (t-s)^{-\frac{1}{2}} 
	(-s)^{-\frac{p+1}{2(p-1)}} ds\right)\\
	&\quad+
	C\tilde \eps^\frac{p-1}{p+1} \int_{-1/9}^t (-s)^{-1} 
	\|\tilde w(\cdot,s)\|_{L^\infty(\R^n)} ds. 
\end{aligned}
\]
Then by an analog of \eqref{eq:GK89l22}, we have 
\[
	\|\tilde w(\cdot,t)\|_{L^\infty(\R^n)} 
	\leq 
	C \tilde \eps^\frac{1}{p+1} (-t)^{-\frac{1}{p-1} } 
	+ C\tilde \eps^\frac{p-1}{p+1} \int_{-1/9}^t (-s)^{-1} 
	\|\tilde w(\cdot,s)\|_{L^\infty(\R^n)} ds. 
\]
From similar computations as for $w$  
with $\tilde \eps_1$ ($>\tilde \eps$)  replaced by a smaller constant 
depending only on $n$ and $p$ if necessary, 
it follows that 
$\|v_{y_i}\|_{L^\infty(Q(1/4))} 
\leq C\tilde \eps^{1/(p+1)}$ and that 
$\|\nabla u\|_{L^\infty(Q(\delta/4))} 
\leq C\tilde \eps^{1/(p+1)} \delta^{-(p+1)/(p-1)}$. 
The proof is complete. 
\end{proof}

\begin{remark}
Note that with a slight adjustment to the test function $\phi$ 
and the use of another test function $\tilde{\phi}$ 
in defining $\tilde{w}:=\tilde{\phi} v_{y_{i}}$, 
we could get that $v$ is bounded in $Q(1/3)$ and could 
then infer that for $-1/9<t<0$, 
\[
	\|\tilde{w}(\cdot,t)\|_{L^{\infty}(\R^n)}
	\leq C\tilde \eps^\frac{1}{p+1}  (-t)^{-\frac{1}{p-1}}. 
\]
We could subsequently bootstrap this to get boundedness of $\tilde{w}$ in $Q(1/4)$. We do not pursue this for notational convenience.
\end{remark}

We prove the desired quantitative $\eps$-regularity. 
This is achieved by using the Giga-Kohn monotonicity formula (Lemma \ref{lem:Emono}), 
uniform bounds on the Giga-Kohn energy (Lemma \ref{lem:Ebdd}) and the previous lemma. 
See also related arguments in \cite{MTap, MTpre, MTno}.

\begin{proof}[Proof of Proposition \ref{pro:epsreg}]
Assume \eqref{eq:1scale}. 
By translation invariance and scaling, it suffices to prove the case $(x_0,t_0)=(0,0)$. 
For $r>0$, $\tilde x\in\R^n$ and $-1<\tilde t\leq 0$ 
with $Q((\tilde x,\tilde t),r)\subset Q(\delta/2)$, 
we note that $|\tilde x|\leq \delta/2$, $-\delta^2/4\leq \tilde t\leq 0$ 
and $r\leq \delta/2$ hold. 
Define $I_r(\tilde x,\tilde t)$ as in \eqref{eq:Irdef}. 
By using the backward similarity variables, 
\[
\begin{aligned}
	I_r(\tilde x,\tilde t) &\leq 
	C \int_{\tilde t-r^2/4}^{\tilde t-r^2/16} 
	(\tilde t-t)^\frac{2}{p-1} 
	\int_{B(\tilde x,r/2)} |u(x,t)|^{p+1} (\tilde t-t)^{-\frac{n}{2}} 
	e^{-\frac{|x-\tilde x|^2}{4(\tilde t-t)}}
	dxdt \\
	&\leq 
	C \int_{-\log(r^2/4)}^{-\log(r^2/16)} \int_{\R^n} 
	|w_{(\tilde x,\tilde t)}(\eta,\tau)|^{p+1} \rho(\eta) d\eta d\tau
	=: J_r(\tilde x,\tilde t). 
\end{aligned}
\]
From \eqref{eq:w22} and  \eqref{eq:GKmn}, it follows that 
for $w=w_{(\tilde x,\tilde t)}$, 
\[
\begin{aligned}
	J_r&= \frac{p+1}{p-1} \int_{-\log(r^2/4)}^{-\log(r^2/16)} 
	\int_{\R^n} w w_\tau \rho d\eta d\tau 
	+ \frac{2(p+1)}{p-1} \int_{-\log(r^2/4)}^{-\log(r^2/16)} 
	\cE_{(\tilde x,\tilde t)} d\tau  \\
	& \leq 
	C \left(\iint |w|^2 \rho  \right)^\frac{1}{2}
	\left(\int_{-\log(r^2/4)}^{-\log(r^2/16)} 
	\int_{\R^n} |w_\tau|^2 \rho  \right)^\frac{1}{2}
	+ C \int_{-\log(r^2/4)}^{-\log(r^2/16)} 
	\cE_{(\tilde x,\tilde t)} d\tau  \\
	& \leq 
	C  ( J_r )^\frac{1}{p+1} 
	\left(\int_{-\log(r^2/4)}^{-\log(r^2/16)} 
	\int_{\R^n} \rho  \right)^{\frac{1}{2}-\frac{1}{p+1}} 
	\left(\cE(-\log(r^2/4)) - \cE(-\log(r^2/16)) \right)^\frac{1}{2}  \\
	&\quad 
	+ C\cE(-\log(r^2/4))
	\int_{-\log(r^2/4)}^{-\log(r^2/16)}  d\tau. 
\end{aligned}
\]
Here and below, we often suppress the subscript $(\tilde x, \tilde t)$. 
On the other hand, straightforward computations with 
\eqref{eq:Type1bound2.0} yield 
\[
\begin{aligned}
	J_r&= 
	\int_{\tilde t-r^2/4}^{\tilde t-r^2/16} 
	(\tilde t-t)^{\frac{2}{p-1}-\frac{n}{2}} 
	\left\| |u(\cdot,t)|^{p+1}  
	e^{-\frac{|\cdot-\tilde x|^2}{4(\tilde t-t)}}
	\right\|_{L^1(\R^n)}  dt \\
	&\leq 
	C \int_{\tilde t-r^2/4}^{\tilde t-r^2/16} 
	(\tilde t-t)^{\frac{2}{p-1}-\frac{n}{2}} 
	\| u(\cdot,t) \|_{L^{q_c,\infty}(\R^n)}^{p+1}  
	\left\| e^{-\frac{|\cdot-\tilde x|^2}{4(\tilde t-t)}}
	\right\|_{L^{\frac{q_c}{q_c-(p+1)},1}(\R^n)}  dt \\
	&\leq 
	C M^{p+1} \int_{\tilde t-r^2/4}^{\tilde t-r^2/16} 
	(\tilde t-t)^{\frac{2}{p-1}-\frac{n}{2}+\frac{n}{2}( 1-\frac{p+1}{q_c} ) } dt
	= 
	CM^{p+1} \log \frac{r^2/4}{r^2/16}
	\leq CM^{p+1}. 
\end{aligned}
\]
This together with $\cE_{(\tilde x,\tilde t)}(-\log(r^2/16))\geq0$, $M>1$, 
\eqref{eq:EcEeqi} and \eqref{eq:GKmnori} gives 
\[
\begin{aligned}
	J_r &\leq 
	C M\left( \cE(-\log(r^2/4))  \right)^\frac{1}{2}
	+ C \cE( -\log(r^2/4) ) \leq 
	C M f( E(\tilde t - r^2/4) ) \\
	&\leq 
	C M f( E(-\delta^2/2) )
	\leq 
	C Mf\left( \frac{1}{\delta^2} 
	\int_{-\delta^2}^{-\delta^2/2} 
	E(t) dt \right), 
\end{aligned}
\]
where $f(s):=s+s^{1/2}$ for $s\geq 0$.

Recall \eqref{eq:EcEeqi}. Then, 
\[
	\frac{1}{\delta^2} \int_{-\delta^2}^{-\delta^2/2} 
	E(t) dt
	= 
	\int_{-\log(\tilde t+\delta^2)}^{-\log(\tilde t+\delta^2/2)} 
	\frac{e^{-\tau} }{\delta^2} \cE(\tau) d\tau 
	\leq 
	C \int_{-\log(\tilde t+\delta^2)}^{-\log(\tilde t+\delta^2/2)} 
	\cE(\tau) d\tau. 
\]
By \eqref{eq:w22}, we can also see that 
\[
\begin{aligned}
	&\int_{-\log(\tilde t+\delta^2)}^{-\log(\tilde t+\delta^2/2)}  
	\cE(\tau) d\tau 
	= 
	- \frac{1}{2} \iint w w_\tau \rho  
	+ \frac{p-1}{2(p+1)} \iint |w|^{p+1} \rho  \\
	& \leq 
	C \left(\iint |w|^2 \rho  \right)^\frac{1}{2}
	\left(\int_{-\log(\tilde t+\delta^2)}^{-\log(\tilde t+\delta^2/2)} 
	\int_{\R^n} |w_\tau|^2 \rho  \right)^\frac{1}{2}
	+ C \iint|w|^{p+1} \rho  \\
	& \leq 
	C \left(\iint |w|^{p+1} \rho  \right)^\frac{1}{p+1} 
	\left(\int_{-\log(\tilde t+\delta^2)}^{-\log(\tilde t+\delta^2/2)} 
	\int_{\R^n} \rho  \right)^{\frac{1}{2}-\frac{1}{p+1}} \\
	&\quad \times 
	\left(\cE(-\log(\tilde t+\delta^2)) 
	- \cE(-\log(\tilde t+\delta^2/2)) \right)^\frac{1}{2}
	+ C \iint|w|^{p+1} \rho.  
\end{aligned}
\]
Since $0\leq \cE_{(\tilde x,\tilde t)}(-\log(\tilde t+\delta^2)) 
= E_{(\tilde x,\tilde t)}(-\delta^2)\leq CM^{2p}$ by Lemma \ref{lem:Ebdd}, 
we obtain 
\[
\begin{aligned}
	&\int_{-\log(\tilde t+\delta^2)}^{-\log(\tilde t+\delta^2/2)}  
	\cE_{(\tilde x,\tilde t)}(\tau) d\tau 
	\leq 
	C M^p \tilde f\left( 
	\int_{-\log(\tilde t+\delta^2)}^{-\log(\tilde t+\delta^2/2)}  \int_{\R^n} 
	|w_{(\tilde x,\tilde t)}(\eta,\tau)|^{p+1} \rho(\eta) d\eta d\tau \right) \\
	&\leq 
	CM^p \tilde f\left( 
	\int_{-\delta^2}^{-\delta^2/2} 
	(\tilde t-t)^\frac{2}{p-1} 
	\int_{\R^n} |u(x,t)|^{p+1} K_{(\tilde x,\tilde t)}(x,t) dxdt \right), 
\end{aligned}
\]
where $\tilde f(s):=s+s^{1/(p+1)}$ for $s\geq0$. 

Let $0<\eps_0<1$ be a constant chosen later and let $0<\eps<\eps_0$.  
Then \eqref{eq:1scale} yields 
\[
	I_r(\tilde x, \tilde t) 
	\leq CM^{p+1} (f \circ \tilde f)\left( 
	\eps + \delta^\frac{4}{p-1} 
	\int_{-\delta^2}^{-\delta^2/2} 
	\int_{\R^n\setminus B(A\delta)} |u(x,t)|^{p+1} K_{(\tilde x,\tilde t)}(x,t) dxdt 
	\right). 
\]
We use the notation $K_{(\tilde y,\tilde s)}$ 
even for $\tilde s>0$ and observe that 
if $Q( (\tilde x,\tilde t), r)\subset Q(\delta/2)$, then 
\[
	K_{(\tilde x,\tilde t)} (x,t) 
	\leq 
	CK_{(0,\delta^2/2)} (x,t)
	\leq 
	Ce^{-\frac{A^2}{24}} K_{(0,\delta^2)} (x,t)
\]
for $x\in \R^n\setminus B(A\delta)$ and $-\delta^2 < t < -\delta^2/2$, 
and thus 
\[
\begin{aligned}
	I_r 
	&\leq 
	CM^{p+1} (f \circ \tilde f)\left( 
	\eps + \delta^\frac{4}{p-1} e^{-\frac{A^2}{24}}
	\int_{-\delta^2}^{-\delta^2/2} 
	\int_{\R^n} |u|^{p+1} K_{(0,\delta^2)} dxdt 
	\right). 
\end{aligned}
\]
By the H\"older inequality and \eqref{eq:Type1bound2.0}, we have 
\[
\begin{aligned}
	&\int_{-\delta^2}^{-\delta^2/2} 
	\int_{\R^n} |u|^{p+1} K_{(0,\delta^2)} dxdt \\
	&\leq 
	C \int_{-\delta^2}^{-\delta^2/2} 
	(\delta^2-t)^{-\frac{n}{2}}
	\| |u|^{p+1}\|_{L^{\frac{q_c}{p+1},\infty}(\R^n)} 
	\left\| e^{-\frac{|\cdot|^2}{4(\delta^2-t)}} 
	\right\|_{L^{\frac{q_c}{q_c-(p+1)},1}(\R^n)} dt \\
	&\leq 
	CM^{p+1} \int_{-\delta^2}^{-\delta^2/2} 
	(\delta^2-t)^{-\frac{p+1}{p-1}} dt 
	= CM^{p+1} \delta^{-\frac{4}{p-1}}. 
\end{aligned}
\]
Therefore, if $0<\eps<\eps_0$, 
$Q( (\tilde x,\tilde t), r)\subset Q(\delta/2)$ 
and $A$ is chosen as in \eqref{eq:1scale}, then 
\[
\begin{aligned}
	I_r 
	&\leq 
	CM^{p+1} (f \circ \tilde f)( 
	\eps + C   M^{p+1}  e^{-\frac{A^2}{24}}) \leq 
	CM^{p+1} (f \circ \tilde f)( \eps + C  \eps ) \\
	&\leq 
	CM^{p+1} (f \circ \tilde f)( \eps ) 
	\leq \tilde C M^{p+1}\eps^\frac{1}{2(p+1)}  
	\leq 
	\tilde C M^{p+1}\eps_0^\frac{1}{2(p+1)}  
	= \tilde \eps_1 
\end{aligned}
\]
with some constants $\tilde C>0$ depending only on $n$ and $p$. 
Here, $\tilde \eps_1$ is given in Lemma \ref{lem:epspre} and 
$\eps_0$ is chosen such that 
$\eps_0:=\tilde C^{-2(p+1)}M^{-2(p+1)^2} \tilde \eps_1^{2(p+1)}$. 
Thus, we can apply Lemma \ref{lem:epspre} 
with $\tilde \eps=\tilde C M^{p+1}\eps^{1/(2(p+1))}$ 
to obtain 
\[
	\|u\|_{L^\infty( Q(\delta/4)) } \leq 
	C M \eps^\frac{1}{2(p+1)^2} \delta^{-\frac{2}{p-1}}, 
	\quad 
	\|\nabla u\|_{L^\infty( Q(\delta/4)) } \leq 
	C M \eps^\frac{1}{2(p+1)^2} \delta^{-\frac{p+1}{p-1}}. 
\]
The proof is complete, with $\eps_1$ in Proposition \ref{pro:epsreg} 
defined as $\eps_1:=\tilde \eps_1/\tilde C$. 
\end{proof}

\section{Estimates via $\eps$-regularity}\label{sec:estiveps}
Let $u$ be a classical solution of \eqref{eq:fujitaeq} 
satisfying the Lorentz norm bound \eqref{eq:Type1bound2.0}. 
We give several estimates based on the above quantitative $\eps$-regularity. 
We set $\eps_0:= M^{-2(p+1)^2} \eps_1^{2(p+1)}$ 
as in Proposition \ref{pro:epsreg}, where $\eps_1$ depends only on $n$ and $p$. 
In the rest of this paper, we fix 
\begin{equation}\label{eq:epsfixed}
	\eps:=\frac{\eps_{0}}{2}, \quad A:=48\log \frac{M^{p+1}}{\eps}, 
\end{equation}
unless otherwise stated. 

\subsection{Propagation of concentration}
We prove the backward-in-time propagation of concentration 
in the form of its contraposition. 
After that, we give estimates on 
intersecting regions of concentration and quantitative regularity.

\begin{proposition}[Backward propagation]\label{pro:backprop}
Let $p>p_S$ and assume \eqref{eq:epsfixed}. 
There exists $M_*>1$ depending only on $n$ and $p$ such that 
the following statement holds true. 
Assume that $u$ satisfies \eqref{eq:Type1bound2.0} with $M\geq M_{*}$. 
If $t', t''\in (-1,0)$ and $u$ satisfy
\begin{align}
	&\label{eq:wellseptimes}
	-\frac{1}{16}\leq t''\leq M^{(p-1)(2p+4)}  A^{\frac{n(p-1)}{p+1}} t', \\
	&\label{eq:backwardsmallness}
	(-t'')^{\frac{2}{p-1}-\frac{n}{2}}
	\int_{t''}^{t''/2} \int_{ B(x_0, A(-t'')^{\frac{1}{2}} )} 
	|u(x,t)|^{p+1} dxdt\leq \eps, 
\end{align}
then we conclude that 
\begin{equation}\label{eq:forwardsmallness}
	(-t')^{\frac{2}{p-1}-\frac{n}{2}}
	\int_{t'}^{t'/2} \int_{ B(x_0, A(-t')^{\frac{1}{2}}) } 
	|u(x,t)|^{p+1} dxdt\leq \eps. 
\end{equation}
\end{proposition}

\begin{proof}
Without loss of generality set $x_0=0$.
The assumptions \eqref{eq:wellseptimes} and \eqref{eq:backwardsmallness} 
allow us to apply Proposition \ref{pro:epsreg}, which gives that 
\[
	\|u\|_{L^\infty( Q((-t'')^{\frac{1}{2}}/4)) } 
	\leq 
	CM \eps^{\frac{1}{2(p+1)^2}} (-t'')^{-\frac{1}{p-1}}, 
\]
and so
\begin{equation}\label{eq:backproppointwiseest}
	\sup_{Q( \frac{1}{4} (-t'')^{\frac{1}{2}} )}|u(x,t)|^{p+1}
	\leq C^{p+1}M^{p+1} \eps^{\frac{1}{2(p+1)}} (-t'')^{-\frac{p+1}{p-1}}.
\end{equation}
From \eqref{eq:wellseptimes}, we have that for $M$ sufficiently large, 
\[
	B(A(-t')^{\frac{1}{2}})
	\times (t',t'/2)\subset  Q( (-t'')^{\frac{1}{2}}/4).
\]
Thus, \eqref{eq:wellseptimes} and \eqref{eq:backproppointwiseest} 
imply that 
\[
\begin{aligned}
	&(-t')^{\frac{2}{p-1}-\frac{n}{2}}
	\int_{t'}^{t'/2} 
	\int_{B(A(-t')^{\frac{1}{2}})} |u(x,t)|^{p+1} dxdt 
	\leq C (M^{p+1}A^n)\left(\frac{-t'}{-t''}\right)^{\frac{p+1}{p-1}} \\
	&\leq C M^{-(2p+3)(p+1)}\leq \eps, 
\end{aligned}
\]
as required.
\end{proof}

From now on we will use the terminology `$M$ being sufficiently large' 
to mean that $M\geq M_0$, where $M_0>1$ is a sufficiently large constant 
depending only on $n$ and $p$. 
We will also utilize the notation 
$M_{k}:=M^{c_k(n,p)}$ for $k\in\N$  
with $M$ being sufficiently large 
and with $c_k(n,p)>0$ being sufficiently 
larger than $c_{k-1}(n,p)$ etc. 

The rest of this subsection is devoted to 
estimates on intersecting regions of concentration and quantitative regularity. 
A related statement for the higher-dimensional Navier-Stokes equations 
was previously proven in \cite[Proposition 5.1]{Pa22}.

\begin{proposition}[Intersecting regions]\label{pro:counting}
Let $p>p_S$ and assume \eqref{eq:epsfixed}. 
Then there exist $C, \bar{C}>0$ depending only on $n$ and $p$ 
such that the following statement holds for all $M$ sufficiently large: 
Assume that $u$ satisfies \eqref{eq:Type1bound2.0} and that 
$-1/64< t''<0$ satisfies 
\begin{equation}\label{eq:conccounting}
	(-t'')^{\frac{2}{p-1}-\frac{n}{2}}
	\int_{t''}^{t''/2} 
	\int_{
	B(0, A(-t'')^{\frac{1}{2}} )
	} |u(x,t)|^{p+1} dxdt> \eps.
\end{equation}
Then we conclude that there exist backward parabolic cylinders 
$Q( (x_*,t_*- r^2/8), \hat{\delta} r )$ and $Q(z_*, r)$ 
with $z_*=(x_*,t_*)$ satisfying 
\begin{equation}\label{eq:cubeinclusion}
	Q( (x_*,t_*- r^2/8), \hat{\delta} r )
	\subset Q(z_*, r) 
	\subset B(0, 20 M_{1}(-t'')^{\frac{1}{2}})
	\times (2t'', t''/4)
\end{equation}
with $r:=8M_{2}^{-(6n+2)}(-t'')^{1/2}$ and $\hat{\delta}:=M_{2}^{-1}$ 
such that 
\begin{align}
	&\label{eq:regcount}
	\|u\|_{L^{\infty}(Q( z_*, r/2))}
	\leq C M^{-2} r^{-\frac{2}{p-1}},
	\quad
	\|\nabla u\|_{L^{\infty}(Q( z_*, r/2))}
	\leq C M^{-2} r^{-\frac{p+1}{p-1}}, \\
	&\label{eq:conccount}
	\int_{Q((x_*,t_*- r^2/8), \hat{\delta} r)} |u|^2 dxdt
	\geq M_{2}^{-\bar{C}}(\hat{\delta} r)^{n+2-\frac{4}{p-1}}. 
\end{align}
\end{proposition}

To prove this proposition, we prepare several estimates 
based on the following interpolation inequality 
on a domain $\Omega\subset \R^n$.

\begin{lemma}\label{lem:interpolativepropertyLorentz}
Let $1\leq \tilde{p}<\tilde{r}<\tilde{q}\leq\infty$, 
$f\in L^{\tilde{p},\infty}(\Omega)\cap L^{\tilde{q},\infty}(\Omega)$  and
$0<\tilde \theta<1$ be such that
\[
	\frac{1}{\tilde{r}}=\frac{\tilde \theta}{\tilde{p}}
	+\frac{1-\tilde \theta}{\tilde{q}}
	\quad \left(\mbox{equivalently, } 
	\theta = \frac{ \frac{1}{\tilde r} -\frac{1}{\tilde q} }
	{ \frac{1}{\tilde p} - \frac{1}{\tilde q} }\right). 
\]
Then, $f\in L^{\tilde{r}}(\Omega)$ and 
\[
	\|f\|_{L^{\tilde{r}}(\Omega)}
	\leq 
	\left( \frac{\tilde{r}}{\tilde{r}-\tilde{p}}
	+\frac{\tilde{r}}{\tilde{q}-\tilde{r}} \right)^{\frac{1}{\tilde{r}}}
	\|f\|_{L^{\tilde{p},\infty}(\Omega)}^{\tilde \theta}
	\|f\|_{L^{\tilde{q},\infty}(\Omega)}^{1-\tilde \theta}.
\]
\end{lemma} 

\begin{proof}
See \cite[Proposition 1.1.14]{Graf08} for example. 
\end{proof}

By using Lemma \ref{lem:interpolativepropertyLorentz}, we also prepare estimates  
for $f:\Omega\times I\to \R$ with $\Omega\times I\subset \R^n\times \R$ bounded.

\begin{lemma}\label{lem:interpolavivecor}
Let $p>p_{S}$ and define $\phi_i=\phi_i(p,n)$ ($i=1,2,3$) by 
\begin{equation}\label{eq:phi123def}
	\phi_{1}:=\frac{ \frac{1}{p+1} -\frac{1}{q_c}}{ \frac{1}{p-1} - \frac{1}{q_c}}, 
	\quad 
	\phi_{2}:=\frac{\frac{1}{p-1}-\frac{1}{q_c}}{\frac{1}{2}-\frac{1}{q_c}}, 
	\quad 
	\phi_{3}:=\frac{\frac{1}{p+1}-\frac{1}{q_c}}{\frac{1}{3q_c/4}-\frac{1}{q_c}}. 
\end{equation}
Suppose that
\[
	\sup_{t\in I}\|f(\cdot,t)\|_{L^{q_{c},\infty}(\Omega)}\leq M.
\]
Then there exists $C>0$ depending only on $n$ and $p$ 
such that the following inequalities hold 
with $\iint=\iint_{\Omega\times I}$: 
\begin{enumerate}[(i)]
\item
$\displaystyle \iint|f|^{p+1} 
	\leq C M^{(1-\phi_{1})(p+1)}
	\left( \iint  
	|f|^{p-1}\right)^{\frac{(p+1)\phi_{1}}{p-1}}
	|I|^{1-\frac{(p+1)\phi_{1}}{p-1}}$.
\item 
$\displaystyle \iint |f|^{p-1} 
\leq 
\left\{ \begin{aligned}
	&\left( \iint |f|^2 \right)^{\frac{p-1}{2}} 
	|\Omega\times I|^{1-\frac{p-1}{2}} &&\mbox{ if }p-1\leq 2, \\
	&C M^{(1-\phi_{2})(p-1)}
	\left(\iint |f|^2 \right)^{\frac{(p-1)\phi_{2}}{2}} 
	|I|^{1-\frac{(p-1)\phi_{2}}{2}} &&\mbox{ if }p-1> 2. 
\end{aligned} \right.$
\item 
$\displaystyle \iint |f|^{p+1} 
\leq 
	CM^{(1-\phi_{3})(p+1)}
	\left(\iint |f|^{\frac{3q_c}{4}} 
	\right)^{\frac{4(p+1)\phi_{3}}{3q_c}} 
	|I|^{1-\frac{4(p+1)\phi_{3}}{3q_c}}$ 
if $\displaystyle \frac{3q_c}{4} < p+1$. 
\end{enumerate}
\end{lemma}

\begin{proof}
(i) We apply Lemma \ref{lem:interpolativepropertyLorentz} 
with $\tilde{p}=p-1$, $\tilde{q}=q_{c}$ and $\tilde{r}=p+1$ 
for each $t\in I$. 
Raising this to the power $p+1$ and integrating with respect to time gives that 
\begin{equation}\label{eq:fp+1est}
	\int_{I}\|f(\cdot,t)\|_{L^{p+1}(\Omega)}^{p+1} dt
	\leq C M^{(1-\phi_{1})(p+1)}
	\int_{I}\|f(\cdot,t)\|_{L^{p-1}(\Omega)}^{(p+1)\phi_{1}}dt.
\end{equation} 
Noting that 
\[
	\frac{1}{p+1}=\frac{\phi_{1}}{p-1}+\frac{1-\phi_{1}}{q_{c}} 
	> \frac{\phi_1}{p-1}, 
\]
we have that $(p+1)\phi_{1}< p-1$. 
This allows us to apply H\"{o}lder's inequality 
to the right-hand side of \eqref{eq:fp+1est} 
to obtain the desired conclusion. 

(ii) The case $p-1\leq 2$ follows immediately from H\"{o}lder's inequality. 
As for the case $p-1> 2$, 
we apply Lemma \ref{lem:interpolativepropertyLorentz} 
with $\tilde{p}=2$, $\tilde{q}=q_{c}$ and $\tilde{r}=p-1$. 
The rest of the proof is the same as (i). 

(iii) This follows from 
Lemma \ref{lem:interpolativepropertyLorentz} 
with $\tilde{p}=3q_c/4$, $\tilde{q}=q_{c}$ and $\tilde{r}=p+1$. 
\end{proof}

The following corollary of Proposition \ref{pro:epsreg} 
will be convenient in several places.

\begin{corollary}\label{cor:corepsgeneral}
Let $p>p_S$ and $C^{*}>0$. 
Then there exists $M_0> 1$ 
depending only on $n$, $p$ and $C^*$ such that  
the following statement holds for 
any $M\geq M_0$ and $0< \delta<(2M)^{-1}$: 
If $u$ satisfies \eqref{eq:Type1bound2.0} and 
\begin{equation}\label{eq:1scalegeneral} 
\begin{aligned}
	\delta^{\frac{4}{p-1}-n} 
	\int_{t_0-\delta^2}^{t_0} \int_{B(x_0,\delta)} 
	|u|^{p+1} dxdt 
	\leq 
	C^{*} M^{-(n+3-\frac{2(p+1)}{p-1}
	+6(p+1)^2+\frac{2(p+1)^3}{p-1})}, 
\end{aligned}
\end{equation}
then there exists $C>0$ depending only on $n$ and $p$ such that 
\[
	\|u\|_{L^\infty( Q((x_0,t_0), \delta/2)) } \leq 
	CM^{-3} \delta^{-\frac{2}{p-1}}, \quad 
	\|\nabla u\|_{L^\infty( Q((x_0,t_0), \delta/2)) } \leq 
	CM^{-2} \delta^{-\frac{p+1}{p-1}}. 
\]
\end{corollary}

\begin{proof} 
In this proof, let 
\[
	\eps:=M^{-(6(p+1)^2+\frac{2(p+1)^3}{p-1})}, \quad 
	A:=48\log \frac{M^{p+1}}{\eps}. 
\]
Clearly, we have that for $M$ sufficiently large, 
\begin{equation}\label{eq:epsilonepsilon0inequality}
	0<\eps<\eps_{0},
\end{equation}
where $\eps_{0}$ is as in Proposition \ref{pro:epsreg}.
Without loss of generality, let $x_0=0$. 
Note that for $M$ sufficiently large
\begin{equation}\label{eq:Msize}
	A\leq M/2.
\end{equation}

Let us rescale 
\[
	u_{\lambda}(x,t)
	:=\lambda^{\frac{2}{p-1}}u(\lambda x, \lambda^2 t+t_0)
	\quad\mbox{ with }
	\lambda={M}^{-1}.
\]
Then, $u_{\lambda}$ satisfies the scale-invariant bound 
\eqref{eq:Type1bound2.0} on the rescaled time interval. 
Let $\tilde{z}_{0}=(\tilde{t}_{0}, \tilde{x}_{0})
\in Q( (0,t_0), \delta/2)$. 
Using $0<\delta<(2M)^{-1}$ and \eqref{eq:Msize}, we observe that
\[
\begin{aligned}
	& \hat{z}_{0}=(\hat{x}_{0}, \hat{t}_{0})
	:=(\tilde{x}_{0}/\lambda, (\tilde{t}_{0}-t_0)/\lambda^2 )\in Q(M\delta/2), 
	\quad \hat{t}_{0}>-1/16, \\
	& B(\hat{x}_{0},A\delta)
	\times (\hat{t}_0-\delta^2/2, \hat{t}_{0}-\delta^2) \subset Q(M\delta). 
\end{aligned}
\]
These together with \eqref{eq:1scalegeneral}, \eqref{eq:epsilonepsilon0inequality} 
and \eqref{eq:Msize} show that
\[
\begin{aligned}
	&\delta^{\frac{4}{p-1}-n} 
	\int_{\hat{t}_0-\delta^2}^{\hat{t}_0-\delta^2/2} 
	\int_{B(\hat{x}_0, A\delta)} 
	|u_{\lambda}(x,t)|^{p+1} dxdt 
	\leq \delta^{\frac{4}{p-1}-n}\iint_{Q(M\delta)} 
	|u_{\lambda}(x,t)|^{p+1} dxdt\\
	&=M^{n+2-\frac{2(p+1)}{p-1}}\left( 
	\delta^{\frac{4}{p-1}-n} 
	\int_{t_0-\delta^2}^{t_0} \int_{B(\delta)} 
	|u(x,t)|^{p+1} dxdt\right)
	\leq C^{*}M^{-1}\eps.
\end{aligned}
\]
By \eqref{eq:epsilonepsilon0inequality} and $\hat{t}_{0}>-1/16$,  
we can apply Proposition \ref{pro:epsreg} to $u_\lambda$ to see that 
\[
\begin{aligned}
	&|u(\tilde{x}_{0}, \tilde{t}_{0})|
	\leq C M^{\frac{p+1}{p-1}}\eps^\frac{1}{2(p+1)^2} 
	\delta^{-\frac{2}{p-1}}
	= C M^{-3}\delta^{-\frac{2}{p-1}}, \\
	&|\nabla u(\tilde{x}_{0}, \tilde{t}_{0})|
	\leq C M^{\frac{p+1}{p-1}}\eps^\frac{1}{2(p+1)^2} 
	\delta^{-\frac{p+1}{p-1}}
	= C M^{-2}\delta^{-\frac{p+1}{p-1}}. 
\end{aligned}
\]
Since $\tilde{z}_{0}=(\tilde{x}_{0}, \tilde{t}_{0})
\in Q( (0,t_0), \delta/2)$ was taken arbitrarily, this implies that
\[
	\|u\|_{L^{\infty}(Q( (0,t_0), \delta/2))}
	\leq CM^{-3}\delta^{-\frac{2}{p-1}}, \quad 
	\|\nabla u\|_{L^{\infty}(Q( (0,t_0), \delta/2))}
	\leq CM^{-2}\delta^{-\frac{p+1}{p-1}}, 
\]
as required.
\end{proof}

We state a useful version of $\eps$-regularity, 
which will be used in the proof of Proposition \ref{pro:counting}.

\begin{lemma}\label{lem:epsgeneralforcounting}
Let $p>p_S$ and $C^*>0$. Define $\alpha(n,p):=\max(a_{1}(n,p), a_{2}(n,p))$ with 
\[
\begin{aligned}
	&a_{1}(n,p):=\frac{3q_c}{4}
	\left( \frac{n+3}{p+1} -\frac{2}{p-1}+6(p+1)+\frac{2(p+1)^2}{p-1} \right), \\
	&a_{2}(n,p):=\frac{3q_c}{4\phi_{3}}
	\left( \frac{n+3}{p+1} -\frac{2}{p-1}+6(p+1) 
	+\frac{2(p+1)^2}{p-1}+(1-\phi_{3}) \right), 
\end{aligned}
\]
where $\phi_{3}$ is defined by \eqref{eq:phi123def}.
Then there exists $M_0> 1$ 
depending only on $n$, $p$ and $C^*$ such that  
the following statement holds for 
any $M\geq M_0$ and $0<\delta<(2M)^{-1}$: 
If $u$ satisfies \eqref{eq:Type1bound2.0} and 
\[
	\delta^{-(\frac{n}{4}+2)} 
	\int_{t_0-\delta^2}^{t_0} \int_{B(x_0,\delta)} 
	 |u(x,t)|^{\frac{3q_c}{4}} dxdt \leq C^{*} M^{-\alpha(n,p)}, 
\]
then there exists $C>0$ depending only on $n$ and $p$ such that 
\[
	\|u\|_{L^\infty( Q( (x_0,t_0), \delta/2) ) } \leq 
	CM^{-3} \delta^{-\frac{2}{p-1}}, 
	\quad \|\nabla u\|_{L^\infty( Q( (x_0,t_0), \delta/2) ) } \leq 
	CM^{-2} \delta^{-\frac{p+1}{p-1}}.
\]
\end{lemma}

\begin{proof}
In the case that $3q_c/4 \geq p+1$, 
we obtain the conclusion by applying H\"{o}lder's inequality 
and Corollary \ref{cor:corepsgeneral}. 
In the case that $3q_c/4< p+1$, 
we apply Lemma \ref{lem:interpolavivecor} (iii) and Corollary \ref{cor:corepsgeneral}.
\end{proof}

We are now in a position to prove Proposition \ref{pro:counting} 
by using Lemma \ref{lem:interpolavivecor}, 
Lemma \ref{lem:epsgeneralforcounting} and 
performing a counting argument inspired by \cite[Proposition 5.1]{Pa22}. 
The scale-invariant bound \eqref{eq:Type1bound2.0} 
is not countably additive over disjoint spatial scales, 
which requires a further judicious choice of indices in the counting argument.

\begin{proof}[Proof of Proposition \ref{pro:counting}]
In the proof, 
$C_{n}^{*}$, $C_{n}^{**}$ and $C_{n}^{***}$ 
will be used to denote certain constants depending only on $n$. 
We divide the proof into 4 steps. 

\noindent\textbf{Step 1: rescaling and derived functions}

Let $-1/64< t''<0$ and let $u$ satisfy \eqref{eq:conccounting}. 
Set $I:=(-1/64, -1/128)$. 
Define the rescaled function 
$u_{\lambda}:\R^n\times (1/(64t''),0)\rightarrow\R$ by
\begin{equation}\label{eq:rescaleudefcount}
	u_{\lambda}(x,t)
	:= \lambda^{\frac{2}{p-1}} u(\lambda x, \lambda^2 t)
	\quad \mbox{ with }\lambda:= 8(-t'')^{\frac{1}{2}}.
\end{equation}
As $u$ satisfies \eqref{eq:Type1bound2.0} and \eqref{eq:conccounting}, we get that $u_{\lambda}$ satisfies
\begin{align}
	&\label{eq:ulambdatypeI}
	\sup_{-1<t<0}\|u_{\lambda}(\cdot,t)\|_{L^{q_{c},\infty}(\R^n)}\leq M, \\
	&\label{eq:ulambdaconccount}
	\iint_{B(A/8) \times I} 
	|u_{\lambda}(x,t)|^{p+1} dxdt\geq C(p)\eps.
\end{align}
Using \eqref{eq:ulambdatypeI}, \eqref{eq:ulambdaconccount} 
and Lemma \ref{lem:interpolavivecor} (i), we infer that
\begin{equation}\label{eq:ulambdaconccountp-1}
	\iint_{B(A/8)\times I} 
	|u_{\lambda}(x,t)|^{p-1} dxdt
	\geq C(p,n)M^{-\frac{(1-\phi_{1})(p-1)}{\phi_{1}}}
	\eps^{\frac{p-1}{(p+1)\phi_1}}, 
\end{equation}
where $\phi_{1}$ is as in \eqref{eq:phi123def}.
Define $U:\R^n\times (1/(64t''),0)\rightarrow\R$ by
\begin{equation}\label{eq:Udef}
	U(x,t):=|u_{\lambda}(x,t)|^{\frac{p-1}{2}}.
\end{equation}
By \eqref{eq:ulambdatypeI} and \eqref{eq:ulambdaconccountp-1}, 
we get that, for $M$ sufficiently large and appropriate $M_{1}$, 
\begin{align}
	&\label{eq:UtypeI}
	\sup_{-1<t<0}\|U(\cdot,t)\|_{L^{n,\infty}(\R^n)}\leq M^{\frac{p-1}{2}}, \\
	&\label{eq:UL2conc}
	\iint_{B(M_1) \times I} |U(x,t)|^2 dxdt\geq M_{1}^{-1}.
\end{align}

\noindent\textbf{Step 2: families of cubes and the regularity of $u_{\lambda}$}

For appropriately chosen $M_{2}:=M^{c_{2}(n,p)}$,  define 
\begin{equation}\label{eq:deltaldef}
	\hat{\delta}:=M_{2}^{-1}, \quad  l:=M_{2}^{-(6n+2)}.
\end{equation}
With this, we define the family of cubes (backward parabolic cylinders)
\[
	\mathcal{C}_{0}:=
	\{ Q(z_0,l); \; 
	z_0=(x_0,t_0) 
	\in ((\hat{\delta} l)^n \Z^n \times (\hat{\delta} l)^2\Z)
	\cap (B(2M_1) \times 2I) \}.
\]
Here, for $0<r\leq 3$, we write 
\[
	rI:= \left(-\frac{3}{256}-\frac{r}{256}, -\frac{3}{256}+\frac{r}{256} \right)
	\subset [-1,0].
\]
Then there exists $C(n)>1$ such that 
\begin{equation}\label{eq:sizecubesC0}
	\frac{1}{C(n)} M_{1}^n (\hat{\delta} l)^{-(n+2)}
	\leq  \#\mathcal{C}_0 
	\leq C(n) M_{1}^n (\hat{\delta} l)^{-(n+2)}, 
\end{equation}
where $\#\mathcal{C}_{0}$ is the cardinality of $\mathcal{C}_0$. 
Now, let us define the subfamily of cubes 
\[
	\mathcal{C}_{1}:=
	\left\{ Q\in\mathcal{C}_{0};\; 
	\|U\|_{L^{\frac{3n}{4}}(Q)}^{\frac{3n}{4}}>M_{1}^{-1} l^n \right\}.
\]
Here $M_{2}=M^{c_2(p,n)}$ is chosen appropriately 
to ensure that if
\[
	\frac{1}{l^{\frac{n}{4}+2}} 
	\|U\|_{L^{\frac{3n}{4}}(Q)}^{\frac{3n}{4}}
	\leq M_{1}^{-1} l^{\frac{3n}{4}-2}, 
\]
then $u_{\lambda}$ satisfies the hypothesis in Lemma \ref{lem:epsgeneralforcounting} with $\delta=l$, which in turn implies that
\begin{equation}\label{eq:ulambdabddC1}
	\|u_{\lambda}\|_{L^{\infty}(Q(z_0,l/2))}
	\leq C l^{-\frac{2}{p-1}}M^{-2}\quad\mbox{ with } Q=Q(z_0,l).
\end{equation}
Hence, to show that there exists a cube $Q\in\mathcal{C}_{0}$ 
satisfying \eqref{eq:ulambdabddC1}, 
it suffices to show that $\mathcal{C}_{0}\setminus\mathcal{C}_{1}\neq \emptyset$. 
To show this, it is sufficient to show that 
\begin{equation}\label{eq:C0greaterthanC1}
	\#\mathcal{C}_{1}< \#\mathcal{C}_{0}.
\end{equation}

We prove \eqref{eq:C0greaterthanC1}. 
Using the definition of $\mathcal{C}_{1}$, 
the fact that each of the cubes in $\mathcal{C}_1$ overlaps with at most 
$C(n) l^{n+2} (\hat \delta l)^{-2-n}$ ($= C(n) \hat{\delta}^{-2-n}$) 
other cubes in $\mathcal{C}_1$, 
\eqref{eq:UtypeI} and H\"{o}lder's inequality for Lorentz spaces, we have
\[
\begin{aligned}
	&M_{1}^{-1} l^n (\#\mathcal{C}_{1})
	\leq \sum_{Q\in\mathcal{C}_{0}} \|U\|^{\frac{3n}{4}}_{L^{\frac{3n}{4}}(Q)}
	\leq C(n) \hat{\delta}^{-2-n} 
	\iint_{B(3M_{1}) \times 3I} |U|^{\frac{3n}{4}} dxdt\\
	&\leq C(n) M_{1}^{\frac{n}{4}} M^{\frac{3n(p-1)}{8}} \hat{\delta}^{-2-n}.
\end{aligned}
\]
Thus,
\begin{equation}\label{eq:sizeC1upperbound}
	\#\mathcal{C}_{1}
	\leq C(n)  M_{1}^{\frac{n}{4}+1}
	M^{\frac{3n(p-1)}{8}} l^2  (\hat{\delta} l)^{-(n+2)}.
\end{equation}
Recalling \eqref{eq:sizecubesC0}, we see that for $M_{2}=M^{c_{2}(p,n)}$ chosen appropriately, we have that $\#\mathcal{C}_{1}<\#\mathcal{C}_{0}$. 
This implies that there exists a cube $Q\in\mathcal{C}_{0}$ 
such that \eqref{eq:ulambdabddC1} holds.

\noindent\textbf{Step 3: concentration of $u$ on a cube 
descending from $\mathcal{C}_{0}\setminus\mathcal{C}_{1}$}

For $Q=Q((x_0,t_0),l)\in\mathcal{C}_{0}$, 
we define the descendant 
$Q':=Q( (x_0,t_0-l^2/8), \hat{\delta}l)$. 
These cubes are such that 
$B(M_{1})  \times  I \subset \{Q'; Q\in\mathcal{C}_{0}\}$.
Using this and \eqref{eq:UL2conc} gives 
\begin{equation}\label{eq:descendantconclowerbound}
	\sum_{Q\in \mathcal{C}_{0}\setminus \mathcal{C}_{1}} 
	\iint_{Q'}|U|^2 dxdt
	+ \sum_{Q\in \mathcal{C}_1}
	\iint_{Q'}|U|^2 dxdt\geq M_{1}^{-1}.
\end{equation} 
By \eqref{eq:UtypeI} and H\"{o}lder's inequality 
in the Lorentz spaces, we have  
\begin{equation}\label{eq:UL2bound}
	\|U\|_{L^{2}(Q')}^{2} \leq C(n)  M^{p-1}(\hat{\delta} l)^n
\end{equation}
for every $Q\in\mathcal{C}_{0}$. 
Now we are going to show that if $\mathcal{C}_1\neq \emptyset$, then 
\begin{equation}\label{eq:cubecontraassumption}
	\sum_{Q\in \mathcal{C}_1}
	\iint_{Q'}|U|^2 dxdt\geq \frac{1}{2} M_{1}^{-1}
\end{equation}
cannot occur.  Assume for contradiction that \eqref{eq:cubecontraassumption} holds.  
We define
\[
	\mathcal{C}_{2}
	:=\{ Q\in \mathcal{C}_{1};\; 
	\|U\|_{L^{2}(Q')}^{2}> M_{1}^{-5}\hat{\delta}^{2+n}l^{n+1-\frac{8}{3n}} \}.
\]

By \eqref{eq:cubecontraassumption},  we get
\begin{equation}\label{eq:C1descendantconclowerbound}
	\sum_{Q\in \mathcal{C}_{1}\setminus \mathcal{C}_{2}} 
	\iint_{Q'}|U|^2 dxdt
	+ \sum_{Q\in \mathcal{C}_2}
	\iint_{Q'}|U|^2 dxdt
	\geq \frac{1}{2} M_{1}^{-1}.
\end{equation}
Now for $Q\in \mathcal{C}_{1}\setminus \mathcal{C}_{2}$, 
we have that 
\[
	\iint_{Q'} |U|^2 dxdt 
	\leq M_{1}^{-5} \hat{\delta}^{2+n} l^{n+1-\frac{8}{3n}}.  
\]
Using this, \eqref{eq:UL2bound} and \eqref{eq:C1descendantconclowerbound}, we get that
\begin{equation}\label{eq:C2C1sizes}
\begin{aligned}
	\frac{1}{2} M_{1}^{-1}
	&\leq \sum_{Q\in \mathcal{C}_{1}\setminus \mathcal{C}_{2}} 
	\iint_{Q'}|U|^2 dxdt
	+ \sum_{Q\in \mathcal{C}_2}\iint_{Q'}|U|^2 dxdt\\
	&\leq 
	(\#\mathcal{C}_{1}\setminus\mathcal{C}_{2})M_{1}^{-5} 
	\hat{\delta}^{2+n} l^{n+1-\frac{8}{3n}} 
	+ C(n) (\#\mathcal{C}_{2}) M^{p-1}(\hat{\delta} l)^n.
\end{aligned}
\end{equation}
Using \eqref{eq:sizeC1upperbound}
gives that, for appropriately chosen $M_{2}=M^{c_2(p,n)}$, 
\begin{equation*}
\begin{aligned}
	&(\#\mathcal{C}_{1}\setminus\mathcal{C}_{2}) 
	M_{1}^{-5} \hat{\delta}^{2+n} l^{n+1-\frac{8}{3n}}
	\leq M_{1}^{\frac{n}{4}-4} 
	 M^{\frac{3n}{8}(p-1)} 
	l^{1-\frac{8}{3n}} \\
	&=M_{1}^{\frac{n}{4}-4} 
	 M^{\frac{3n}{8}(p-1)} 
	M_{2}^{-(6n+2)(1-\frac{8}{3n})}
	\leq \frac{M_{1}^{-1}}{4}.
\end{aligned}
\end{equation*}
Substituting this into \eqref{eq:C2C1sizes} then gives that
\begin{equation}\label{eq:sizeC2lowerbound}
	\#\mathcal{C}_{2}\geq \frac{M_{1}^{-1}M^{1-p}}{4C(n) (\hat{\delta} l)^n}.
\end{equation}

Next, note that from the definition of $\mathcal{C}_{2}$ 
and H\"{o}lder's inequality, 
we have that for $Q\in \mathcal{C}_{2}$, 
\begin{equation}\label{eq:U3n4lowerC2}
	\iint_{Q'} |U|^{\frac{3n}{4}} dxdt
	\geq C_{n}^{*}\hat{\delta}^{n+2} l^{\frac{8+5n}{8}} M_{1}^{-\frac{15n}{8}}. 
\end{equation}
Noting that 
$\bigcup_{Q\in \mathcal{C}_{2}} Q'\subset B(3M_{1}) \times 3I$ 
and that each descendant $Q'$ from $\mathcal{C}_{2}$ 
intersects with at most $C_{n}^{**}$ descendants, we get
\begin{equation}\label{eq:U3n4upperC2}
	\sum_{Q\in\mathcal{C}_{2}} \iint_{Q'} |U|^{\frac{3n}{4}} dxdt
	\leq C_{n}^{**} \iint_{B(3M_{1}) \times 3I} 
	|U|^{\frac{3n}{4}} dxdt
	\leq C_{n}^{***} M^{\frac{3n(p-1)}{8}} M_{1}^{\frac{n}{4}}, 
\end{equation}
where the last inequality follows from \eqref{eq:UtypeI} 
and H\"{o}lder's inequality for Lorentz spaces. 
Combining \eqref{eq:sizeC2lowerbound}, \eqref{eq:U3n4lowerC2} 
and \eqref{eq:U3n4upperC2}  gives
\begin{equation*}
\begin{aligned}
	&\sum_{Q\in\mathcal{C}_{2}} \iint_{Q'} |U|^{\frac{3n}{4}} dxdt 
	\geq 
	(\#\mathcal{C}_{2})C_{n}^{*}\hat{\delta}^{n+2} 
	l^{\frac{8+5n}{8}} M_{1}^{-\frac{15n}{8}}
	\geq 
	\frac{C_{n}^{*}\hat{\delta}^{n+2} 
	l^{\frac{8+5n}{8}} M_{1}^{-\frac{15n}{8}-1}M^{1-p}}{4C(n)(\hat{\delta} l)^n}, 
	\\
	&\sum_{Q\in\mathcal{C}_{2}} 
	\iint_{Q'} |U|^{\frac{3n}{4}} dxdt
	\leq C_{n}^{***} M^{\frac{3n(p-1)}{8}} M_{1}^{\frac{n}{4}}.
\end{aligned}
\end{equation*}
This, together with \eqref{eq:deltaldef}, implies that
\begin{equation}\label{eq:Meqn}
	C_{n}^{*} (4C(n))^{-1} M_{1}^{-\frac{15n}{8}-1}
	M^{1-p} M_{2}^{\frac{(3n-8)(6n+2)}{8}-2} 
	\leq C_{n}^{***} M^{\frac{3n(p-1)}{8}} M_{1}^{\frac{n}{4}}.
\end{equation}
Noting that $((3n-8)(6n+2)/8)-2\geq 1/2$ for $n\geq 3$ 
and that $M_{2}\geq 1$, we see that \eqref{eq:Meqn} implies that
\begin{equation}\label{eq:Mcontradiction}
	C_{n}^{*} (4C(n))^{-1} 
	M_{1}^{-\frac{15n}{8}-1}M^{1-p}M_{2}^{\frac{1}{2}}
	\leq C_{n}^{***} M^{\frac{3n(p-1)}{8}} M_{1}^{\frac{n}{4}}.
\end{equation}
Thus, for $M_{2}=M^{c_{2}(p,n)}$ chosen appropriately, 
\eqref{eq:Mcontradiction} gives a contradiction. 
Hence \eqref{eq:cubecontraassumption} cannot occur.

\noindent\textbf{Step 4: conclusion}

As \eqref{eq:cubecontraassumption} cannot occur, 
by \eqref{eq:descendantconclowerbound}, we must have 
\begin{equation}\label{eq:lowerboundC0minusC1}
	\sum_{Q\in \mathcal{C}_{0}\setminus\mathcal{C}_1}
	\iint_{Q'}|U|^2 dxdt\geq \frac{1}{2} M_{1}^{-1}. 
\end{equation}
By \eqref{eq:sizecubesC0}, we have that
\[
	\#\mathcal{C}_{0}\leq \frac{C_{n} M_{1}^n}{(\hat{\delta} l)^{n+2}}. 
\]
From this, \eqref{eq:lowerboundC0minusC1}, the pigeonhole principle 
and the definition of $U$ in \eqref{eq:Udef}, 
we see that there exists $Q\in \mathcal{C}_{0}\setminus\mathcal{C}_1$ such that
\[
	\iint_{Q'} |u_{\lambda}|^{p-1} dxdt
	\geq \frac{M_{1}^{-n-1}(\hat{\delta} l)^{n+2}}{2 C_{n}}.
\]
From this, \eqref{eq:ulambdatypeI}, Lemma \ref{lem:interpolavivecor} (ii) 
and \eqref{eq:deltaldef}, it follows that
there exists $\bar{C}>0$ 
depending only on $n$ and $p$ and 
satisfying  
\begin{equation}\label{eq:ulambdaL2lower}
	\iint_{Q'} |u_{\lambda}|^2 dxdt
	\geq M^{-\bar{C}}
	\quad\mbox{ with }
	Q\in \mathcal{C}_{0}\setminus\mathcal{C}_{1}
\end{equation}
for $M$ sufficiently large. 
We also get from 
$ Q\in \mathcal{C}_{0}\setminus\mathcal{C}_{1}$ and 
\eqref{eq:ulambdabddC1} that
\begin{equation}\label{eq:ulambdabddfinal}
	\|u_{\lambda}\|_{L^{\infty}(Q(z_0, l/2))}
	\leq C l^{-\frac{2}{p-1}}M^{-2}\quad\mbox{ with }Q=Q(z_0,l).
\end{equation}
Moreover, from the definition of $\mathcal{C}_{0}$, 
we see that for $M$ sufficiently large, 
\begin{equation}\label{eq:cubeinclusionrescale}
	Q=Q(z_0,l)\subset B(5M_1/2)\times (-1/32, -1/256).
\end{equation}
Using \eqref{eq:ulambdaL2lower}, \eqref{eq:ulambdabddfinal} 
and \eqref{eq:cubeinclusionrescale}, 
we can undo the rescaling given by \eqref{eq:rescaleudefcount} 
to get the desired conclusion for the original function $u$. 
\end{proof}

\subsection{Annuli and slices of regularity}
We give annuli of quantitative regularity 
by using Proposition \ref{pro:epsreg} and arguments in \cite[Section 6]{BP21}. 
Compared to the arguments for the $3$-dimensional 
Navier-Stokes equations in \cite{BP21}, 
the lack of analogous energy structure for the energy supercritical 
nonlinear heat equation prevents us transferring 
\eqref{eq:Type1bound3.0} to a countably additive supercritical bound. 
We overcome this by applying a Calder\'{o}n type splitting \cite{Ca90}. 
This enables us to apply the $\varepsilon$-regularity criterion 
below the critical exponent (Proposition \ref{pro:epsreg}), 
together with the pigeonhole principle.

\begin{proposition}[Annuli of regularity]\label{pro:annulus}
Let $p>p_S$, $T_1>0$ and $\lambda>2$. 
Define $\beta(n,p):=(4(p+1))^{-1} \max(a_{3}(n,p), a_{4}(n,p))$ with 
\[
\begin{aligned}
	&a_{3}(n,p):=4n(6p^3+6p^2-6p-8)+4n^2(p-1), \\
	&a_{4}(n,p):=n^{2}(p-1)+n(6p^3+6p^2-7p-9)+2(6p^3+18p^2+19p+7).
\end{aligned}
\]
Suppose that $u$ is a classical solution to 
\eqref{eq:fujitaeq} on $\R^n\times [-T_{1},0]$ with
\begin{equation}\label{eq:Type1bound3.0}
	\sup_{-T_{1}<t<0}\|u(\cdot,t)\|_{L^{q_{c},\infty}(\R^n)}\leq M
\end{equation} 
and $M$ sufficiently large. 
Then, for any $R\geq 2M$, there exists
\[
	\hat{R}\in [R,R^{\lambda^{M^{q_{c}+1+\beta(n,p)}}}]
\]
such that 
\[
\begin{aligned}
	&\|u\|_{L^{\infty}\left( 
	\left\{ x\in\R^n ;\; 2T_1^{\frac{1}{2}}\hat{R}<|x|
	<\tfrac{1}{2}T_1^{\frac{1}{2}}(\hat{R})^{\lambda}  \right\}
	\times \left( -\tfrac{1}{2}T_1,0 \right)
	\right)}
	\leq M^{-1} T_{1}^{-\frac{1}{p-1}}, \\
	&\|\nabla u\|_{L^{\infty}\left( 
	\left\{ x\in\R^n;\; 2T_1^{\frac{1}{2}}\hat{R}<|x|
	<\tfrac{1}{2}T_1^{\frac{1}{2}}(\hat{R})^{\lambda}  \right\}
	\times \left( -\tfrac{1}{2}T_1,0 \right)
	\right)}
	\leq M^{-1}T_{1}^{-\frac{p+1}{2(p-1)}}.
\end{aligned}
\]
\end{proposition}

\begin{proof}
By means of the rescaling 
\[
	u_{T_1}(x,t):=T_{1}^{\frac{1}{p-1}} u(T_1^{\frac{1}{2}} x, T_1 t), 
\]
we can assume without loss of generality that $T_1=1$.

Let us decompose
\begin{equation}\label{eq:udecomp}
	u=u_{-}+u_{+}\quad\mbox{ with }
	u_{-}(x,t):=\chi_{ \{ (x,t);\; |u(x,t)|\leq 1\} }(x,t) u(x,t).
\end{equation}
Then by $(q_c+p+1)/2<q_c<2q_c$ and 
\cite[p.22]{Graf08}, we get that  
\[
	\|u_{-}(\cdot,t)\|_{L^{2q_{c}}(\R^n)}^{2q_{c}}
	+\|u_{+}(\cdot,t)\|_{L^{\frac{q_{c}+p+1}{2}}(\R^n)}^{\frac{q_{c}+p+1}{2}}
	\leq C M^{q_{c}}
\]
for $t\in [-1,0]$. 
Thus,
\[
	\sum_{k=0}^{\infty}\int_{-1}^{0} \int_{R^{\lambda^k}<|x|<R^{\lambda^{k+1}}} 
	|u_{-}(x,t)|^{2q_{c}}+|u_{+}(x,t)|^{\frac{q_{c}+p+1}{2}} 
	dxdt\leq C M^{q_{c}}.
\]
By the pigeonhole principle, there exists 
$k_0\in\{0,1\ldots, \lceil M^{q_{c}+1+\beta(n,p)} \rceil\}$ 
such that for $M$ sufficiently large, 
\begin{equation}\label{eq:uannuluspigeon}
	\int_{-1}^{0}\int_{R^{\lambda^{k_0}}<|x|<R^{\lambda^{k_0+1}}} 
	|u_{-}(x,t)|^{2q_{c}}+|u_{+}(x,t)|^{\frac{q_{c}+p+1}{2}} dxdt
	\leq M^{-\beta(n,p)}.
\end{equation}
Here $\lceil c \rceil$ is the least integer with $\lceil c \rceil \geq c$. 
Now fix 
\[
	y_0\in \{2R^{\lambda^{k_0}}<|x|<(1/2) R^{\lambda^{k_0+1}} \}, 
	\quad 
	s_0\in (-1/2,0). 
\]
As $M$ is sufficiently large, $\lambda>2$ and $R\geq 2M$, we see that
\[
	B(y_0,M) \subset 
	\{ R^{\lambda^{k_0}}<|x|<R^{\lambda^{k_0+1}}\}
	\quad \mbox{ for }
	y_0\in\{2R^{\lambda^{k_0}}<|x|<(1/2) R^{\lambda^{k_0+1}}\}. 
\]
Using this and \eqref{eq:uannuluspigeon}, we get that
\[
	\int_{-1}^{0} \int_{B(y_0,M)}|u_{-}(x,t)|^{2q_{c}}
	+|u_{+}(x,t)|^{\frac{q_{c}+p+1}{2}} dxdt\leq M^{-\beta(n,p)}. 
\]
Applying H\"{o}lder's inequality and using \eqref{eq:udecomp} gives that
\[
	\int_{-1}^{0} \int_{B(y_0,M)}  |u|^{p+1} dxdt
	\leq 
	C \left(
	M^{n-(1+\frac{\beta(n,p)}{n})\frac{p+1}{p-1}}
	+M^{\frac{-4\beta(n,p) (p+1)+n(n(p-1)-2(p+1))}{n(p-1)+2(p+1)}} \right).
\]

Define the rescaled function $\bar{u}:\R^n\times [-1,0]\rightarrow \R$ by
\begin{equation}\label{eq:urescaledannulus}
	\bar{u}(x,t):=2^{-\frac{1}{p-1}}u(2^{-\frac{1}{2}}x+y_0, 2^{-1}{t}+s_0).
\end{equation}
Then $\bar{u}$ satisfies \eqref{eq:Type1bound2.0} and 
\[
\begin{aligned}
	&{5^{n-\frac{4}{p-1}}} 
	\int_{-1/25}^{-1/50} \int_{B(M)} |\bar{u}|^{p+1} dxdt \\
	&\leq C \left( M^{n-(1+\frac{\beta(n,p)}{n}) \frac{p+1}{p-1}}
	+M^{\frac{-4\beta(n,p)(p+1)+n(n(p-1)-2(p+1))}{n(p-1)+2(p+1)}} \right).
\end{aligned}
\]
Using this and the definition of $\beta(n,p)$ 
in Proposition \ref{pro:annulus}, 
we get that for $M$ sufficiently large
\[
	5^{n-\frac{4}{p-1}} \int_{-1/25}^{-1/50} \int_{B(M)} 
	|\bar{u}|^{p+1} dxdt
	\leq M^{-6(p+1)^2}.
\]
For $M$ sufficiently large this enables us to apply Proposition \ref{pro:epsreg} 
with $(x_0,t_0)=(0,0)$, $\delta=1/5$, 
$\eps= M^{-6(p+1)^2}$ and $A=2(24\log (M^{p+1}/\eps))^{1/2}$. 
This then gives that for $M$ sufficiently large, 
\[
	|\bar{u}(0,0)|\leq CM^{-2}, \quad 
	|\nabla \bar{u}(0,0)|\leq CM^{-2}.
\]
Bearing in mind the definition of $\bar{u}$ in \eqref{eq:urescaledannulus}, this gives for $M$ sufficiently large, 
\[
	|{u}(y_0,s_0)|\leq M^{-1}, \quad 
	|\nabla {u}(y_0,s_0)|\leq M^{-1}.
\]
As $y_0$ and $s_0$ were taken arbitrarily 
in $\{2R^{\lambda^{k_0}}<|x|<(1/2) R^{\lambda^{k_0+1}} \} \times (-1/2,0)$, 
we then obtain the desired conclusion with $\hat{R}:=R^{\lambda^{k_0}}$. 
\end{proof}

We give a slice of regularity 
based on Corollary \ref{cor:corepsgeneral}, 
which is inspired by \cite[Proposition 3.5]{Pa22}. 
Compared to \cite{Pa22}, 
we cannot use countable additivity of the scale-invariant bound 
to obtain a quantity to which we can obtain a slice of regularity 
via the pigeonhole principle and $\varepsilon$-regularity. 
As was the case in Proposition \ref{pro:annulus}, 
we overcome this by applying a Calder\'{o}n type splitting \cite{Ca90}.

\begin{proposition}[Slice of regularity]\label{pro:slices}
Let $p>p_S$, $T_1>0$, $z_0=(x_0,t_0)\in \R^n\times [-T_1/2,0]$ 
and $R\leq (T_1/4)^{1/2}$.
Define $\gamma:=\max(2,a_5(n,p), a_{6}(n,p)) $ with 
\[
\begin{aligned}
	&a_5(n,p):=\frac{p-1}{p+3}\left( 
	q_{c}+n+7+12p+6p^2+\frac{2(p+1)^3-4}{p-1}\right), \\
	&a_{6}(n,p):=\frac{\frac{(p+1)q_c}{n(p-1)-(p+1)} 
	+n+7+12p+6p^2+\frac{2(p+1)^3-4}{p-1}}{\frac{2(p+1)}{p-1}-\frac{p+1}{n(p-1)-(p+1)}}.
\end{aligned}
\]
Suppose that $u$ is a classical solution to 
\eqref{eq:fujitaeq} on $\R^n\times [-T_{1},0]$ 
satisfying \eqref{eq:Type1bound3.0} and with $M$ sufficiently large. 
Then there exist a direction $\theta\in {\bf{S}}^{n-1}$ 
and a time interval $I\subset [t_0-R^2, t_0]$ 
with $|I|=(R/M^\gamma)^2$ such that within the slice 
\[
\begin{aligned}
	S&=
	\left\{ x\in\R^n;\; 
	\begin{aligned}
	&\dist(x,x_0+\R_{+}\theta)\leq 10M^{-\gamma} |(x-x_0)\cdot\theta|,  \\
	&|x-x_0|\geq 20 R 
	\end{aligned}
	\right\} \times I\\
	&\subset \R^n\times [-T_{1},0], 
\end{aligned}
\]
the following estimates holds: 
\begin{equation}\label{eq:uboundslice}
	\|u\|_{L^{\infty}(S)}
	\leq M^{-1} \left( \frac{R}{M^{\gamma}} \right)^{-\frac{2}{p-1}}, 
	\quad 
	\|\nabla u\|_{L^{\infty}(S)}
	\leq M^{-1}\left( \frac{R}{M^{\gamma}} \right)^{-\frac{p+1}{p-1}}.
\end{equation}
\end{proposition}
\begin{proof}
By means of the rescaling 
\[
	u_{R}(y,t)=R^{\frac{2}{p-1}}u(Ry+x_0, R^2 t+t_0), 
\]
we can take $z_0=(x_0,t_0)=(0,0)$ and $R=1$ without loss of generality.

Let us decompose
\begin{equation}\label{eq:udecompslices}
	u=u_{-}+u_{+}\quad \mbox{ with } 
	u_{-}(x,t):=\chi_{ \{(x,t); |u(x,t)|\leq 1\} }(x,t) u(x,t).
\end{equation}
Then by $p+1<q_c<n(p-1)-(p+1)$ and  \cite[p.22]{Graf08}, 
we get that for $t\in [-1,0]$, 
\[
	\|u_{-}(\cdot,t)\|_{L^{n(p-1)-(p+1)}(\R^n)}^{n(p-1)-(p+1)}
	+\|u_{+}(\cdot,t)\|_{L^{p+1}(\R^n)}^{p+1}\leq C M^{q_{c}}. 
\]
Thus,
\begin{equation}\label{eq:udecompintegestslices}
	\int_{-1}^{0}\int_{\R^n} 
	( |u_{-}(x,t)|^{n(p-1)-(p+1)}+|u_{+}(x,t)|^{p+1} ) dxdt
	\leq C M^{q_{c}}.
\end{equation}
Consider the collection of space-time slices of the form
\begin{equation}\label{eq:slicecollection}
\begin{aligned}
	S^*&=\{x\in\R^n;\; 
	\dist(x,\R_{+}\theta)\leq 20M^{-\gamma}|x\cdot\theta|, \; |x|\geq 10\}  \\
	&\quad \times [-10M^{-2\gamma}k,-10M^{-2\gamma}(k-1) ]
\end{aligned}
\end{equation}
with $\theta\in \textbf{S}^{n-1}$ 
and for $k\in\{1,\ldots, \lfloor (1/10)M^{2\gamma}\rfloor\}$, 
where $\lfloor c \rfloor$ is the greatest integer with $\lfloor c \rfloor \leq c$. 
We may find a collection $\mathcal{S}$ of disjoint slices of the form \eqref{eq:slicecollection} within $\R^n\times [-1,0]$ such that
\[
	\#\mathcal{S} \geq  
	\frac{M^{2\gamma} M^{\gamma (n-1)}}{C(n)}
	= \frac{M^{\gamma(n+1)}}{C(n)}.
\]
Thus \eqref{eq:udecompintegestslices} and the pigeonhole principle implies that there exists a slice $S^*\in\mathcal{S}$ of the form
\[
\begin{aligned}
	&S^*=S^*_x\times [-10M^{-2\gamma}k,-10M^{-2\gamma}(k-1) ], \\
	&S^*_{x}:=\{x\in\R^n; \; \dist(x,\R_{+}\theta)
	\leq 20M^{-\gamma}|x\cdot\theta|, \; |x|\geq 10\}
\end{aligned}
\]
such that 
\[
	\int_{-10M^{-2\gamma}k}^{-10M^{-2\gamma}(k-1)}\int_{S^{*}_{x}}
	( |u_{-}|^{n(p-1)-(p+1)}+|u_{+}|^{p+1} ) dxdt\leq C M^{q_{c}-(n+1)\gamma}. 
\]

Now fix 
\begin{equation}\label{eq:ysfixed}
\begin{aligned}
	&y\in \{x\in\R^n;\;  \dist(x,\R_{+}\theta)
	\leq 10M^{-\gamma}|x\cdot\theta|,\; |x|\geq 20\}, \\
	&s\in (-10M^{-2\gamma}(k-1)-M^{-2\gamma},-10M^{-2\gamma}(k-1)).
\end{aligned}
\end{equation}
For such $(y,s)$, we have 
\[
	Q( (y,s), M^{-\gamma})
	\subset  S^{*}_{x}
	\times (-10M^{-2\gamma}(k-1)-M^{-2\gamma},-10M^{-2\gamma}(k-1)).
\]
Thus,
\[
	\iint_{Q( (y,s), M^{-\gamma})} 
	(|u_{-}|^{n(p-1)-(p+1)}+|u_{+}|^{p+1}) dxdt
	\leq C M^{q_{c}-(n+1)\gamma}.
\]
Using H\"{o}lder's inequality, \eqref{eq:udecompslices} and the definition of $\gamma$ in the statement of the proposition gives that for $M$ sufficiently large, 
\begin{equation}\label{eq:upplus1cube} 
\begin{aligned}
	&M^{(n-\frac{4}{p-1})\gamma}
	\iint_{Q( (y,s), M^{-\gamma})} |u|^{p+1} dxdt\\
	&\leq C \left( M^{q_{c}-\frac{(p+3)\gamma}{p-1}}
	+M^{\frac{(p+1)(q_c+\gamma)}{n(p-1)-(p+1)}-\frac{2(p+1)\gamma}{p-1}} 
	\right)\\
	&\leq C M^{-(n+7+12p+6p^2+\frac{2(p+1)^3-4}{p-1})}. 
\end{aligned}
\end{equation}
Noting that, as $\gamma\geq 2$, 
taking $M^{-\gamma}<(2M)^{-1}$ 
for $M$ sufficiently large, 
Thus, for $M$ sufficiently large, 
\eqref{eq:upplus1cube} allows us to apply Corollary \ref{cor:corepsgeneral} 
with $\delta=M^{-\gamma}$ giving 
\[
	|u(y,s)|\leq M^{-1}\left( \frac{1}{M^{\gamma}}\right)^{-\frac{2}{p-1}}, 
	\quad 
	|\nabla u(y,s)|\leq M^{-1}\left( \frac{1}{M^{\gamma}} \right)^{-\frac{p+1}{p-1}}. 
\]
As $(y,s)$ are chosen to be any space-time points satisfying \eqref{eq:ysfixed}, we obtain the desired conclusion \eqref{eq:uboundslice} 
for $z_0=(x_0,t_0)=(0,0)$ and $R=1$.
\end{proof}

\section{Carleman inequalities}\label{sec:Carl}
The quantitative unique continuation Carleman inequality we use 
is a higher-dimensional analogue of that given 
by Tao \cite[Proposition 4.3]{Ta21}, in a lower regularity setting. 
The quantitative backward uniqueness Carleman inequality that we use 
is from Palasek \cite[Proposition 9]{Pa21}, in a lower regularity setting 
(and with a certain adjustment of indices). 
We require a lower regularity setting ($C^{2,1}$) 
compared to previous works, 
since bounded solutions to the nonlinear heat equation 
with non-smooth nonlinearity 
are not necessarily $C^{\infty}$, 
in contrast with the Navier-Stokes equations, for example.

The proof of these Carleman inequalities 
(in both the smooth and lower regularity settings) 
hinges on the following general inequality proven by Tao \cite[Lemma 4.1]{Ta21}.
Throughout this section, let $L$ be the backward heat operator
\[
	L:=\partial_{t}+\Delta.
\]

\begin{lemma}(General Carleman inequality, \cite[Lemma 4.1]{Ta21})\label{lem:gencarl}
Let $g:\R^n \times [t_1,t_2]\rightarrow\R $ be smooth and let $D^2 g$ be the bilinear form expressed in coordinates (with usual summation convention) as
 \begin{equation}\label{eq:bilineargdef}
 D^2 g(a,b):=(\partial_{i}\partial_{j} g) a_{ik} b_{jk}.
 \end{equation} 
 Let $F:\R^n \times [t_1,t_2]\rightarrow\R $ denote the function
\begin{equation}\label{eq:Fdef}
	F:=\partial_{t}g-\Delta g-|\nabla g|^2.
\end{equation}
Then, for any vector-valued function 
$v\in C^{\infty}_{0}(\R^n\times [t_1,t_2]; \R^m)$, 
\[
\begin{aligned}
	&\frac{d}{dt}\int_{\R^n}\left( 
	|\nabla v|^2+\frac{1}{2} F|v|^2 \right) 
	e^{g} dx \\
	&\geq \int_{\R^n} \left( 
	\frac{1}{2}(LF)|v|^2+2D^2 g(\nabla v,\nabla v)-\frac{1}{2}|L v|^2
	\right) e^g dx
\end{aligned}
\]
for all $t\in [t_1,t_2]$. 
Moreover, $v$ also satisfies 
\[
\begin{aligned}
	&\int_{t_1}^{t_2} \int_{\R^n} 
	\left( \frac{1}{2}(LF)|v|^2+2D^2 g(\nabla v,\nabla v) \right) e^g dxdt \\
	&\leq \frac{1}{2}\int_{t_1}^{t_2} \int_{\R^n} 
	|L v|^2 e^g dxdt
	+ \left. \int_{\R^n}\left( 
	|\nabla v|^2+\frac{1}{2} F|v|^2 \right)
	e^g dx \right|^{t=t_2}_{t=t_1}.
\end{aligned}
\]
\end{lemma}
We state the quantitative backward uniqueness 
and unique continuation Carleman inequalities, 
and then we show an iterated quantitative unique continuation.

\begin{proposition}[Backward uniqueness Carleman inequality]\label{pro:backuniqueness}
Let $0<r_-<r_+<\infty$ and $0<T_1<\infty$. 
Define a space-time annulus $\mathcal{A}$ by 
\begin{equation*}
	\mathcal{A}
	:=\{(x,t)\in \R^n\times \R; \; r_-\leq |x|\leq r_+, \; t\in[-T_1,0] \}.
\end{equation*}
Let $w\in C^{2,1}(\mathcal A)$ satisfy the differential inequality 
\begin{equation}\label{e.diffineq}
	|\partial_t w - \Delta w| 
	\leq \frac{|w(x,t)|}{C_{\Carl} T_1}
	+\frac{|\nabla w (x,t)|}{(C_{\Carl} T_1)^{1/2}}
	\quad\mbox{ on } \mathcal{A}
\end{equation}
with a constant $C_{\Carl}\geq 3n$. 
Assume
\begin{equation}\label{e.lowerr-}
r_-^2\geq 4C_{\Carl}T_1.
\end{equation}
Then there exists $C(n)>0$ depending only on $n$ such that 
\[
\begin{aligned}
	&\int_{-T_1/4}^{0}\int_{10r_-\leq|x|\leq r_+/2}
	\left( \frac{|w(x,t)|^2}{T_1} +|\nabla w(x,t)|^2 \right) dxdt \\
	&\leq C(n) C_{\Carl}e^{-\frac{r_-\cdot r_+}{4C_{\Carl}T_1}} X
	+ C(n) C_{\Carl}
	e^{-\frac{r_-\cdot r_+}{4C_{\Carl}T_1} + \frac{2r_+^2}{C_{\Carl}T_1}} Y,
\end{aligned}
\]
where 
\[
\begin{aligned}
	&X:=\iint_{\mathcal A}e^{\frac{2|x|^2}{C_{\Carl}T_1}}
	\left( \frac{|w(x,t)|^2}{T_1} +|\nabla w(x,t)|^2 \right) dxdt, \\
	&Y:=\int_{r_-\leq |x|\leq r_+}|w(x,0)|^2 dx.
\end{aligned}
\]
\end{proposition}

\begin{proof}
Let us briefly outline 
how the corresponding proof in \cite[Proposition 9]{Pa21} also applies 
in our less regular setting.
Without loss of generality, suppose $20r_{-}\leq r_{+}$.
By the change of variables, we have 
\[
	X=\int_{0}^{T_1}\int_{r_{-}\leq |x|\leq r_{+}}e^{\frac{2|x|^2}{C_{\Carl}T_1}}
	\left( \frac{|w(x,-t)|^2}{T_1} +|\nabla w(x,-t)|^2 \right) dxdt. 
\] 
Then by the pigeonhole principle, 
there exists $T_0\in [T_1/2, 3T_1/4]$ such that
\begin{equation}\label{eq:pigeonbackuniqueCarl}
	\int_{r_{-}\leq |x|\leq r_{+}}e^{\frac{2|x|^2}{C_{\Carl}T_1}}
	\left( \frac{|w(x,-T_0)|^2}{T_1} +|\nabla w(x,-T_0)|^2 \right) dx\leq 4 T_{1}^{-1} X.
\end{equation}
Let 
\begin{equation}\label{eq:gdefbackwarduniqueness}
	g(x,t):=\frac{r_{+}(T_0-t)|x|}{2C_{\Carl}{(T_1)}^{2}}+\frac{|x|^2}{C_{\Carl}T_1}.
\end{equation}

By using the definitions \eqref{eq:bilineargdef} and \eqref{eq:Fdef}, 
together with $r_{-}^{2}\geq 4 C_{\Carl} T_{1}$ and $C_{\Carl}\geq 3n$, 
direct computations readily show that 
for $(x,t)\in \{r_{-}\leq |x|\leq r_{+}\}\times [0, T_0]$, 
\begin{equation}\label{eq:FGsingedbackCarl}
	F(x,t)\leq 0, 
	\quad 
	D^2 g(x,t)\geq \frac{2}{C_{\Carl} T_{1}}
	Id, \quad 
	LF(x,t)\geq \frac{4}{C_{\Carl} (T_{1})^2}, 
\end{equation}
where $Id$ is the identity matrix. 
Let $\psi$ be a smooth spatial cut-off 
compactly supported in $\{r_{-}\leq |x|\leq r_{+}\}$ 
that equals $1$ in $\{2r_{-}\leq |x|\leq r_{+}/2\}$ and obeys 
\begin{equation}\label{eq:psiderest}
	|\nabla^{j} \psi(x)|\leq C(n,j) r_{-}^{-j} \quad \mbox{ for }j=0,1,2.
\end{equation}
The main difference in the less regular setting compared 
with \cite[Proposition 9]{Pa21}, 
is that we do not have sufficient regularity 
to directly apply Lemma \ref{lem:gencarl} (with the weight $g$ as above) 
to $\psi(x) w(x,-t)$  on $\R^n\times [0, T_0]$. 
Instead, we apply Lemma \ref{lem:gencarl} 
with the weight $g$ to $\psi(x) v_{1/M, 1/N}(x,t)$ 
on  $\R^n\times [\varepsilon, T_0]$. 
Here, 
\[
	v_{\frac{1}{M}, \frac{1}{N}}(x,t):=\int_{0}^{T_1} 
	\omega_{\R, \frac{1}{M}} (t-s)\int_{r_{-}\leq |y|\leq r_{+}} 
	\omega_{\R^n,\frac{1}{N}} (x-y) w(y,-s) dyds 
\]
and $\omega_{\R, 1/M}$ (resp. $\omega_{\R^n, 1/N}$) 
denotes a standard mollifier on $\R$ (resp. $\R^n$) 
supported in $[-1/M, 1/M]$ (resp. $\{|x|\leq 1/N\}$).

Applying Lemma \ref{lem:gencarl}, 
together with using \eqref{eq:FGsingedbackCarl}, gives
\[
\begin{aligned}
    & \int_{\varepsilon}^{T_0} \int_{\R^n} \left( 
    \frac{2}{C_{\Carl} (T_{1})^2}|\psi v_{\frac{1}{M},\frac{1}{N}}|^2
    +\frac{4}{C_{\Carl} T_{1}} |\nabla(\psi v_{\frac{1}{M},\frac{1}{N}})|^2 
    \right) e^g dxdt \\
    &\leq \frac{1}{2}\int_{\varepsilon}^{T_0}
    \int_{\R^n} |L(\psi v_{\frac{1}{M},\frac{1}{N}})|^2 e^g dxdt
    +\int_{\R^n}|\nabla(\psi v_{\frac{1}{M},\frac{1}{N}})(x, T_0)|^2 e^{g(x,T_0)} dx\\
    &\quad 
    + \frac{1}{2}\int_{\R^n}|F(x,\varepsilon)||\psi v_{\frac{1}{M},\frac{1}{N}}(x, \varepsilon)|^2 e^{g(x,\varepsilon)} dx.
\end{aligned}
\]
Now we take the limit $N\rightarrow \infty$, 
then $M\rightarrow \infty$ and finally $\varepsilon\rightarrow 0$.
Using $w\in C^{2,1}(\mathcal{A})$, $T_0\in[T_1/2, 3T_1/4]$  
and standard properties of mollification yield that 
\[
\begin{aligned}
    & \int_{0}^{T_0}\int_{\R^n} 
    \left( \frac{2}{C_{\Carl} (T_{1})^2}|\psi w(x,-t)|^2+\frac{4}{C_{\Carl} T_{1}} |\nabla(\psi w(x,-t))|^2 \right) e^g dxdt\\
    &\leq \frac{1}{2}\int_{0}^{T_0}\int_{\R^n} |L(\psi w(x, -t))|^2 e^g dxdt
    +\int_{\R^n}|\nabla(\psi w)(x, -T_0)|^2 e^{g(x,T_0)} dx\\
    &\quad + \int_{\R^n}|F(x,0)||\psi w(x, 0)|^2 e^{g(x,0)} dx.
\end{aligned}
\]
This together with \eqref{e.diffineq} (with $C_{\Carl}\geq 3n$), 
\eqref{e.lowerr-}, \eqref{eq:pigeonbackuniqueCarl}, 
\eqref{eq:gdefbackwarduniqueness} and \eqref{eq:psiderest} 
allows us to apply the same arguments 
used in \cite[Proposition 9]{Pa21} 
(particularly, \cite[pp. 1521--1523]{Pa21}), 
which then yields the desired conclusion. 
In particular, the reader can readily check that 
$w\in C^{2,1}(\mathcal{A})$ is sufficient regularity for this final step.

We remark that the setting in \cite[Proposition 9]{Pa21} 
concerns $x\in \R^{d_1+d_2}$, 
with the (anisotropic) annuli 
in the  quantitative backward uniqueness Carleman inequality 
being in $\R^{d_1}$. 
This is more general than we require 
and we can consider the reasoning in \cite[Proposition 9]{Pa21} 
with just $d_1=n$ and $d_2=0$.  
Moreover, we just require that the time derivative 
and second spatial derivatives are square integrable to perform the final step, 
along with \eqref{e.diffineq}. 
\end{proof}

\begin{proposition}[Unique continuation Carleman inequality]\label{pro:uniquecont}
Let $r>0$, $0<T_1<\infty$ and $1\leq C_{0}<\infty$. 
Define a space-time cylinder $\mathcal{C}$ by 
\begin{equation*}
\mathcal C:=\{(x,t)\in \R^n \times \R; \; |x|\leq r, \; t\in [-T_{1},0]\}.
\end{equation*}
Let $w\in C^{2,1}(\mathcal C)$ satisfy the differential inequality
\begin{equation}\label{e.diffineqC}
	|\partial_{t}w - \Delta w|
	\leq \frac{|w(x,t)|}{C_0 T_{1}}
	+\frac{|\nabla w (x,t)|}{(C_0 T_{1})^{1/2}}\quad\mbox{ on }\mathcal C.
\end{equation}
Assume 
\begin{equation}\label{e.lowerr}
	r^2\geq 16000T_{1}.
\end{equation}
Then there exists $C(n)>0$ depending only on $n$ such that 
\[
\begin{aligned}
	&\int_{-2\overline s}^{-\overline s}\int_{|x|\leq r/2}
	\left( \frac{|w(x,t)|^2}{T_1} +|\nabla w(x,t)|^2 \right) 
	e^{\frac{|x|^2}{4t}} dxdt \\
	&\leq C(n) e^{-\frac{r^2}{500\overline s}}X
	+C(n) (\overline s)^\frac{n}{2}
	\left( \frac{3e\overline s}{\underline s} \right)^{\frac{r^2}{200\overline s}}Y
\end{aligned}
\]
for all $0<\underline{s}\leq\overline s<T_1/(10000n)$,  where 
\[
\begin{aligned}
	&X:=\int_{-T_{1}}^0\int_{|x|\leq r}
	\left( \frac{|w(x,t)|^2}{T_1} + |\nabla w(x,t)|^2 \right) dxdt, \\
	&Y:=\int_{|x|\leq r}|w(x,0)|^2
	(\underline s)^{-\frac{n}{2}}e^{-\frac{|x|^2}{4\underline s}} dx.
\end{aligned}
\]
\end{proposition}

\begin{proof}
Let us briefly outline the modification required 
to the proof of \cite[Proposition 4.3]{Ta21} 
to our setting with general dimension and lower regularity.
Clearly, 
\[
	X=\int_{0}^{T_1}\int_{ |x|\leq r}
	\left( \frac{|w(x,-t)|^2}{T_1} +|\nabla w(x,-t)|^2 \right) dxdt
\]
and by the pigeonhole principle 
there exists $S_0\in [T_1/200,T_1/100]$ such that
\begin{equation}\label{eq:pigeonuniquecontCarl}
	\int_{ |x|\leq r}
	\left( \frac{|w(x,-S_0)|^2}{T_1} +|\nabla w(x,-S_0)|^2 \right) dx\leq 200 T_{1}^{-1} X.
\end{equation}

Let $\alpha:= r^2/(400\bar{s})$ and  
\[
	g(x,t):=\frac{-|x|^2}{4(t+\underline{s})}
	-\frac{n}{2}\log(t+\underline{s})
	-\alpha\log\left( \frac{t+\underline{s}}{S_0+\underline{s}}\right)
	+\alpha\frac{t+\underline{s}}{S_0+\underline{s}}, 
\]
where we note that in \cite[Proposition 4.3]{Ta21}, 
this weight is considered in the case $n=3$.
Using the definitions \eqref{eq:bilineargdef} and \eqref{eq:Fdef}, 
direct computations readily give that
\begin{equation}\label{eq:FGuniquecontCarl}
\begin{aligned}
	&F(x,t)=\frac{\alpha}{S_0+\underline{s}}-\frac{\alpha}{t+\underline{s}},
	\quad D^2 g(x,t)=-\frac{1}{2(t+\underline{s})} Id, \\
	&LF(x,t)=\frac{\alpha}{(t+\underline{s})^2}.
\end{aligned}
\end{equation}
Let $\psi\in C_{0}^{\infty}(B(r))$ 
be a smooth spatial cut-off such that 
$\psi=1$  on $B(0,r/2)$ and 
\[
	|\nabla^{j} \psi(x)|\leq C(n,j) r^{-j} \quad \mbox{ for }j=0,1,2.
\]
The main difference in the less regular setting 
compared with \cite[Proposition 4.3]{Ta21}, 
is that we do not have sufficient regularity 
to directly apply Lemma \ref{lem:gencarl} (with the weight $g$ as above) 
to $\psi(x) w(x,-t)$ on $\R^n\times [0, S_0]$. 
Instead, we apply Lemma \ref{lem:gencarl} 
with the weight $g$ to $\psi v_{1/M, 1/N}(x,t)$ on $\R^n\times [0, S_0]$. 
Here,
\[
	v_{\frac{1}{M}, \frac{1}{N}}(x,t)
	:=\int_{0}^{T_1} \omega_{\R, \frac{1}{M}} (t-s)
	\int_{|y|\leq r} \omega_{\R^n,\frac{1}{N}} (x-y) w(y,-s) dyds.
\]

Applying Lemma \ref{lem:gencarl}, together with using \eqref{eq:FGuniquecontCarl}, gives that for $t\in [0,S_0]$, 
\begin{align}
	&\label{eq:generalCarluniquecont}
	\frac{dE_{M,N}}{dt}(t) 
	\geq \int_{\R^n}\left( 
	\frac{\alpha|\psi v_{\frac{1}{M},\frac{1}{N}}|^2}{2(t+\underline{s})^2} 
	-\frac{|\nabla(\psi v_{\frac{1}{M},\frac{1}{N}})|^2}{t+\underline{s}}
	-\frac{|L(\psi v_{\frac{1}{M},\frac{1}{N}})|^2}{2}\right) e^{g} dx,  \\
	&\notag
	E_{M,N}(t):=\int_{\R^n} \left( 
	|\nabla(\psi v_{\frac{1}{M},\frac{1}{N}})|^2
	+\frac{1}{2}\left( \frac{\alpha}{S_0+\underline{s}}
	-\frac{\alpha}{t+\underline{s}} \right)| 
	\psi v_{\frac{1}{M},\frac{1}{N}}|^2 \right) e^{g} dx. 
\end{align}
Using \eqref{eq:generalCarluniquecont}, $S_0\in [T_1/200, T_1/100]$ 
and $0<\underline{s}\leq\overline s<T_1/(10000n)$, we get that for $t\in [0,S_0]$, 
\begin{equation}\label{eq:EMNweightedinequality}
\begin{aligned}
	& \frac{d}{dt}\left( \left( 
	(t+\underline{s})+\frac{(t+\underline{s})^2}{10 S_0}\right) 
	E_{M,N}(t)\right) \\
	&\geq \int_{\R^n} \left( 
	\frac{\alpha |\psi v_{\frac{1}{M},\frac{1}{N}}|^2}{10(S_0+\underline{s})} 
	+ \frac{(t+\underline{s}) |\nabla(\psi v_{\frac{1}{M},\frac{1}{N}})|^2}
	{10(S_0+\underline{s})}
	-(t+\underline{s})|L(\psi v_{\frac{1}{M},\frac{1}{N}})|^2\right) e^{g} dx.
\end{aligned}
\end{equation}
Now take $\varepsilon\in (0,S_0]$. 
Using this, 
\eqref{eq:EMNweightedinequality} and the fundamental theorem of calculus, one can deduce that
\[
\begin{aligned}
	&\int_{\varepsilon}^{S_0} \int_{\R^n}
	\left( \frac{\alpha|\psi v_{\frac{1}{M},\frac{1}{N}}|^2}{10(S_0+\underline{s})} 
	+\frac{(t+\underline{s}) |\nabla(\psi v_{\frac{1}{M},\frac{1}{N}})|^2 }
	{10(S_0+\underline{s})}\right) e^{g} dxdt\\
	&\leq \int_{\varepsilon}^{S_0}\int_{\R^n}
	(t+\underline{s})|L(\psi v_{\frac{1}{M},\frac{1}{N}})|^2 e^{g} dx dt
	+\left. \left( \left( (t+\underline{s})+\frac{(t+\underline{s})^2}{10 S_0}
	\right)E_{M,N}(t)\right) \right|^{t=S_0}_{t=\varepsilon}\\
	&\leq \int_{\varepsilon}^{S_0}\int_{\R^n} 
	(t+\underline{s})|L(\psi v_{\frac{1}{M},\frac{1}{N}})|^2 e^{g} dx
	+C(n)S_0\int_{\R^n} |\nabla(\psi v_{\frac{1}{M},\frac{1}{N}})(x,S_0)|^2 
	e^{g(x,S_0)}dx\\
	&\quad 
	+C(n)\alpha \int_{\R^n} |\psi v_{\frac{1}{M},\frac{1}{N}}(x,\varepsilon)|^2 
	e^{g(x,\varepsilon)}dx.
\end{aligned}
\]

We now take the limit of the above inequality 
as $N\rightarrow \infty$, $M\rightarrow \infty$ and finally $\varepsilon\rightarrow 0$. 
Using that $w\in C^{2,1}(\mathcal{C})$, $S_0\in[T_1/200, T_{1}/100]$ 
and standard properties of mollification yields that 
\[
\begin{aligned}
	& \int_{0}^{S_0}\int_{\R^n}\left( 
	\frac{\alpha}{10(S_0+\underline{s})} |\psi w(x,-t)|^2 
	+\frac{t+\underline{s}}{10(S_0+\underline{s})}|\nabla(\psi w(x,-t))|^2 
	\right) e^{g} dxdt\\
	&\leq \int_{0}^{S_0}\int_{\R^n}(t+\underline{s})|L(\psi w(x,-t))|^2 
	e^{g} dx
	+C(n)S_0\int_{\R^n} |\nabla(\psi w)(x,-S_0)|^2 e^{g(x,S_0)}dx\\
	&\quad 
	+ C(n)\alpha \int_{\R^n} |\psi w(x,0)|^2 e^{g(x,0)}dx.
\end{aligned}
\]
This together with \eqref{e.diffineqC} (with $C_{0}\geq 1$), 
\eqref{e.lowerr}, $0<\underline{s}\leq\overline s<T_1/(10000n)$ 
and \eqref{eq:pigeonuniquecontCarl} allows us to apply 
similar arguments used in \cite[Proposition 4.3]{Ta21} 
(particularly, from \cite[(4.17)]{Ta21}), 
albeit in the setting with general dimension. 
Then we obtain the desired conclusion. 
In particular, 
the reader can readily check that $w\in C^{2,1}(\mathcal{C})$ 
is sufficient regularity for this final step.
\end{proof}

We will use the terminology `$\eta$ being sufficiently small' 
to mean that $0<\eta\leq \eta_{0}(n,p)$, 
where $\eta_0$ is positive and sufficiently small. 
Let us prove the iterated quantitative unique continuation, 
which is based upon \cite[Proposition 4.2]{Pa22} that corresponds to the case $p=2$.

\begin{proposition}[Iterated unique continuation]
\label{pro:iterateduniquecontinuation}
Let $p>1$, $0<\eta\leq 1$, $0<T_1<\infty$, $1\leq C_0<\infty$ 
and $\theta\in \bS^{n-1}$. 
Define 
\begin{equation}\label{eq:sliceuniquecont}
	S:=\{x\in \R^n: |x|>10T_{1}^{\frac{1}{2}}, \; \dist(x,\R_{+}\theta)
	\leq \eta |x\cdot\theta|\} 
	\times [-\eta T_1,0]. 
\end{equation}
Let $w\in C^{2,1}(S)$ satisfy 
\begin{equation}\label{eq:diffinequality}
\begin{aligned}
	&\|\nabla^j w\|_{L^{\infty}(S)}
	\leq ( \eta T_1)^{-\tfrac{1}{p-1}-\tfrac{j}{2}}\quad \mbox{ for }j=0,1,  \\
	&|\partial_{t} w -\Delta w|
	\leq \frac{|w(x,t)|}{C_0\eta T_1}
	+\frac{|\nabla w (x,t)|}{(C_0\eta T_1)^{1/2}}\quad \mbox{ on }S.
\end{aligned}
\end{equation}
Assume that for every $t\in [-\eta T_1,0]$,  
\[
	\int_{B(R_0\theta, \eta^5 R_0)} |w(x,t)|^2 dx
	\geq \epsilon T_{1}^{\tfrac{n}{2}-\tfrac{2}{p-1}}, 
\]
where  $20T_{1}^{1/2}\leq R_0\leq \eta^{-2} T_1^{1/2}$ 
and $0<\epsilon\leq \eta^8$. 
Then, for $\eta$ sufficiently small, 
\begin{equation}\label{eq:witeratedconc}
	\int_{B(R\theta, \eta^2 R)} |w(x,t)|^2 dx
	\geq  \epsilon^{\left( \frac{R}{R_0} \right)^{\eta^{-4}}} 
	T_{1}^{\tfrac{n}{2}-\tfrac{2}{p-1}}
\end{equation}
for all $t\in [-\eta T_1/2,0]$ and $R\geq R_0$. 
\end{proposition}

To prove Proposition \ref{pro:iterateduniquecontinuation}, 
we prepare the following lemma:

\begin{lemma}\label{lem:uniquecontslice}
Let $p>1$, $0<\eta\leq 1$, $0<T_1<\infty$, 
$1\leq C_0 <\infty$, $\theta\in \bS^{n-1}$ and $a\in (\eta/2, 3\eta/4]$.
Define $S$ by \eqref{eq:sliceuniquecont} and let 
$w\in C^{2,1}(S)$ satisfy \eqref{eq:diffinequality}. 
Assume that  
\begin{equation}\label{eq:winitialconclem}
\begin{aligned}
	&\int_{B(R\theta, \eta^5 R)} 
	|w(x,t)|^2 dx\geq \epsilon_{0} T_{1}^{\tfrac{n}{2}-\tfrac{2}{p-1}} 
	\quad \mbox{ for every } t\in [-aT_{1},0], \\
	&\mbox{where }R\geq 20T_{1}^{\frac{1}{2}} \mbox{ and } 0<\epsilon_{0}\leq 
	\min\left( \eta^8,\left( \frac{R^2}{T_1}\right)^{-5000n\eta},
	e^{-\frac{40000n\eta^4 R^2}{T_1}}\right). 
\end{aligned}
\end{equation}
Then, for $\eta$ sufficiently small and $R'=(1+\eta^3)R$, 
\[
	\int_{B(R'\theta,\eta^5 R')} |w(x,t)|^2 dx
	\geq \epsilon_{0}^{{\eta^{-2}}} T_{1}^{\tfrac{n}{2}-\tfrac{2}{p-1}}
\]
for all $t\in [-aT_{1}+2\eta^5 R^2 (\log(1/\epsilon_0))^{-1},0]$. 
\end{lemma}

\begin{proof}
Lemma \ref{lem:uniquecontslice} is based on \cite[Lemma 4.3]{Pa22}. 
However, we required minor differences on the condition 
for $a$ and $\epsilon_{0}$ for the statement and proof to hold. 
For the readers' convenience, 
we include the proof (following \cite[Lemma 4.3]{Pa22}) with these adjustments.

First note that by means of the rescaling 
\[
	w_{T_{1}}(x,t):=T_1^{\frac{1}{p-1}} w(T_1^{\frac{1}{2}} x, T_{1} t), 
\]
we can take $T_{1}=1$ without loss of generality. 
Let $R$ be given in \eqref{eq:winitialconclem}. 
Note also that for $0<\eta\leq 1$, 
the fact that 
$0<\epsilon_{0}\leq \min(\eta^8,e^{{-40000n\eta^4 R^2}} )$ 
and $a\in(\eta/2, 3\eta/4]$  ensures that 
\[
	-\frac{3}{4}\eta \leq -a 
	\leq -a + \frac{2\eta^5 R^2}{ \log (1 / \epsilon_0)}<0. 
\]

Now fix $t'\in  [-a + 2\eta^5 R^2 (\log(1/\epsilon_0))^{-1} ,0]$ and define 
\[
\begin{aligned}
	r:=\eta^2 R, \quad T_{c}:=\min(\eta/4, \eta^5 R^2), \quad 
	t_0:=\frac{\eta^5 R^2}{\log (1/\epsilon_0)}, \quad 
	t_1:=\frac{\eta^{15} R^2}{\log(1/\epsilon_0)}. 
\end{aligned}
\]
As $0<\epsilon_{0}\leq \eta^8$, 
$0<a\leq 3\eta/4$ and $R'=(1+\eta^3)R$, for $\eta$ sufficiently small, 
\begin{align}
	&\label{eq:domaininclusioninS}
	B(R'\theta, r) \times [t'-T_{c}, t']\subset S, \\
	&\label{eq:timeintervalinclusion}
	\begin{aligned}
	&B(R'\theta,r/2)\supset B(R\theta,\eta^5 R), 
	\quad [t'-2t_0, t'-t_0]\subset [-a,0], \\
	&|x-R'\theta|\leq 2\eta^3 R
	\quad \mbox{ for }x\in B(R\theta,\eta^5 R). 
	\end{aligned}
\end{align}
As $0<\epsilon_{0}\leq \min(\eta^8,e^{{-40000n\eta^4 R^2}})$, 
we also get that for $\eta$ sufficiently small, 
\begin{equation}\label{eq:t0t1rCarl}
	r^2\geq 16000T_{c}, \quad 0<t_{1}<t_0\leq \frac{T_{c}}{10000n}.
\end{equation}
Furthermore, from \eqref{eq:diffinequality}, 
\eqref{eq:domaininclusioninS} and $T_{c}\leq \eta/4$, 
it follows that 
\begin{equation}\label{eq:diffinequalityTc}
	|\partial_{t}w -\Delta w |
	\leq \frac{|w(x,t)|}{4 C_0 T_c} + \frac{|\nabla w(x,t)|}{(4 C_0 T_c)^{1/2}}
	\quad \mbox{ on } B(R'\theta,r) \times [t'-T_{c}, t'].
\end{equation}
Then \eqref{eq:t0t1rCarl} and \eqref{eq:diffinequalityTc}, 
along with the fact that $1\leq C_0<\infty$, 
allows us to apply Proposition \ref{pro:uniquecont} 
with $T_{c}$, $\bar{s}=t_0$ and $\underline{s}=t_{1}$ to 
\[
	(x,t)\mapsto w(x+R'\theta, t+t')
	\quad\mbox{ for }(x,t)\in B(r)\times [-T_{c}, 0].
\]
This yields 
\begin{equation}\label{e.carlslice}
	Z \leq C(n) e^{-\frac{r^2}{500t_0}}X 
	+C(n) (t_0)^\frac{n}{2}
	\left( \frac{3e t_0}{t_1} \right)^{\frac{r^2}{200 t_0}}Y,
\end{equation}
where 
\[
\begin{aligned}
	&X:=\int_{t'-T_{c}}^{t'}\int_{B(R'\theta,r)}
	(T_{c}^{-1}|w|^2+|\nabla w|^2) dxdt, \\
	&Y:=\int_{B(R'\theta,r)} |w(x,t')|^2 (t_1)^{-\frac{n}{2}}
	e^{-\frac{|x-R'\theta|^2}{4 t_1}} dx, \\
	&Z:=\int_{t'-2 t_0}^{t'-t_0}\int_{B(R'\theta,r/2)}
	(T_{c}^{-1}|w|^2+|\nabla w|^2)e^{-\frac{|x-R'\theta|^2}{4(t'-t)}} dxdt. 
\end{aligned}
\]

Using \eqref{eq:winitialconclem} and \eqref{eq:timeintervalinclusion}, 
recalling that $t_0:=\eta^5 R^2 (\log(1/\epsilon_0))^{-1}$  
and $T_{c}\leq \eta^5 R^2$, we see that
\[
	Z\geq t_0 T_{c}^{-1} e^{-\frac{\eta^6 R^2}{t_0}}  \epsilon_{0} 
	\geq t_0(\eta^5 R^2)^{-1} e^{-\frac{\eta^6 R^2}{t_0}}  \epsilon_{0}
	= \frac{\epsilon_{0}^2}{\epsilon_{0}^{1-\eta} \log (1/\epsilon_0))}. 
\]
Then, for $0<\eta<1/2$ and $0<\epsilon_{0}\leq \eta^8$, we have 
\[
	Z\geq \frac{\epsilon_{0}^2}{\epsilon_{0}^{1/2} \log(1/\epsilon_0)}. 
\]
As $\lim_{x\to 0} x^{1/2} \log x=0$ 
and $0<\epsilon_{0}\leq \eta^8$, we infer that for $\eta$ sufficiently small, 
\begin{equation}\label{eq:Zlowerbound}
	Z\geq \epsilon_{0}^2.
\end{equation}

Now let us show that $e^{-r^2/(500 t_0)}X$ is much smaller than the lower bound of $Z$. Using \eqref{eq:diffinequality}, \eqref{eq:domaininclusioninS} and 
$T_{c}\in [0,\eta/4]$, we get that
\[
	e^{-\frac{r^2}{500 t_0}}X 
	\leq C(n) \epsilon_{0}^{\frac{1}{500\eta}} \eta^{2n}R^n\eta^{-\frac{2}{p-1}}.
\]
As $\epsilon_{0}\leq \min(\eta^8, R^{-10000n\eta})$ and $R\geq 20$, we see that for $\eta$ sufficiently small, 
\begin{equation}\label{eq:Xfirstest}
\begin{aligned}
	e^{-\frac{r^2}{500 t_0}}X
	&\leq C(n) \epsilon_{0}^{\frac{1}{500\eta}} \eta^{2n}R^n\eta^{-\frac{2}{p-1}}
	=C(n) \epsilon_{0}^{\frac{1}{1000\eta}}
	(\epsilon_{0}^{\frac{1}{2000\eta}}\eta^{2n-\frac{2}{p-1}})
	\epsilon_{0}^{\frac{1}{2000\eta}}R^n \\
	&\leq C(n) \epsilon_{0}^{\frac{1}{1000\eta}}
	(\eta^{\frac{1}{250\eta}+2n-\frac{2}{p-1}})R^{-4n}
	\leq C(n) \epsilon_{0}^{\frac{1}{1000\eta}}
	(\eta^{\frac{1}{250\eta}+2n-\frac{2}{p-1}}).
\end{aligned}
\end{equation}
As $\epsilon_{0}\leq \eta^8$, for $\eta$ sufficiently small, 
we have that 
\[
	\epsilon_{0}\in (0,1), \quad
	\frac{1}{1000\eta}\geq 2, \quad
	\frac{1}{250\eta}+2n-\frac{2}{p-1} \geq 3. 
\]
Using these and \eqref{eq:Xfirstest}, 
we get for $\eta$ sufficiently small that 
\[
	e^{-\frac{r^2}{500 t_0}}X 
	\leq C(n) \epsilon_{0}^{\frac{1}{1000\eta}}
	(\eta^{\frac{1}{250\eta}+2n-\frac{2}{p-1}})
	\leq C(n)\epsilon_{0}^{2}\eta^3 \leq \epsilon_{0}^{2}\eta^2.
\]
From this and \eqref{eq:Zlowerbound}, 
we see that for $\eta$ sufficiently small, 
$e^{-r^2/500 t_0}X$ is negligible compared to $Z$. 
Thus \eqref{e.carlslice} gives that
\[
	\epsilon_{0}^2\leq C(n) 
	(3e \eta^{-10})^{\tfrac{n}{2}+\frac{\log(1/\epsilon_0)}{200\eta}}
	\int_{B(R'\theta,r)}|w(x,t')|^2 e^{-\frac{|x-R'\theta|^2}{4 t_1}} dx . 
\]

Using $\epsilon_{0}\leq \eta^8$ gives for $\eta$ sufficiently small, 
\begin{equation}\label{eq:Yafterabsorp}
\begin{aligned}
	\epsilon_{0}^2
	&\leq C(n) \eta^{-8n-\frac{8\log(1/\epsilon_{0})}{50\eta}}
	\int_{B(R'\theta,r)}|w(x,t')|^2 
	e^{-\frac{|x-R'\theta|^2}{4 t_1}} dx\\
	&\leq C(n) 
	\left(\frac{1}{\epsilon_{0}}\right)^{n+\frac{1}{50\eta}\log(\frac{1}{\eta^8})}
	\int_{B(R'\theta,r)} |w(x,t')|^2 e^{-\frac{|x-R'\theta|^2}{4 t_1}} dx. 
\end{aligned}
\end{equation}
For $\eta$ sufficiently small we have that 
$n+(50\eta)^{-1} \log(1/\eta^8) \leq 1/(2\eta^2)$. 
From this, \eqref{eq:Yafterabsorp} and $r=\eta^2 R$, 
we get that 
\begin{equation}\label{eq:Ysplit}
	\epsilon_{0}^2 \leq C(n) {\epsilon_{0}}^{-\frac{1}{2\eta^2}} Y_{1}
	+C(n) {\epsilon_{0}}^{-\frac{1}{2\eta^2}} Y_{2}, 
\end{equation}
where 
\[
\begin{aligned}
	&Y_{1}:=
	\int_{\eta^5 R'\leq |x-R'\theta|<\eta^2 R} |w(x,t')|^2 
	e^{-\frac{|x-R'\theta|^2}{4 t_1}}  dx, \\
	&Y_{2}:=
	\int_{ |x-R'\theta|<\eta^5 R'}|w(x,t')|^2 
	e^{-\frac{|x-R'\theta|^2}{4 t_1}} dx.
\end{aligned}
\]
Using \eqref{eq:diffinequality} and $\epsilon_{0}\leq \min(\eta^8,R^{-10000n\eta})$, we get that 
\[
\begin{aligned} 
	C(n) \epsilon_{0}^{-\frac{1}{2\eta^2}} Y_{1}
	&\leq C(n) \epsilon_{0}^{-\frac{1}{2\eta^2}+\frac{(1+\eta^3)^2}{4\eta^5}} 
	\eta^{2n-\frac{2}{p-1}} R^{n}
	\leq C(n) \epsilon_{0}^{-\frac{1}{2\eta^2}
	+\frac{(1+\eta^3)^2}{4\eta^5}-\frac{1}{10000\eta}} 
	\eta^{2n-\frac{2}{p-1}}\\
	&= C(n) \epsilon_{0}^{-\frac{1}{2\eta^2}
	+\frac{(1+\eta^3)^2}{8\eta^5}-\frac{1}{10000\eta}} 
	\epsilon_{0}^{\frac{(1+\eta^3)^2}{8\eta^5}} \eta^{2n-\frac{2}{p-1}} \\
	&\leq C(n) {\epsilon_{0}}^{-\frac{1}{2\eta^2}
	+\frac{(1+\eta^3)^2}{8\eta^5}-\frac{1}{10000\eta}}
	\eta^{2n-\frac{2}{p-1}+\frac{(1+\eta^3)^2}{\eta^5}}.
\end{aligned}
\]
Using this and the fact that for  $\eta$ sufficiently small, 
\[
	-\frac{1}{2\eta^2}+\frac{(1+\eta^3)^2}{8\eta^5}-\frac{1}{10000\eta}
	\geq 2, \quad 
	2n-\frac{2}{p-1}+\frac{(1+\eta^3)^2}{\eta^5}\geq 3, 
\] 
we infer that
\[
	C(n) \epsilon_{0}^{-\frac{1}{2\eta^2}} Y_{1}
	\leq C(n)  \epsilon_{0}^2 \eta^3
	\leq \epsilon_{0}^2 \eta^2. 
\]
Thus, for $\eta$ sufficiently small, 
$C(n) \epsilon_{0}^{-1/(2\eta^2)} Y_{1}$ 
can be absorbed into the left-hand-side of \eqref{eq:Ysplit}. Then,  
\begin{equation}\label{eq:Y2est}
	\int_{ |x-R'\theta|<\eta^5 R'}|w(x,t')|^2  dx
	\geq 
	\frac{1}{C(n)} \epsilon_0^{ 2+\frac{1}{2\eta^2} }
	= \frac{1}{ C(n) \epsilon_0^{(1/(2\eta^2))-2}} 
	\epsilon_{0}^{\eta^{-2}}.
\end{equation}
As $\epsilon_{0}\leq \eta^8$, we see that for $\eta$ sufficiently small, 
\[
	\frac{1}{ C(n) \epsilon_0^{(1/(2\eta^2))-2}} 
	\geq 
	\frac{1}{ C(n) \epsilon_0} \geq 1. 
\]
Using this and \eqref{eq:Y2est}, we obtain that 
for $\eta$ sufficiently small, 
\[
	\int_{ |x-R'\theta|<\eta^5 R'}|w(x,t')|^2  dx
	\geq \epsilon_{0}^{{\eta^{-2}}}.
\]
As $t'\in  [-a + 2\eta^5 R^2 (\log(1/\epsilon_0))^{-1} ,0]$
was taken arbitrarily, we obtain the desired conclusion for $T_{1}=1.$ 
\end{proof}

Let us prove the iterated unique continuation.

\begin{proof}[Proof of Proposition \ref{pro:iterateduniquecontinuation}]
Proposition \ref{pro:iterateduniquecontinuation} 
is based on \cite[Proposition 4.2]{Pa22} however 
we required the final conclusion \eqref{eq:witeratedconc} 
to be over a larger domain of integration for the statement and proof to hold. 
The different condition on $a$ in Lemma \ref{lem:uniquecontslice} 
also necessitates  adjustments in the iterated parameters \eqref{eq:paracarliteratedef}. 
For the readers' convenience, 
we include the proof (following \cite[Proposition 4.2]{Pa22}) with these adjustments. 

By means of the rescaling 
\[
	w_{T_{1}}(x,t):=T_1^{\frac{1}{p-1}} w(T_1^{\frac{1}{2}} x, T_{1} t), 
\] 
we can assume without loss of generality that $T_{1}=1$.
For $k\in\N$, define
\begin{equation}\label{eq:paracarliteratedef}
\begin{aligned}
	&\epsilon_{k}:=\epsilon^{\eta^{-2k}},\quad 
	&&R_{k}:= R_{0}(1+\eta^3)^{k}, \\ 
	&a_{0}=\frac{3\eta}{4},  
	&&a_{k}:=\frac{3\eta}{4}-2\sum_{i=0}^{k-1} 
	\frac{\eta^5 R_{i}^{2}}{ \log( 1/\epsilon_i )}, 
\end{aligned}
\end{equation}
where $R_0$ is given in 
the assumption of Proposition \ref{pro:iterateduniquecontinuation}. 
Fix $R\geq R_{0}$. Then there exists $N\in\N$ such that  
\begin{equation}\label{eq:RcfRn}
	R_{N}\leq R<R_{N+1}\quad\mbox{ with }
	N=\lfloor\log_{1+\eta^3} \frac{R}{R_0}\rfloor.
\end{equation}
Let $x\in B(R_{N}\theta,\eta^5 R_{N})$. 
Then, from \eqref{eq:paracarliteratedef} and \eqref{eq:RcfRn}, 
it follows that
\[
\begin{aligned}
	|x-R\theta| &\leq |x-R_N\theta| + R-R_{N} 
	< \eta^5 R_N + R_{N+1}-R_N  \\
	&= (\eta^5+\eta^3) R_N\leq (\eta^5+\eta^3) R.  
\end{aligned}
\]
Then, for $\eta$ sufficiently small, 
\begin{equation}\label{eq:balliterateinclusion}
B(R\theta,\eta^{2} R)\supseteq B(R_{N}\theta,\eta^5 R_{N}).
\end{equation}

By L'H\^opital's rule
\[
	\lim_{x\rightarrow 0} \frac{x^{\frac{4}{3}} \log x}{\log(1+x)}=0.
\]
Thus, for all $\eta$ sufficiently small, 
\[
	\frac{2\eta^4 \log(\eta^3)}{3\log(1+\eta^3)}\geq -1, 
\]
and so
\[
\begin{aligned}
	\frac{2}{\log_{\eta}(1+\eta^3)}
	=\frac{2\log(\eta^3)}{3\log(1+\eta^3)}\geq -\eta^{-4}.
\end{aligned}
\]
From this and $R\geq R_0$, we infer that for $\eta$ sufficiently small, 
\[
	\eta^{-2N}\leq \eta^{-2\log_{1+\eta^3} \frac{R}{R_0}}
	=\eta^{-2\frac{\log_{\eta}(R/R_0)}{\log_{\eta}(1+\eta^3)}}
	=\left( \frac{R}{R_0}\right)^{-\frac{2}{\log_{\eta}(1+\eta^3)}}
	\leq \left( \frac{R}{R_0} \right)^{\eta^{-4}}. 
\]
As $\epsilon\leq\eta^8$, we infer from this that for $\eta$ sufficiently small, 
\begin{equation}\label{eq:epsilonNlower}
	\epsilon_{N}:=\epsilon^{\eta^{-2N}}
	\geq \epsilon^{\left( \frac{R}{R_0} \right)^{\eta^{-4}}}.
\end{equation}
Note that as $\epsilon\leq \eta^8$ and $R_0\leq \eta^{-2}$, 
we get for $\eta$ sufficiently small, 
\[
	0\leq 2\eta^5 \sum_{i=0}^{k-1} \frac{ R_{i}^{2} }{ \log( 1/\epsilon_{i} ) }
	\leq \frac{2\eta}{\log (1/\eta^8)}
	\sum_{i=0}^{k-1}(\eta+\eta^4)^{2i}
	\leq \frac{{2\eta} (\log(1/\eta^8))^{-1}}{1-(\eta+\eta^4)^2}
	<\frac{\eta}{4} 
\]
for $k=1,\ldots, N$. 
Thus, for $\eta$ sufficiently small, 
$a_{k}$ given by \eqref{eq:paracarliteratedef} satisfies
\begin{equation}\label{eq:aklemma}
	a_{k}\in ( \eta/2, 3\eta/4 ]
	\quad\mbox{ for all }k= 0,\ldots, N.
\end{equation}

From \eqref{eq:balliterateinclusion}, 
\eqref{eq:epsilonNlower} and \eqref{eq:aklemma}, 
we see that to show \eqref{eq:witeratedconc}, it suffices to show
\begin{equation}\label{eq:iteratedconcreduce}
	\int_{B(R_{N}\theta, \eta^5 R_{N})} |w(x,t)|^2 dx
	\geq \epsilon_{N}
	=\epsilon^{\eta^{-2N}} 
	\quad \mbox{ for all } t\in [-a_{N},0]. 
\end{equation} 
This is achieved by iteratively applying Lemma \ref{lem:uniquecontslice} 
for $k=0,\ldots, N$, 
with the parameters defined in \eqref{eq:paracarliteratedef}. 
To show that this can be done, 
it suffices to show that 
$\epsilon_{k}$ and $R_{k}$ satisfy 
the hypothesis  in Lemma \ref{lem:uniquecontslice} for each $k=0,\ldots, N$, 
where we note that \eqref{eq:aklemma} shows that 
$a_{k}$ satisfies the hypothesis in Lemma \ref{lem:uniquecontslice}. 
First it is clear that as $0<\epsilon\leq\eta^8$, 
for $\eta$ sufficiently small, 
\begin{equation}\label{eq:epsilonk1st}
	\epsilon_{k}:=\epsilon^{\eta^{-2k}}
	\leq \eta^8\quad \mbox{ for all }k=0,\ldots, N.
\end{equation}

Next, let us verify that
\begin{equation}\label{eq:epsilonk2nd}
	\epsilon_{k}\leq e^{-{40000n\eta^4 R_{k}^2}}
	\quad \mbox{ for all } k=0,\ldots, N.
\end{equation}
To see this, first note that for $\eta$ sufficiently small (independent of $k$), 
we have that
\[
	40000n(\eta+\eta^4)^{2k}
	\leq 8\log(1/\eta) \quad\mbox{ for all }k=0,\ldots, N.
\]
Using this, $R_{k}:= R_{0}(1+\eta^3)^{k}$, $R_0\leq \eta^{-2}$ 
and $\epsilon\leq\eta^8$, we get that for all $k=0,\ldots, N$, 
\[
	40000n \eta^4 R_{k}^2
	\leq 40000n (1+\eta^3)^{2k}
	\leq 8\eta^{-2k}\log(1/\eta)
	\leq \eta^{-2k}\log(1/\epsilon).
\]
Thus, for all $k=0,\ldots, N$, we obtain 
\[
	\epsilon_{k}
	:=\epsilon^{\eta^{-2k}}
	=e^{-\eta^{-2k}\log(1/\epsilon)}
	\leq e^{-40000n \eta^4 R_{k}^2}. 
\]

Finally, let us verify that 
\begin{equation}\label{eq:epsilonk3rd}
	\epsilon_{k}
	\leq (R_{k}^2)^{-5000n\eta}\qquad\forall k=0,\ldots, N.
\end{equation}
First note that for $\eta$ sufficiently small (independent of $k$), 
for all $k=0,\ldots, N$, 
\[
	(1+\eta^3)^{2k}\leq \frac{1}{\eta^{2k}}, \quad 
	\log\left( \frac{1}{\eta^{2k}} \right)\leq \frac{1}{\eta^{2k}}. 
\]
Using this and $\epsilon\leq\eta^8$, 
we infer that for $\eta$ sufficiently small (independent of $k$) 
and for all $k=0,\ldots, N$, 
\[
\begin{aligned}
	&5000 n\eta^{2k+1}\log(\eta^{-4}(1+\eta^3)^{2k}) \\
	&= 20000n\eta^{2k+1}\log ( 1/\eta )
	+5000n\eta^{2k+1}\log((1+\eta^3)^{2k}) \\
	&\leq 20000n\eta^{2k+1}\log(1/\eta)
	+5000n\eta^{2k+1}\log(1/(\eta^{2k})) \\
	&\leq 
	20000n\eta\log(1/\eta) +5000n\eta
	\leq 8\log(1/\eta)
	\leq \log(1/\epsilon).
\end{aligned}
\]
Using this, $R_{k}:= R_{0}(1+\eta^3)^{k}$ 
and the fact that $R_0\leq \eta^{-2}$, 
for $\eta$ sufficiently small (independent of $k$) 
we have that for all $k=0,\ldots, N$, 
\[
	5000 n\eta \log( R_{k}^2)
	\leq 5000 n\eta \log(\eta^{-4}(1+\eta^3)^{2k})
	\leq \eta^{-2k} \log(1/\epsilon). 
\]
Thus, for $\eta$ sufficiently small (independent of $k$), 
we have that for all $k=0,\ldots, N$
\[
	\epsilon_{k}
	:=\epsilon^{\eta^{-2k}}
	=e^{-\eta^{-2k}\log(1/\epsilon)}
	\leq e^{-5000n\eta\log(R_k^2)}=(R_{k}^2)^{-5000n\eta}, 
\]
as required. 
Combining \eqref{eq:aklemma} 
with \eqref{eq:epsilonk1st}, \eqref{eq:epsilonk2nd} and \eqref{eq:epsilonk3rd} 
means that we can iteratively apply Lemma \ref{lem:uniquecontslice} 
for $k=0,\ldots, N$, 
with the parameters defined in \eqref{eq:paracarliteratedef}. 
This gives \eqref{eq:iteratedconcreduce}, 
which implies \eqref{eq:witeratedconc} with $T_{1}=1$ 
by the aforementioned reasoning. 
The proof of Proposition \ref{pro:iterateduniquecontinuation} is complete. 
\end{proof}

\section{Proof of main quantitative estimate}\label{sec:prmain}
To prove our main quantitative estimate (Proposition \ref{pro:quantregheat}), 
we show the following proposition by applying the propositions 
in Sections \ref{sec:estiveps} and \ref{sec:Carl}.

\begin{proposition}\label{pro:quantestheat}
Let $n\geq3$, $p>p_S$ and $u$ be a classical solution 
of \eqref{eq:fujitaeq}. 
Then there exist constants $C_0, M_0>0$ depending only on $n$ and $p$ 
such that the following statement holds true. 
If there exist 
$M\geq M_{0}$, $x_0\in\R^n$ and $t'\in (-M^{-p(2p+3)},0)$ such that 
$u$ satisfies \eqref{eq:Type1bound2.0} and 
\begin{equation}\label{eq:initialconct'}
	(-t')^{\frac{2}{p-1}-\frac{n}{2}}
	\int_{t'}^{t'/2} 
	\int_{B(x_0, A(-t')^{\frac{1}{2}})} |u(x,t)|^{p+1} dxdt
	> \eps, 
\end{equation}
where $\eps$ and $A$ are constants given by \eqref{eq:epsfixed}, then 
we conclude that 
\begin{align}
	&\label{eq:lowerboundconcN}
	N:=\int_{ B(x_0,e^{e^{M^{C_0}}}) } 
	|u(x,0)|^{q_c} dx \geq M e^{-{e^{e^{M^{C_0}}}}}, \\
	&\label{eq:lowerboundconctime} 
	-t' \geq \exp\left( - N e^{{e^{e^{M^{C_0}}}}} \right). 
\end{align}
\end{proposition}

\begin{remark}
\eqref{eq:lowerboundconcN} is equivalent to 
$\exp( - N e^{{e^{e^{M^{C_0}}}}} ) \leq e^{-M}$, 
and so the lower bound \eqref{eq:lowerboundconctime} in the conclusion 
is consistent with the upper bound $-t'<M^{-p(2p+3)}$ in the assumption 
of Proposition \ref{pro:quantestheat}. 
\end{remark}

\begin{proof}[Proof of Proposition \ref{pro:quantestheat}]
Without loss of generality, we can assume $x_0=0$. 
Let $t'\in (-M^{-p(2p+3)}, 0)$ satisfy 
\eqref{eq:initialconct'}, where $M$ is sufficiently large. 
For the meaning of the terminology `$M$ being sufficiently large', 
see just before Proposition \ref{pro:counting}. 
We divide the proof into 12 steps. 

\noindent\textbf{Step 1: backward propagation}

Let $t''\in (-1/64,0)$ satisfy 
\begin{equation}\label{eq:t''def}
	-M^{-3} \leq t''\leq M^{(p-1)(2p+4)}  A^{\frac{n(p-1)}{p+1}} t'. 
\end{equation}
We note that \eqref{eq:wellseptimes} is satisfied. 
Now we have \eqref{eq:initialconct'} with $x_0=0$, 
the opposite inequality of 
\eqref{eq:forwardsmallness} with $x_0=0$.
Then the contraposition of Proposition \ref{pro:backprop} shows that 
the opposite inequality of \eqref{eq:backwardsmallness} with $x_0=0$
also holds, 
that is, the estimate 
\begin{equation}\label{eq:conccounting2}
	(-t'')^{\frac{2}{p-1}-\frac{n}{2}}
	\int_{t''}^{t''/2} 
	\int_{ B(0, A(-t'')^{\frac{1}{2}} )} |u(x,t)|^{p+1} dxdt> \eps 
\end{equation}
holds, where \eqref{eq:conccounting2} is the same as \eqref{eq:conccounting}.

\noindent\textbf{Step 2: preparation of slice of regularity}

By \eqref{eq:conccounting2} (or \eqref{eq:conccounting}) 
and Proposition \ref{pro:counting}, there exist a point 
\[
	z_*=(x_*,t_*)\in B(0, 20M_1(-t'')^\frac{1}{2})\times (2t'',t''/4)
\]
and  backward parabolic cylinders 
$Q( (x_*,t_*- r^2/8), \hat{\delta} r )$ and $Q(z_*, r)$ satisfying 
\begin{equation}\label{eq:cubeinclusion2}
	Q( (x_*,t_*- r^2/8), \hat{\delta} r )
	\subset Q(z_*, r) 
	\subset B(0, 20 M_{1}(-t'')^{\frac{1}{2}})
	\times (2t'', t''/4)
\end{equation}
with $r:=8M_{2}^{-(6n+2)}(-t'')^{1/2}$ and $\hat{\delta}:=M_{2}^{-1}$ 
such that 
\begin{align}
	&\label{eq:regcount2}
	\|u\|_{L^{\infty}(Q( z_*, r/2))}
	\leq C M^{-2} r^{-\frac{2}{p-1}},
	\quad
	\|\nabla u\|_{L^{\infty}(Q( z_*, r/2))}
	\leq C M^{-2} r^{-\frac{p+1}{p-1}}, \\
	&\label{eq:conccount2}
	\int_{Q((x_*,t_*- r^2/8), \hat{\delta} r)} |u|^2 dxdt
	\geq M_{2}^{-\bar{C}}(\hat{\delta} r)^{n+2-\frac{4}{p-1}}. 
\end{align}
We note that \eqref{eq:cubeinclusion2}, \eqref{eq:regcount2} and \eqref{eq:conccount2} 
are the same as \eqref{eq:cubeinclusion}, \eqref{eq:regcount} and \eqref{eq:conccount}.

In steps 2-4 we follow the arguments of \cite{Pa22} 
for the higher-dimensional Navier-Stokes equations. 
This enables the transfer of the lower bound \eqref{eq:conccount2} 
to large scales in a slice of regularity by using \eqref{eq:regcount2}, 
quantitative unique continuation and its iterated counterpart.

By Proposition \ref{pro:slices} with 
$R=\hat \delta r$, $T_1=1$, 
$x_0=x_*$, $t_0=t_*-(r^2/8)+100(\hat \delta r)^2$, 
there exists a direction $\theta\in \bS^{n-1}$ 
and a time interval 
$I\subset [t_*-(r^2/8)+99(\hat \delta r)^2, t_*-(r^2/8)+100(\hat \delta r)^2]$ 
with $|I|=(\hat \delta r/M^\gamma)^2$ such that 
within the slice 
\[
\begin{aligned}
	&S=S_*\times I \subset \R^n\times [-1,0], \\
	&S_*:=\{x\in\R^n;  
	\dist (x,x_*+\R_+ \theta )
	\leq 10 M^{-\gamma} | (x-x_*)\cdot \theta |, \; 
	|x-x_*| \geq 20 \hat \delta r\}, 
\end{aligned}
\]
we have 
\begin{equation}\label{eq:slrSS}
	\|u\|_{L^\infty(S)} \leq 
	M^{-1} \left( \frac{\hat \delta r}{M^\gamma} \right)^{-\frac{2}{p-1}}, 
	\quad 
	\|\nabla u\|_{L^\infty(S)} \leq 
	M^{-1} \left( \frac{\hat \delta r}{M^\gamma} \right)^{-\frac{p+1}{p-1}}. 
\end{equation}
Here, $\gamma=\gamma(n,p)$ is 
defined in Proposition \ref{pro:slices}.

\noindent\textbf{Step 3: unique continuation Carleman inequality}

Fix $\tilde t\in I$. 
Let us apply Proposition \ref{pro:uniquecont} for 
\[
	u_1(x,t):= u(x+x_*+30\hat \delta r \theta, t+ \tilde t), 
	\quad 
	x\in \overline{B(r/4)}, \; t\in [-1500000 n (\hat \delta r)^2, 0], 
\]
where $B(r/4)=B(0,r/4)$. Since 
\[
\begin{aligned}
	&x+x_*+30\hat \delta r \theta \in B(x_*,r/2) 
	&&\mbox{ for }x\in \overline{B(r/4)}, \\
	&t+ \tilde t \in (t_*-(r/2)^2, t_*) 
	&&\mbox{ for }t\in [-1500000 n (\hat \delta r)^2, 0], 
\end{aligned}
\]
we can see from \eqref{eq:regcount2} that 
\[
	|\partial_t u_1 - \Delta u_1|
	\leq \frac{C_1}{M^{2(p-1)} r^2} |u_1| 
	\quad 
	\mbox{ in }
	\overline{B(r/4)}\times  [-1500000 n (\hat \delta r)^2, 0]
\]
for some $C_1=C_1(n,p)>0$. Notice that 
$(r/4)^2 \geq 16000 \times  1500000 n (\hat \delta r)^2$. 
By Proposition \ref{pro:uniquecont} with $T_1=1500000 n (\hat \delta r)^2$, 
$C_0=M^{2(p-1)}/(1500000 C_1 n \hat \delta ^2)\geq1$ and 
\[
	\underline{s}= \left( \frac{\hat \delta^4}{M_3^4} \right) (\hat \delta r)^2 
	\leq \overline{s}= 75(\hat \delta r)^2 < 
	\frac{1500000 n (\hat \delta r)^2}{10000 n}, 
\]
we have $Z_1\leq X_1 + Y_1$, where 
\[
\begin{aligned}
	&X_1 := C(n)
	e^{-\frac{(r/4)^2}{500 \times 75(\hat \delta r)^2}} 
	\int_{-1500000 n (\hat \delta r)^2}^0 
	\int_{B(r/4)} 
	\left( \frac{|u_1|^2}{1500000 n (\hat \delta r)^2} 
	+ |\nabla u_1|^2 \right) dxdt, \\
	&Y_1 := C(n)
	\left( 75(\hat \delta r)^2 \right)^\frac{n}{2} 
	\left( \frac{ 3e \times 75(\hat \delta r)^2}
	{(\hat \delta^4/M_3^4) 
	(\hat \delta r)^2}  \right)^\frac{(r/4)^2}{200\times 75(\hat \delta r)^2} \\
	&\quad \times \int_{B(r/4)} 
	|u_1(x,0)|^2 \left( \frac{\hat \delta^4}{M_3^4} (\hat \delta r)^2 
	\right)^{-\frac{n}{2}}
	e^{-\frac{|x|^2}{4 (\hat \delta^4/M_3^4) (\hat \delta r)^2 }} dx, \\
	&Z_1:= 
	\int_{-2\times 75(\hat \delta r)^2}^{-75(\hat \delta r)^2} 
	\int_{B(r/8)} 
	\left( \frac{|u_1|^2}{1500000 n (\hat \delta r)^2} 
	+ |\nabla u_1|^2 \right) e^{-\frac{|x|^2}{4t}} dxdt. 
\end{aligned}
\]

From the change of variables and \eqref{eq:conccount2}, it follows that 
\[
\begin{aligned}
	Z_1&\geq  
	\frac{1}{C(n)(\hat \delta r)^2} \int_{-150(\hat \delta r)^2}^{-75(\hat \delta r)^2} 
	\int_{B(r/8)} 
	|u(x+x_*+30\hat \delta r \theta, t+ \tilde t)|^2  e^{-\frac{|x|^2}{4t}} dxdt \\
	&\geq 
	\frac{1}{C(n)(\hat \delta r)^2} 
	e^{-\frac{31^2 (\hat \delta r)^2}{4(-150)(\hat \delta r)^2}}
	\int_{-150(\hat \delta r)^2}^{-75(\hat \delta r)^2} 
	\int_{B(31\hat \delta r)} 
	|u(x+x_*+30\hat \delta r \theta, t+ \tilde t)|^2  dxdt \\
	&= 
	\frac{e^\frac{31^2}{600}}{C(n)(\hat \delta r)^2} 
	\int_{\tilde t-150(\hat \delta r)^2}^{\tilde t-75(\hat \delta r)^2} 
	\int_{B(x_*+30\hat \delta r\theta, 31\hat \delta r)} |u(y,s)|^2  dyds \\
	&\geq 
	\frac{1}{C(n)(\hat \delta r)^2} 
	\int_{t_*-(r^2/8)-(\hat \delta r)^2}^{t_* - (r^2/8)}  
	\int_{B(x_*,\hat \delta r)} |u(y,s)|^2  dyds \\
	&\geq 
	\frac{1}{C(n)(\hat \delta r)^2} M_2^{-\overline{C}(n,p)} 
	(\hat \delta r)^{n+2-\frac{4}{p-1}}. 
\end{aligned}
\]
Then by $\hat \delta=M_2^{-1}$, we obtain 
\begin{equation}\label{eq:Z1lowM2}
\begin{aligned}
	Z_1 \geq \frac{1}{C(n)} M_2^{\frac{4}{p-1}-n-\overline{C}(n,p)} 
	r^{n-\frac{4}{p-1}} 
	\geq e^{-M_2} r^{n-\frac{4}{p-1}}. 
\end{aligned}
\end{equation}
By \eqref{eq:regcount2} and $M^4 \geq \hat \delta$, 
we see that 
\[
\begin{aligned}
	X_1 &\leq 
	C e^{-\frac{1}{600000\hat \delta^2}} 
	\int_{-1500000 n (\hat \delta r)^2}^0 
	\int_{B(r/4)} 
	\left( \frac{ r^{-\frac{4}{p-1}} M^{-4} }{1500000 n (\hat \delta r)^2} 
	+ r^{-\frac{2(p+1)}{p-1}} M^{-4} \right) dxdt \\
	&\leq 
	\frac{C}{M^4} \left( 1 + \hat \delta^2 \right)
	r^{n-\frac{4}{p-1}}  
	e^{-\frac{M_2^2}{600000}} 
	\leq 
	e^{-\frac{M_2^2}{600000}} 
	r^{n-\frac{4}{p-1}}. 
\end{aligned}
\]
Then, $X_1$ can be absorbed in the right-hand side of \eqref{eq:Z1lowM2}, 
and so 
\begin{equation}\label{eq:Y1low2M2}
	Y_1 \geq e^{-2 M_2} r^{n-\frac{4}{p-1}}. 
\end{equation}

On the other hand, by splitting 
\[
\begin{aligned}
	Y_1 &\leq 
	C(n) \left( \frac{M_3^4}{\hat \delta^4} 
	 \right)^{\frac{n}{2} + \frac{1}{240000 \hat \delta^2}}
	\int_{B(\hat \delta r/ M_3)} 
	|u_1(x,0)|^2 e^{-\frac{M_3^4 |x|^2}{4 \hat \delta^4 (\hat \delta r)^2 }} dx \\
	&\quad 
	+ C(n) \left( \frac{M_3^4}{\hat \delta^4} 
	 \right)^{\frac{n}{2} + \frac{1}{240000 \hat \delta^2}}
	\int_{B(r/4)\setminus B(\hat \delta r/ M_3)} 
	|u_1(x,0)|^2 e^{-\frac{M_3^4 |x|^2}{4 \hat \delta^4 (\hat \delta r)^2 }} dx
	=: Y_1' + Y_1''
\end{aligned}
\]
and by using \eqref{eq:regcount2}, 
$(n/2) + (M_2^2/240000) \leq M_2^2/5$ and 
$M_3^{M_2^2}\leq e^{M_3^2}$, we see that 
\[
\begin{aligned}
	Y_1''&\leq 
	C(n) M^{-4} 
	\left( \frac{225e M_3^4}{\hat \delta^4} 
	 \right)^{\frac{n}{2} + \frac{1}{240000 \hat \delta^2}}
	r^{-\frac{4}{p-1}}
	\int_{B(r/4)\setminus B(\hat \delta r/ M_3)} 
	e^{-\frac{M_3^4 |x|^2}{4 \hat \delta^4 (\hat \delta r)^2 }} dx \\
	&\leq 
	\frac{C(n)}{M^4} ( 225e M_2^4 M_3^4 )^{\frac{n}{2} 
	+ \frac{M_2^2}{240000}}
	r^{n-\frac{4}{p-1}} e^{-\frac{M_3^2}{4 \hat \delta^4}} \\
	&\leq 
	( M_3^5 )^\frac{M_2^2}{5} 
	e^{-\frac{M_2^4 M_3^2}{8}}
	e^{-\frac{M_2^4 M_3^2}{8}} r^{n-\frac{4}{p-1}} 
	\leq 
	\frac{M_3^{M_2^2}}{ e^{M_3^2} } e^{-M_3^2} r^{n-\frac{4}{p-1}} 
	\leq 
	e^{-M_3^2} r^{n-\frac{4}{p-1}}. 
\end{aligned}
\]
Thus, $Y_1''$ can be absorbed in the right-hand side of \eqref{eq:Y1low2M2}. 
Then we have 
\[
\begin{aligned}
	e^{-3 M_2} r^{n-\frac{4}{p-1}} 
	&\leq Y_1' \leq 
	C ( M_2^4 M_3^4)^{\frac{n}{2} + \frac{1}{240000 \hat \delta^2}}
	\int_{B(\hat \delta r/ M_3)}  
	|u_1(x,0)|^2  dx \\ 
	&\leq 
	M_3^{M_2^2}
	\int_{B(x_*+30\hat \delta r \theta, \hat \delta r/ M_3)}  
	|u(y,\tilde t)|^2  dy. 
\end{aligned}
\]
Since $e^{-3 M_2} M_3^{-M_2^2} \geq e^{-M_3}$, we obtain 
\begin{equation}\label{eq:lowtildet}
	\int_{B(x_*+30\hat \delta r \theta, \hat \delta r/ M_3)}  
	|u(y,\tilde t)|^2  dy \geq e^{-M_3} r^{n-\frac{4}{p-1}}
	\quad \mbox{ for any }\tilde t\in I. 
\end{equation}

\noindent\textbf{Step 4: iterated unique continuation on thinner slice}

Let us apply Proposition \ref{pro:iterateduniquecontinuation} to 
\[
	u_2(x,t):= u(x+x_*, t+\tau), \quad 
	(x,t) \in \tilde S. 
\]
Here we set $\tau:= \sup I\leq 0$ and 
\[
\begin{aligned}
	&\tilde S:=\tilde S_0\times [-M_2^{-3}(\hat \delta r)^2,0] 
	\subset \R^n\times [-1,0], \\
	&\tilde S_0:=\{x\in\R^n; \; \dist (x,\R_+ \theta )
	\leq M_2^{-3} | x\cdot \theta |, \; |x| \geq 30 \hat \delta r\}. 
\end{aligned}
\]
For $x\in \tilde S_0$ and $t\in [-M_2^{-3}(\hat \delta r)^2,0]$, 
we can check that $x+x_*\in S_*$ and $t+\tau\in I$. 
Particularly, $[\tau-M_2^{-3} (\hat \delta r)^2,\tau]\subset I$. 
Indeed, 
\[
\begin{aligned}
	&\dist(x+x_*,x_*+\R_+\theta) 
	= \dist(x,\R_+\theta) 
	\leq M_2^{-3} | x\cdot \theta |
	\leq 10M^{-\gamma} | (x+x_*-x_*) \cdot \theta |, \\
	&|x+x_*-x_*| = |x| \geq 30\hat \delta r \geq 20\hat \delta r, \quad 
	t+\tau \leq \tau = \sup I, \\
	&t+\tau \geq \sup I - M_2^{-3}(\hat \delta r)^2 
	\geq \sup I - (\hat \delta r/M^\gamma)^2 = 
	\sup I - |I| = \inf I. 
\end{aligned}
\]
Then by \eqref{eq:slrSS}, we observe that 
\[
\begin{aligned}
	&\|u_2\|_{L^\infty(\tilde S)}
	\leq 
	\|u\|_{L^\infty(S)} \leq 
	M^{\frac{2\gamma}{p-1}-1} ( \hat \delta r )^{-\frac{2}{p-1}}
	\leq \left( M_2^{-3} \times (\hat \delta r)^2  \right)^{-\frac{1}{p-1}}
	, \\
	&\|\nabla u_2\|_{L^\infty(\tilde S)}
	\leq 
	\|\nabla u\|_{L^\infty(S)} \leq 
	M^{\frac{2\gamma}{p-1}-1} ( \hat \delta r )^{-\frac{2}{p-1}-1} 
	\leq \left( M_2^{-3} \times (\hat \delta r)^2 
	 \right)^{-\frac{1}{p-1}-\frac{1}{2}}, 
\end{aligned}
\]
and that the condition on the first line in \eqref{eq:diffinequality} is satisfied 
with $T_1= (\hat \delta r)^2$ and $\eta=M_2^{-3}$. 
Moreover, since 
\[
	|\partial_t u_2 - \Delta u_2| 
	\leq |u_2|^{p-1} |u_2| 
	\leq \frac{|u_2|}{M^{-(2\gamma-(p-1))} (\hat \delta r)^2 } 
	\leq \frac{|u_2|}{M_2^{-3} \times (\hat \delta r)^2 } 
\]
in $\tilde S$, the condition on the second line in \eqref{eq:diffinequality} 
is also satisfied with $C_0=1$. 
Notice that $0<\eta \leq \min( C_0^{-1},1) =1$. 
From \eqref{eq:lowtildet} and $[\tau-M_2^{-3} (\hat \delta r)^2,\tau]\subset I$, 
it follows that 
\[
	\int_{B(30\hat \delta r \theta, M_3^{-1} \hat \delta r)}  
	|u(x+x_*,t+\tau)|^2  dx \geq e^{-M_3} r^{n-\frac{4}{p-1}} 
\]
for $t\in [-M_2^{-3} (\hat \delta r)^2,0]$. 
In particular, 
\[
	\int_{B(30\hat \delta r \theta, 30 M_2^{-15} \hat \delta r)}  
	|u_2(x,t)|^2  dx 
	\geq e^{-M_3} r^{n-\frac{4}{p-1}} 
	= \hat\delta^{\frac{4}{p-1}-n} e^{-M_3} (\hat \delta r)^{n-\frac{4}{p-1}} 
\]
with 
$20 (\hat \delta r) \leq 30 \hat \delta r 
\leq M_2^6 (\hat \delta r)$ 
and $0< \hat\delta^{(4/(p-1))-n} e^{-M_3} \leq M_2^{-24}$. 
Then, Proposition \ref{pro:iterateduniquecontinuation} 
with $R_0=30\hat \delta r$, $\epsilon=\hat\delta^{(4/(p-1))-n} e^{-M_3}$ 
and $4/(p-1)<n$ yields 
\[
\begin{aligned}
	\int_{B(\tilde R_0 \theta, M_2^{-6} \tilde R_0)}  
	|u_2(x,t)|^2  dx 
	&\geq \left( \hat\delta^{\frac{4}{p-1}-n} e^{-M_3}
	\right)^{\left( \frac{\tilde R_0}{ 30\hat \delta r } \right)^{ M_2^{12} }}
	(\hat \delta r)^{n-\frac{4}{p-1}} \\
	&\geq e^{-M_3\left( \frac{\tilde R_0}{ 30\hat \delta r } \right)^{ M_2^{12} }}
	(\hat \delta r)^{n-\frac{4}{p-1}} 
\end{aligned}
\]
for all $t\in [-M_2^{-3} (\hat \delta r)^2/2,0]$ 
and $\tilde R_0\geq 30\hat \delta r$. 
Since $r=8M_2^{-(6n+2)}(-t'')^{1/2}$ 
and $\hat \delta=M_2^{-1}$, 
we have 
\begin{equation}\label{eq:u2lowee} 
\begin{aligned}
	&\int_{B(\tilde R_0 \theta, M_2^{-6} \tilde R_0)}  
	|u(x+x_*,t+\tau)|^2  dx \\
	&\geq 
	e^{-M_3\left( \frac{M_2^{6n+3} \tilde R_0}{ 240 (-t'')^{1/2} } \right)^{ M_2^{12} }} 
	(8 M_2^{-(6n+3)})^{n-\frac{4}{p-1}} (-t'')^{\frac{n}{2}-\frac{2}{p-1}} \\
	&\geq 
	e^{-2M_3 \left( 
	\frac{M_2^{6n+3} \tilde R_0}{ 240 (-t'')^{1/2} } \right)^{ M_2^{12} } }
	(-t'')^{\frac{n}{2}-\frac{2}{p-1}}
\end{aligned}
\end{equation}
for all $t\in [32M_2^{-(12n+9)} t'',0]$ 
and $\tilde R_0\geq 240M_2^{-(6n+3)}(-t'')^{1/2}$.

\noindent\textbf{Step 5: preparation of annulus of regularity}

By Proposition \ref{pro:annulus} with $\lambda=M_4$, $R=M_4$ 
and $T_1= M^3 (-t'')$, 
there exists 
\begin{equation}\label{eq:recallRhatnotat}
	M_4 \leq R_1 \leq M_4^{{M_4}^{M^{q_c + 1+ \beta}}}
\end{equation}
such that 
\begin{equation}\label{eq:annconcapp}
\begin{aligned}
	&\| u \|_{L^{\infty}( \cA)} 
	\leq 
	T_1^{-\frac{1}{p-1}} M^{-1}, 
	\quad 
	\| \nabla u \|_{L^{\infty}(\cA)} 
	\leq 
	T_1^{-\frac{p+1}{2(p-1)}}  M^{-1}, \\
	&\cA:= 
	\left\{ (x,t); \; 
	2 T_1^\frac{1}{2} R_1 <|x| 
	<\frac{1}{2} T_1^\frac{1}{2} R_1^{M_4}, 
	\; t\in \left[- \frac{1}{2} T_1, 0 \right]
	\right\}. 
\end{aligned}
\end{equation}
Based on \cite{Pa22} for the higher-dimensional Navier-Stokes equations, 
we apply quantitative backward uniqueness in a sub-annulus of $\mathcal{A}$ 
and use the lower bound \eqref{eq:u2lowee}. 
However, in \cite{Pa22}, the annulus used for backward uniqueness 
is similar to the annulus of regularity, 
which seems to create a gap in the arguments involving subsequent 
uses of unique continuation remaining inside the annulus of regularity. 
As in \cite{Ta21} and subsequently \cite{BP21}, 
we instead apply backward uniqueness Carleman inequalities 
in a much smaller annulus than the annulus of regularity.

Let $T_2:=M^{-2} T_1 = M(-t'')$. 
We define a smaller annulus $\tilde \cA$ by 
\[
\begin{aligned}
	\tilde \cA&:= 
	\left\{ (x,t); \; 
	M_3 T_2^\frac{1}{2}R_1 <|x| < M_3^{-1} T_2^\frac{1}{2} R_1^{M_4}, 
	\; t\in \left[- T_2, 0 \right]
	\right\} \\
	&\subset 
	\left\{ (x,t); \; 
	2 M T_2^\frac{1}{2} R_1 <|x| 
	<\frac{1}{2} M T_2^\frac{1}{2} R_1^{M_4}, 
	\; t\in \left[- \frac{1}{2} M^2 T_2, 0 \right]
	\right\} 
	= \cA. 
\end{aligned}
\]
Note that the estimates in \eqref{eq:annconcapp} 
also hold in $\tilde \cA$. 
We substitute $20M_3 T_2^{1/2} R_1 (= 20M^{1/2}M_3R_1(-t'')^{1/2} )$ 
into $\tilde R_0$ in \eqref{eq:u2lowee}. 
Then the change of variables yields 
\[
\begin{aligned}
	Z_2&:=\int_{B(x_* + 20M_3 T_2^{1/2} R_1 \theta, 
	20M_2^{-6} M_3 T_2^{1/2} R_1)}  
	|u(y,s)|^2  dy  \\
	&\geq 
	e^{-2M_3 \left( 
	\frac{M_2^{6n+3} 20M^{1/2}M_3R_1(-t'')^{1/2}}{ 240 (-t'')^{1/2} } 
	\right)^{ M_2^{12} } }
	(-t'')^{\frac{n}{2}-\frac{2}{p-1}} \\
	&\geq 
	e^{-2M_3 \left( 
	\frac{M_2^{6n+3} M^{1/2}M_3}{ 12 } R_1 
	\right)^{ M_2^{12} } }
	(-t'')^{\frac{n}{2}-\frac{2}{p-1}}
\end{aligned}
\]
for all 
$t\in [\tau+32M_2^{-(12n+9)} t'',\tau]
= [\tau- 32M^{-1} M_2^{-(12n+9)} T_2, \tau]$. 
Since $R_1 \geq M_4$, we see that 
\[
	Z_2\geq 
	e^{- R_1^{ 4M_2^{12} } }
	(-t'')^{\frac{n}{2}-\frac{2}{p-1}}
	\geq 
	e^{- R_1^{ 4M_2^{12} } }
	(M^{-1} T_2)^{\frac{n}{2}-\frac{2}{p-1}}
	\geq 
	e^{- R_1^{ M_3 } }
	T_2^{\frac{n}{2}-\frac{2}{p-1}}. 
\]
Hence by integrating $Z_2$ over $s\in [\tau- 32M^{-1} M_2^{-(12n+9)} T_2, \tau]$, 
we obtain 
\begin{equation}\label{eq:lowfarRr}
\begin{aligned}
	&\int_{\tau- \frac{ 32T_2 }{ M M_2^{12n+9}} }^\tau 
	\int_{B(x_* + 20M_3 T_2^\frac{1}{2} R_1 \theta, 
	20 M_3 M_2^{-6} T_2^\frac{1}{2} R_1)}  
	|u(x,t)|^2  dx dt \\
	&\geq 
	\frac{ 32 }{ M M_2^{12n+9}} e^{- R_1^{ M_3 } }
	T_2^{\frac{n}{2}-\frac{2}{p-1}+1}
	\geq 
	e^{-2R_1^{ M_3 }} T_2^{\frac{n}{2}-\frac{2}{p-1}+1}. 
\end{aligned}
\end{equation}

\noindent\textbf{Step 6: backward uniqueness Carleman inequality---two cases}

Let us apply Proposition \ref{pro:backuniqueness} to $u$. 
By \eqref{eq:annconcapp}, $M>1$ 
and $T_2=M^{-2} T_1$, we have 
\[
	|\partial_t u -\Delta u| 
	\leq 
	\frac{|u|}{ T_1 M^{p-1} } 
	\leq \frac{|u|}{ T_1 } 
	= \frac{|u|}{ M^2 \times T_2 } 
	\quad \mbox{ in }\tilde \cA. 
\]
We apply Proposition \ref{pro:backuniqueness} on $\tilde \cA$ 
with $r_- = M_3 T_2^{1/2} R_1$, $r_+ = M_3^{-1} T_2^{1/2} R_1^{M_4}$ 
and $C_\Carl= M^2$. 
Here we note that $C_\Carl \geq 3n$ 
and $r_-^2 \geq 4 C_\Carl T_2$ hold 
by taking $M$ so large that $M^2 \geq 3n$ and 
by using $R_1\geq M_4$, respectively.  
Then, we have $Z_3 \leq X_3 + Y_3$, where 
\[
\begin{aligned}
	&X_3:= C(n)  M^2 
	e^{-\frac{R_1^{1+M_4}}{4 M^2 } } 
	\iint_{\tilde \cA} e^\frac{2|x|^2}{ T_1 } 
	\left( \frac{|u|^2}{ T_2 } + |\nabla u|^2 \right) dxdt, \\
	&\begin{aligned}
	Y_3&:= 
	C(n)  M^2
	e^{-\frac{R_1^{1+M_4}}{4 M^2 } } 
	e^\frac{2 R_1^{2M_4}}{M^2 M_3^2 }  
	\int_{M_3 T_2^{1/2} R_1 \leq |x| \leq \frac{T_2^{1/2} R_1^{M_4} }{M_3} } 
	|u(x,0)|^2 dx, 
	\end{aligned} \\
	&Z_3 :=
	\int_{-\frac{T_2}{4}}^0 
	\int_{10 M_3 T_2^{1/2} R_1 \leq |x| \leq \frac{ T_2^{1/2}}{2M_3} R_1^{M_4} } 
	\left( \frac{|u|^2}{ T_2 } + |\nabla u|^2 \right) dxdt. 
\end{aligned}
\]

We estimate $Z_3$ by using \eqref{eq:lowfarRr}. 
Recall that 
$x_* \in B(0, 20M_1(-t'')^{1/2})$
in the first part of this section, 
and so $|x_*|\leq 20M^{-1/2} M_1 T_2^{1/2}$ since $T_2=M(-t'')$. 
For $x\in  B(x_* + 20M_3 T_2^{1/2} R_1 \theta, 20 M_2^{-6} M_3  T_2^{1/2} R_1)$, 
we observe that 
\[
\begin{aligned}
	&|x| \leq |x-x_* - 20M_3 T_2^\frac{1}{2} R_1 \theta| 
	+ |x_*| + 20M_3 T_2^\frac{1}{2} R_1 
	\leq 22 M_3 T_2^\frac{1}{2} R_1 
	\leq \frac{ T_2^\frac{1}{2}}{2M_3} R_1^{M_4}, \\
	&\begin{aligned}
	|x| &\geq 
	|20M_3 T_2^\frac{1}{2} R_1 \theta + x_* | 
	- |x-x_*-20M_3 T_2^\frac{1}{2} R_1 \theta| \\
	&\geq 20M_3 T_2^\frac{1}{2} R_1 
	- 20M^{-\frac{1}{2}} M_1 T_2^\frac{1}{2}
	 - 20 M_2^{-6} M_3  T_2^\frac{1}{2} R_1 
	 \geq 10 M_3 T_2^\frac{1}{2} R_1. 
	\end{aligned}
\end{aligned}
\]
Thus, 
\[
\begin{aligned}
	B(x_* + 20M_3 T_2^\frac{1}{2} R_1 \theta, 
	20 M_2^{-6} M_3  T_2^\frac{1}{2} R_1) 
	\subset 
	\left\{ 10 M_3 T_2^\frac{1}{2} R_1
	\leq |x| \leq \frac{ T_2^\frac{1}{2} }{2M_3} R_1^{M_4} \right\}. 
\end{aligned}
\]
Recall Step 2 that 
\[
\begin{aligned}
	&r=8M_2^{-(6n+2)} (-t'')^\frac{1}{2}=8M^{-1/2} M_2^{-(6n+2)} T_2^\frac{1}{2}, 
	\quad \hat \delta=M_2^{-1}, \\
	&\tau=\sup I\geq t_*-(r^2/8)+99(\hat \delta r)^2. 
\end{aligned}
\]
Moreover, the inclusions in \eqref{eq:cubeinclusion2} particularly show that 
$t_*-(r^2/8)-(\hat \delta r)^2\geq 2t''$. 
Then we see that  
\[
\begin{aligned}
	0&\geq \tau \geq \tau - \frac{ 32T_2 }{ M M_2^{12n+9}} 
	\geq 2t'' +100(\hat \delta r)^2 - \frac{ 32T_2 }{ M M_2^{12n+9}} \\
	&= 2t'' + \frac{6400 T_2}{M M_2^{12n+6}} - \frac{ 32T_2 }{ M M_2^{12n+9}}
	\geq 2t''= -\frac{2T_2}{M} \geq -\frac{T_2}{4}. 
\end{aligned}
\]
Then, $[\tau - 32 M^{-1} M_2^{-(12n+9)} T_2,\tau]\subset [-T_2/4,0]$. 
By \eqref{eq:lowfarRr}, we see that 
\[
\begin{aligned}
	Z_3 
	&\geq 
	\int_{\tau- \frac{ 32T_2 }{ M M_2^{12n+9}} }^\tau 
	\int_{B(x_* + 20M_3 T_2^\frac{1}{2} R_1 \theta, 
	20\frac{M_3}{M_2^6} T_2^\frac{1}{2} R_1)}  
	\frac{|u|^2}{ T_2 }  dxdt \\
	&\geq 
	\frac{1}{T_2} e^{-2R_1^{ M_3 }} T_2^{\frac{n}{2}-\frac{2}{p-1}+1}
	\geq 
	e^{-2R_1^{ M_3 }} T_2^{\frac{n}{2}-\frac{2}{p-1}}. 
\end{aligned}
\]
Hence, at least, one of the following holds: 
\begin{align}
	&\label{eq:Y3lowcl}
	Y_3\geq e^{-3R_1^{ M_3 }} T_2^{\frac{n}{2}-\frac{2}{p-1}}, \\
	&\label{eq:X3lowcl} 
	X_3\geq e^{-3R_1^{ M_3 }} T_2^{\frac{n}{2}-\frac{2}{p-1}}. 
\end{align}
Eventually, in both cases, we can derive the same estimate 
as can be seen in \eqref{eq:laY3case} and \eqref{eq:laX3casefinal} below. 
However, the second case requires more involved arguments than the first case.

\noindent\textbf{Step 7: estimate for first case}

If \eqref{eq:Y3lowcl} is satisfied, we see that 
\[
\begin{aligned}
	e^{-3R_1^{ M_3 }} T_2^{\frac{n}{2}-\frac{2}{p-1}} 
	&\leq Y_3 \leq 
	C(n) M^2 
	e^\frac{2 R_1^{2M_4} }{M^2 M_3^2 }  
	\int_{M_3 T_2^{1/2} R_1 \leq |x| \leq \frac{T_2^{1/2} R_1^{M_4} }{M_3} } 
	|u(x,0)|^2 dx \\
	&\leq 
	e^{R_1^{3M_4}}
	\int_{M_3 T_2^{1/2} R_1 \leq |x| \leq \frac{T_2^{1/2} R_1^{M_4} }{M_3} } 
	|u(x,0)|^2 dx. 
\end{aligned}
\]
Hence by expanding the integral domain, we obtain 
\begin{equation}\label{eq:laY3case}
	\int_{T_2^\frac{1}{2} R_1 \leq |x| \leq T_2^\frac{1}{2} R_1^{M_4}  } 
	|u(x,0)|^2 dx
	\geq 
	e^{-R_1^{ 4M_4 }} T_2^{\frac{n}{2}-\frac{2}{p-1}}. 
\end{equation}

\noindent\textbf{Step 8: estimate for second case 
by triple pigeonhole---first pigeonhole in space variable}

In steps 8-11, we show that \eqref{eq:X3lowcl} also implies \eqref{eq:laY3case} 
by utilizing arguments in \cite{Ta21} 
involving the pigeonhole principle 
and a further use of quantitative unique continuation.

We examine the case \eqref{eq:X3lowcl} and 
deduce the same inequality as \eqref{eq:laY3case}. 
By $T_2=M^{-2} T_1$ and $R_1\geq M_4$, we have 
\[
\begin{aligned}
	X_3 &\leq 
	C(n)  M^2 
	e^{-\frac{R_1^{1+M_4} }{4 M^2} } 
	\int_{-T_2}^0 
	\int_{M_3 T_2^\frac{1}{2} R_1 \leq |x| 
	\leq \frac{ T_2^\frac{1}{2}}{M_3} R_1^{M_4} } 
	e^\frac{2|x|^2}{ T_1 } 
	\left( \frac{|u|^2}{ T_2 } + |\nabla u|^2 \right) dxdt \\
	&\leq 
	e^{-\frac{R_1^{1+M_4} }{5 M^2} } 
	\int_{-\frac{T_1}{M^2}}^0 
	\int_{\frac{M_3}{M} T_1^\frac{1}{2} R_1 
	\leq |x| \leq \frac{ T_1^\frac{1}{2}}{MM_3} R_1^{M_4} } 
	e^\frac{2|x|^2}{ T_1 } 
	\left( \frac{M^2 |u|^2}{ T_1 } + |\nabla u|^2 \right) dxdt \\
	&\leq 
	e^{-\frac{R_1^{1+M_4} }{6 M^2} } 
	\int_{-\frac{T_1}{M^2}}^0 
	\int_{ M_3^\frac{1}{2} T_1^\frac{1}{2} R_1 
	\leq |x| \leq M_3^{-\frac{1}{2}} T_1^\frac{1}{2} R_1^{M_4} } 
	e^\frac{2|x|^2}{ T_1 } 
	\left( \frac{ |u|^2}{ T_1 } + |\nabla u|^2 \right) dxdt. 
\end{aligned}
\]
This together with \eqref{eq:X3lowcl} gives 
\[
\begin{aligned}
	X_3' &:= 
	\int_{-\frac{T_1}{M^2}}^0 
	\int_{ M_3^\frac{1}{2} T_1^\frac{1}{2} R_1 
	\leq |x| \leq M_3^{-\frac{1}{2}} T_1^\frac{1}{2} R_1^{M_4} } 
	e^\frac{2|x|^2}{ T_1 } 
	\left( \frac{ |u|^2}{ T_1 } + |\nabla u|^2 \right) dxdt \\
	&\geq 
	e^\frac{R_1^{1+M_4} }{6 M^2} 
	e^{-3R_1^{ M_3 }} M^{\frac{4}{p-1} - n } T_1^{\frac{n}{2}-\frac{2}{p-1}}
	\geq 
	e^\frac{R_1^{\frac{1}{2}+M_4} }{6 M^2}  T_1^{\frac{n}{2}-\frac{2}{p-1}}. 
\end{aligned}
\]

We claim that there exists 
$ M_3^{1/2} T_1^{1/2} R_1 \leq R_3 
\leq M_3^{-1/2} T_1^{1/2} R_1^{M_4}$ satisfying 
\[
	\int_{-\frac{T_1}{M^2}}^0 
	\int_{ R_3\leq |x| \leq 2R_3 } 
	e^\frac{2|x|^2}{ T_1 } 
	\left( \frac{ |u|^2}{ T_1 } + |\nabla u|^2 \right) dxdt
	\geq 
	e^{ \frac{R_1^{M_4} }{ 6 M^2} } T_1^{\frac{n}{2}-\frac{2}{p-1}}. 
\]
To obtain a contradiction, we suppose that 
\begin{equation}\label{eq:lowclaR3u2}
	\int_{-\frac{T_1}{M^2}}^0 
	\int_{ R_3\leq |x| \leq 2R_3 } 
	e^\frac{2|x|^2}{ T_1 } 
	\left( \frac{ |u|^2}{ T_1 } + |\nabla u|^2 \right) dxdt
	< 
	e^{ \frac{R_1^{M_4} }{ 6 M^2} } T_1^{\frac{n}{2}-\frac{2}{p-1}}. 
\end{equation}
for any $ M_3^{1/2} T_1^{1/2} R_1 \leq R_3 
\leq M_3^{-1/2} T_1^{1/2} R_1^{M_4}$. 
We take the positive integer $k'$ satisfying  
$2^{k'} \leq R_1^{M_4-1}/M_3 \leq 2^{k'+1}$. 
Note that 
\[
	M_3^{-\frac{1}{2}} T_1^\frac{1}{2} R_1^{M_4} 
	\leq 2^{k'+1} M_3^\frac{1}{2} T_1^\frac{1}{2} R_1
	\leq 2M_3^{-\frac{1}{2}} T_1^\frac{1}{2} R_1^{M_4}. 
\]
Then from \eqref{eq:lowclaR3u2}, it follows that 
\[
\begin{aligned}
	X_3' &\leq \sum_{k=0}^{k'} 
	\int_{-\frac{T_1}{M^2}}^0 
	\int_{2^k M_3^\frac{1}{2} T_1^\frac{1}{2} R_1 
	\leq |x| \leq 2^{k+1} M_3^\frac{1}{2} T_1^\frac{1}{2} R_1} 
	e^\frac{2|x|^2}{ T_1 } 
	\left( \frac{ |u|^2}{ T_1 } + |\nabla u|^2 \right) dxdt \\
	&< (k'+1) e^{ \frac{R_1^{M_4} }{ 6 M^2} } T_1^{\frac{n}{2}-\frac{2}{p-1}}. 
\end{aligned}
\]
Since $k'\leq (\log(R_1^{M_4-1}/M_3)/\log 2) \leq 2M_4 \log R_1 \leq R_1^2$ 
by $R_1\geq M_4$, 
we obtain 
\[
	X_3' < (R_1^2+1) e^{ \frac{R_1^{M_4} }{ 6 M^2} } T_1^{\frac{n}{2}-\frac{2}{p-1}}
	\leq e^{ \frac{R_1^{\frac{1}{2}+M_4} }{ 6 M^2} } T_1^{\frac{n}{2}-\frac{2}{p-1}} 
	\leq X_3', 
\]
a contradiction. Hence the claim follows. 
In particular, by $e^{ 2|x|^2/ T_1} \leq e^{8R_3^2/T_1}$ for 
$|x|\leq 2 R_3$ and by $e^{R_1^{M_4}/(6 M^2)}\geq 1$, we see that 
\begin{equation}\label{eq:tilX3low}
\begin{aligned}
	\tilde X_3&:=\int_{-\frac{T_1}{M^2}}^0 
	\int_{ R_3\leq |x| \leq 2R_3 } 
	\left( \frac{ |u|^2}{ T_1 } + |\nabla u|^2 \right) dxdt \\
	&\geq 
	e^{ \frac{R_1^{M_4} }{ 6 M^2} - \frac{8R_3^2}{T_1}} 
	T_1^{\frac{n}{2}-\frac{2}{p-1}} 
	\geq 
	e^{ - \frac{8R_3^2}{T_1}} 
	T_1^{\frac{n}{2}-\frac{2}{p-1}}. 
\end{aligned}
\end{equation}

\noindent\textbf{Step 9: estimate for second case 
by triple pigeonhole---second pigeonhole in time variable}

We decompose 
\[
	\tilde X_3 = 
	\int_{- \frac{T_1}{M^2} }^{- T_1 e^{ - \frac{10R_3^2}{T_1}} } 
	\int_{ R_3\leq |x| \leq 2R_3 }  (\cdots) 
	 + \int_{ - T_1 e^{ - \frac{10R_3^2}{T_1}} }^0 
	\int_{ R_3\leq |x| \leq 2R_3 }  (\cdots) 
	=: \tilde X_3' + \tilde X_3''. 
\]
Note that 
\[
	\{R_3\leq |x| \leq 2R_3 \} 
	\subset 
	\left\{ 2 T_1^\frac{1}{2} R_1 <|x| 
	<\frac{1}{2} T_1^\frac{1}{2} R_1^{M_4} \right\}. 
\]
Thus, we can apply \eqref{eq:annconcapp} to see that 
\[
	\tilde X_3'' 
	\leq 
	\frac{C(n)}{M^2} R_3^n T_1^{1-\frac{p+1}{p-1}}  
	e^{ - \frac{10R_3^2}{T_1}}  
	= 
	\frac{C(n)}{M^2} 
	\left( \frac{R_3^2}{ T_1 }\right)^\frac{n}{2} 
	e^{ - \frac{10R_3^2}{T_1}} T_1^{\frac{n}{2}-\frac{2}{p-1}} 
	\leq 
	e^{ - \frac{9R_3^2}{T_1}} T_1^{\frac{n}{2}-\frac{2}{p-1}}. 
\]
Hence $\tilde X_3''$ can be absorbed by the right-hand side of \eqref{eq:tilX3low}, 
and then 
\[
	\tilde X_3' 
	= \int_{- \frac{T_1}{M^2} }^{- T_1 e^{ - \frac{10R_3^2}{T_1}} } 
	\int_{ R_3\leq |x| \leq 2R_3 }  
	\left( \frac{ |u|^2}{ T_1 } + |\nabla u|^2 \right) dxdt
	\geq e^{ - \frac{9R_3^2}{T_1}} 
	T_1^{\frac{n}{2}-\frac{2}{p-1}}. 
\]

We claim that there exists 
$T_1 e^{ - 10R_3^2/T_1} \leq -t_3 \leq T_1/M^2$ 
satisfying 
\begin{equation}\label{eq:hX3lowb}
	\hat X_3 
	:= \int_{t_3}^{\frac{t_3}{2}} 
	\int_{R_3 \leq |x| \leq 2 R_3}
	\left( \frac{|u|^2}{ T_1 } + |\nabla u|^2 \right) dxdt 
	\geq e^{ - \frac{10R_3^2}{T_1}} 
	T_1^{\frac{n}{2}-\frac{2}{p-1}}. 
\end{equation}
To obtain a contradiction by the same argument as in Step 8, 
we suppose that 
\[
	\int_{t_3}^{\frac{t_3}{2}} 
	\int_{R_3 \leq |x| \leq 2 R_3}
	\left( \frac{|u|^2}{ T_1 } + |\nabla u|^2 \right) dxdt 
	< e^{ - \frac{10R_3^2}{T_1}} 
	T_1^{\frac{n}{2}-\frac{2}{p-1}}. 
\]
for any $T_1 e^{ - 10R_3^2/T_1} \leq -t_3 \leq T_1/M^2$. 
We take the positive integer $\tilde k'$ satisfying 
\[
	2^{-\tilde k'-1} \leq M^2 e^{ - \frac{10R_3^2}{T_1}}
	\leq 2^{-\tilde k'}. 
\]
Then, we can observe that 
\[
	- T_1 e^{ - \frac{10R_3^2}{T_1}} \leq 
	- 2^{-\tilde k'-1} M^{-2} T_1
	\leq - 2^{-1} T_1 e^{ - \frac{10R_3^2}{T_1} }. 
\] 
In addition, $\tilde k'\log 2 \leq 10 R_3^2/T_1 - \log M^2$. 
Thus, 
\[
\begin{aligned}
	\tilde X_3' &\leq \sum_{k=0}^{\tilde k'}
	\int_{-2^{-k} \frac{T_1}{M^2} }^{-2^{-k-1} \frac{T_1}{M^2}} 
	\int_{R_3 \leq |x| \leq 2 R_3} 
	\left( \frac{ |u|^2}{ T_1 } + |\nabla u|^2 \right) dxdt
	< (\tilde k'+1) e^{ - \frac{10R_3^2}{T_1}} 
	T_1^{\frac{n}{2}-\frac{2}{p-1}} \\
	&\leq 
	\left( \frac{10 }{\log 2}\left( \frac{R_3^2}{T_1} \right) 
	+1 \right) e^{ - \frac{10R_3^2}{T_1}} 
	T_1^{\frac{n}{2}-\frac{2}{p-1}}
	\leq 
	e^{ - \frac{9R_3^2}{T_1}} T_1^{\frac{n}{2}-\frac{2}{p-1}}
	\leq \tilde X_3', 
\end{aligned}
\]
a contradiction. The claim follows.

\noindent\textbf{Step 10: estimate for second case 
by triple pigeonhole---third pigeonhole by covering }

We take a covering 
$\{B(\tilde x_i, (-t_3)^{1/2} )\}_{i=1,\ldots,N_\cov}$ of 
$\{R_3 \leq |x| \leq 2 R_3\}$ such that 
$R_3 \leq |\tilde x_i| \leq 2 R_3$, where 
$N_\cov \leq C_\cov ( R_3/(-t_3)^{1/2})^n$ 
with some constant $C_\cov>0$ depending only on $n$. 
We claim that there exists 
$x_3:=\tilde x_{i'}$ ($1\leq i'\leq N_\cov$) satisfying 
\begin{equation}\label{eq:tripigeonT0x0}
	\int_{t_3}^{\frac{t_3}{2}} 
	\int_{B(x_3,  (-t_3)^\frac{1}{2})} 
	\left( \frac{|u|^2}{ T_1 } + |\nabla u|^2 \right) dxdt
	\geq e^{ - \frac{10n R_3^2}{T_1}} 
	T_1^{\frac{n}{2}-\frac{2}{p-1}}. 
\end{equation}
To obtain a contradiction, suppose that 
\[
	\int_{t_3}^{\frac{t_3}{2}} 
	\int_{B(\tilde x_i, (-t_3)^\frac{1}{2})} 
	\left( \frac{|u|^2}{ T_1 } + |\nabla u|^2 \right) dxdt
	< e^{ - \frac{10n R_3^2}{T_1}} 
	T_1^{\frac{n}{2}-\frac{2}{p-1}}. 
\]
for any $1\leq i\leq N_\cov$. 
By $-t_3 \geq T_1 e^{ - 10R_3^2/T_1}$, $n\geq3$ 
and \eqref{eq:hX3lowb}, 
\[
\begin{aligned}
	\hat X_3 &\leq 
	\sum_{i=1}^{N_\cov}
	\int_{t_3}^{\frac{t_3}{2}} 
	\int_{B(\tilde x_i, (-t_3)^\frac{1}{2})} 
	\left( \frac{|u|^2}{ r^2 } + |\nabla u|^2 \right) dxdt 
	< N_\cov e^{ - \frac{10n R_3^2}{T_1}} 
	T_1^{\frac{n}{2}-\frac{2}{p-1}} \\
	&\leq 
	C_\cov \left( \frac{R_3^2}{T_1} \right)^\frac{n}{2} 
	e^{ - \frac{5n R_3^2}{T_1}} 
	T_1^{\frac{n}{2}-\frac{2}{p-1}}
	\leq  
	e^{ - \frac{4n R_3^2}{T_1}} 
	T_1^{\frac{n}{2}-\frac{2}{p-1}} 
	\leq e^{ - \frac{12 R_3^2}{T_1}} 
	T_1^{\frac{n}{2}-\frac{2}{p-1}}
	\leq \hat X_3, 
\end{aligned}
\]
a contradiction. The claim follows. 

\noindent\textbf{Step 11: final estimate for second case 
by unique continuation Carleman inequality again}

We apply Proposition \ref{pro:uniquecont} to 
\[
	u_4(x,t):= u(x+x_3, t), \quad 
	x\in B(0, 100n T_1^{-\frac{1}{2}} R_3(-t_3)^\frac{1}{2} ),\;  t\in [10000n t_3,0].  
\]
We recall that 
\begin{equation}\label{eq:recM32M3x_0}
\begin{aligned}
	&R_3 \leq |x_3| \leq 2 R_3, 
	&& M_3^\frac{1}{2} T_1^\frac{1}{2} R_1 
	\leq R_3 \leq M_3^{-\frac{1}{2}} T_1^\frac{1}{2} R_1^{M_4}, \\
	&T_1 e^{ - \frac{10R_3^2}{T_1} } \leq -t_3 \leq M^{-2} T_1, 
	&&R_1\geq M_4. 
\end{aligned}
\end{equation}
Then, for $x\in B(0,  100n T_1^{-1/2} R_3(-t_3)^{1/2} )$, 
\[
\begin{aligned}
	&|x+x_3| \leq \left( 2+100n \left( \frac{-t_3}{T_1} \right)^{1/2} \right)R_3 
	\leq \left( \frac{2}{M_3^\frac{1}{2}}
	+ \frac{100n}{M M_3^\frac{1}{2}}  \right) T_1^\frac{1}{2} R_1^{M_4} 
	\leq \frac{1}{2} T_1^\frac{1}{2} R_1^{M_4}, 
	\\
	&|x+x_3| \geq \left( 1-100n \left( \frac{-t_3}{T_1} \right)^{1/2} \right)R_3 
	\geq \left( M_3^\frac{1}{2} 
	- \frac{100n M_3^\frac{1}{2}}{M} \right)
	T_1^\frac{1}{2} R_1 
	\geq 2 T_1^\frac{1}{2} R_1. 
\end{aligned}
\]
Moreover, for $t\in [10000n t_3,0]$, 
\[
	0\geq t\geq  10000n t_3 
	\geq - 10000n M^{-2} T_1
	\geq - \frac{1}{2} T_1. 
\]
Thus, $(x+x_3,t)\in \cA$ for 
$(x,t) \in B(0, 100n T_1^{-1/2} R_3(-t_3)^{1/2} ) \times [10000n t_3,0]$, 
where $\cA$ is defined in \eqref{eq:annconcapp}. 
By \eqref{eq:annconcapp}, we see that 
\[
	|\partial_t u_4 - \Delta u_4| 
	\leq \frac{|u_4|}{ T_1 M^{p-1} } 
	\leq \frac{|u_4|}{ T_1 }
	= \frac{|u_4|}{(10000n)^{-1} ( T_1/(-t_3)) \times 10000n (-t_3)} 
\]
in $B(0; 100n T_1^{-1/2} R_3(-t_3)^{1/2} )\times [10000n t_3,0]$. 
Since 
\[
\begin{aligned}
	&\frac{T_1/(-t_3)}{10000n} 
	\geq \frac{M^2}{10000n} \geq 1,  \\
	&(100n T_1^{-\frac{1}{2}} R_3(-t_3)^\frac{1}{2})^2 
	\geq 
	10000n^2 M_3  R_1^2(-t_3) 
	\geq 16000 \times 10000n(-t_3), 
\end{aligned}
\]
we can apply Proposition \ref{pro:uniquecont} 
to $u_4$ on $\mathcal{C}=B(0,100n T_1^{-1/2} R_3(-t_3)^{1/2})
\times [10000n t_3,0]$ 
with $\underline{s} =  \overline{s} = (-t_3)/2$. 
Then we obtain $Z_4 \leq X_4 + Y_4$, where 
\[
\begin{aligned}
	&\begin{aligned}
	X_4 &:= C(n) e^{-\frac{(100n T_1^{-1/2} R_3(-t_3)^{1/2})^2}
	{500 (-t_3/2) }} \\
	&\quad \times 
	\int_{10000n t_3}^0 \int_{ B(0, 100n T_1^{-\frac{1}{2}} R_3(-t_3)^\frac{1}{2}) }
	\left( \frac{|u_4|^2}{ 10000n (-t_3) } + |\nabla u_4|^2 
	\right) dxdt, 
	\end{aligned} \\
	&\begin{aligned}
	Y_4 &:= C(n) ( -t_3/2 )^\frac{n}{2} 
	( 3e )^\frac{(100n T_1^{-1/2} R_3(-t_3)^{1/2})^2}{200 (-t_3/2)} \\
	&\quad \times \int_{ B(0, 100n T_1^{-\frac{1}{2}} R_3(-t_3)^\frac{1}{2} ) }
	|u_4(x,0)|^2 ( -t_3/2 )^{-\frac{n}{2}} 
	e^{-\frac{|x|^2}{4 (-t_3/2) }} dx, 
	\end{aligned}\\
	&Z_4 := 
	\int_{t_3}^{\frac{t_3}{2}}
	\int_{ B(0, 50n T_1^{-\frac{1}{2}} R_3(-t_3)^\frac{1}{2}) } 
	\left( \frac{|u_4|^2}{ 10000n (-t_3) } + |\nabla u_4|^2 \right) 
	e^\frac{|x|^2}{4t} dxdt. 
\end{aligned}
\]

From the change of variables, the relation 
\[
	50n T_1^{-\frac{1}{2}} R_3 (-t_3)^\frac{1}{2}\geq 
	50n M_3^\frac{1}{2} (-t_3)^\frac{1}{2} \geq (-t_3)^\frac{1}{2}
\]
and \eqref{eq:tripigeonT0x0} with $-t_3 \leq M^{-2} T_1$, 
it follows that 
\begin{equation}\label{eq:Z4low9n}
\begin{aligned}
	Z_4 &\geq 
	\int_{t_3}^{\frac{t_3}{2}}
	\int_{ B(0, 50n T_1^{-\frac{1}{2}} R_3(-t_3)^\frac{1}{2}) } 
	\frac{|u_4|^2}{ 10000n (-t_3) } e^\frac{|x|^2}{4t} dxdt \\
	&\geq 
	\frac{T_1}{10000n (-t_3)} e^\frac{-t_3}{4(t_3/2)} 
	\int_{t_3}^{\frac{t_3}{2}}
	\int_{ B(x_3, (-t_3)^\frac{1}{2}) } 
	\frac{|u|^2}{T_1} dxdt \\
	&\geq 
	\frac{T_1 }{10000n e^{1/2}(-t_3)} 
	e^{ - \frac{10n R_3^2}{T_1}} 
	T_1^{\frac{n}{2}-\frac{2}{p-1}} 
	\geq 
	\frac{M^2}{10000ne^{1/2}} 
	e^{ - \frac{10n R_3^2}{T_1}} 
	T_1^{\frac{n}{2}-\frac{2}{p-1}} \\
	&\geq 
	e^{ - \frac{10n R_3^2}{T_1}} 
	T_1^{\frac{n}{2}-\frac{2}{p-1}}. 
\end{aligned}
\end{equation}
On the other hand, by \eqref{eq:annconcapp}, 
$T_1 e^{ - 10R_3^2/T_1} \leq -t_3 \leq M^{-2} T_1$ and $n\geq3$, 
we see that 
\[
\begin{aligned}
	X_4 &\leq 
	C(n) e^{-\frac{10000n^2 R_3^2}{250T_1}} 
	(10000n (-t_3)) ( 100n T_1^{-\frac{1}{2}} R_3(-t_3)^\frac{1}{2})^n \\
	&\quad \times \left( \frac{(T_1^{-\frac{1}{p-1}}M^{-1})^2}{ 10000n (-t_3) } 
	+ (T_1^{-\frac{p+1}{2(p-1)}}M^{-1})^2 
	\right) 
	 \\
	&\leq  
	C(n) M^{-4}  e^{-\frac{40 n^2 R_3^2}{T_1}} 
	\left( \frac{R_3^2}{T_1} \right)^\frac{n}{2} 
	\left( e^{\frac{10R_3^2}{T_1}}+ 1 \right) 
	T_1^{-\frac{2}{p-1}} (M^{-2} T_1)^\frac{n}{2} \\
	&\leq 
	C(n)  M^{-n-4} 
	e^{-\frac{20n^2 R_3^2}{T_1}} 
	T_1^{\frac{n}{2}-\frac{2}{p-1}}.  
\end{aligned}
\]
From $R_3^2/T_1 \geq M_3 R_1^2 \geq M_3 M_4^2$, it follows that 
\[
\begin{aligned}
	X_4 &\leq C(n)  M^{-n-4} 
	e^{-\frac{10n^2 R_3^2}{T_1}} 
	e^{-\frac{10n^2 R_3^2}{T_1}} 
	T_1^{\frac{n}{2}-\frac{2}{p-1}}  \\
	&\leq
	C(n)  M^{-n-4} 
	e^{-10n^2 M_3 M_4^2} 
	e^{-\frac{10 n^2 R_3^2}{T_1}} 
	T_1^{\frac{n}{2}-\frac{2}{p-1}} 
	\leq
	e^{-\frac{10n^2 R_3^2}{T_1}} 
	T_1^{\frac{n}{2}-\frac{2}{p-1}}. 
\end{aligned}
\]
Therefore, $X_4$ can be absorbed by the right-hand side of \eqref{eq:Z4low9n}, 
and so 
\begin{equation}\label{eq:Y4low20}
	Y_4\geq 
	e^{ - \frac{20n R_3^2}{T_1}} 
	T_1^{\frac{n}{2}-\frac{2}{p-1}}
	\geq e^{ - \frac{20n^2 R_3^2}{T_1}} 
	T_1^{\frac{n}{2}-\frac{2}{p-1}}. 
\end{equation}

For $Y_4$, using $R_3^2/T_1 \geq M_3 M_4^2$, we compute that 
\[
\begin{aligned}
	Y_4 
	&\leq  
	C(n) ( 3e )^\frac{100n^2 R_3^2}{ T_1} 
	\int_{ B(0, 100n T_1^{-\frac{1}{2}} R_3(-t_3)^\frac{1}{2} ) }
	|u_4(x,0)|^2 dx \\
	&\leq 
	C(n) ( 3e )^{-\frac{100n^2 R_3^2}{T_1}} 
	( e^3 )^\frac{200n^2 R_3^2}{T_1}  
	\int_{ B(0, 100n T_1^{-\frac{1}{2}} R_3(-t_3)^\frac{1}{2} ) }
	|u_4(x,0)|^2 dx \\
	&\leq 
	e^\frac{600n^2 R_3^2}{T_1} 
	\int_{ B(0, 100n T_1^{-\frac{1}{2}} R_3(-t_3)^\frac{1}{2} ) }
	|u_4(x,0)|^2 dx. 
\end{aligned}
\]
From this and \eqref{eq:Y4low20} (with a change of variables), 
we infer that
\[
	Y_4':=\int_{ B(x_3, 100n T_1^{-\frac{1}{2}} R_3(-t_3)^\frac{1}{2} ) }
	|u(x,0)|^2 dx 
	\geq 
	e^{-\frac{620n^2 R_3^2}{T_1}} 
	T_1^{\frac{n}{2}-\frac{2}{p-1}}. 
\]
Recall \eqref{eq:recM32M3x_0} and $T_2=M^{-2} T_1$. 
Then, for $x\in B(x_3, 100n T_1^{-1/2} R_3(-t_3)^{1/2} )$, 
\[
\begin{aligned}
	&\begin{aligned}
	|x| &\leq |x_3| + |x-x_3| \leq 
	\left( 2+ 100n \left(\frac{-t_3}{T_1}\right)^{1/2} \right)R_3
	\leq 
	\left( 2+ \frac{100n}{M} \right)R_3 \\
	&\leq 3 M_3^{-\frac{1}{2}} T_1^\frac{1}{2} R_1^{M_4} 
	= 3 M_3^{-\frac{1}{2}} M T_2^\frac{1}{2} R_1^{M_4} 
	\leq T_2^\frac{1}{2} R_1^{M_4} , 
	\end{aligned}\\
	&\begin{aligned}
	|x| &\geq |x_3| - |x-x_3| 
	\geq \left( 1- 100n T_1^{-\frac{1}{2}} (-t_3)^\frac{1}{2} \right)R_3 
	\geq \left( 1- \frac{100n}{M} \right)R_3 \\
	&\geq 
	\frac{1}{2} M_3^\frac{1}{2} T_1^\frac{1}{2} R_1 
	= \frac{1}{2} M_3^\frac{1}{2} M T_2^\frac{1}{2} R_1 
	\geq T_2^\frac{1}{2} R_1. 
	\end{aligned}
\end{aligned}
\]
and so $B(x_3, 100n T_1^{-1/2} R_3(-t_3)^{1/2} )\subset 
\{ T_2^{1/2} R_1 \leq |x| \leq T_2^{1/2} R_1^{M_4}\}$. 
From these together with 
$R_3 \leq M_3^{-1/2} T_1^{1/2} R_1^{M_4}$ and $R_1\geq M_4$, 
it follows that 
\begin{equation}\label{eq:laX3casefinal}
\begin{aligned}
	&\int_{ T_2^\frac{1}{2} R_1 \leq |x| \leq T_2^\frac{1}{2} R_1^{M_4} } 
	|u(x,0)|^2  dx 
	\geq Y_4' 
	\geq 
	e^{-\frac{620n^2 R_3^2}{T_1}} 
	T_1^{\frac{n}{2}-\frac{2}{p-1}} \\
	&= 
	M^{n-\frac{4}{p-1}} 
	e^{-\frac{620n^2 R_3^2}{T_1}} 
	T_2^{\frac{n}{2}-\frac{2}{p-1}} 
	\geq 
	M^{n-\frac{4}{p-1}} e^{-\frac{620n^2}{M_3} R_1^{2M_4} } 
	T_2^{\frac{n}{2}-\frac{2}{p-1}} \\
	&\geq 
	M^{n-\frac{4}{p-1}} e^{- R_1^{2M_4} } 
	T_2^{\frac{n}{2}-\frac{2}{p-1}} 
	\geq 
	e^{- R_1^{4M_4} } T_2^{\frac{n}{2}-\frac{2}{p-1}}. 
\end{aligned}
\end{equation}
Therefore, \eqref{eq:laY3case} also holds in this case. 
In particular, \eqref{eq:laY3case} holds in both cases \eqref{eq:Y3lowcl} and 
\eqref{eq:X3lowcl}. 
By applying the H\"older inequality, we obtain 
\[
\begin{aligned}
	&e^{-R_1^{ 4M_4 }} T_2^{\frac{n}{2}-\frac{2}{p-1}}
	\leq 
	\int_{T_2^\frac{1}{2} R_1 \leq |x| \leq T_2^\frac{1}{2} R_1^{M_4}  } 
	|u(x,0)|^2 dx  \\
	&\leq 
	C \left( \int_{T_2^\frac{1}{2} R_1 \leq |x| \leq T_2^\frac{1}{2} R_1^{M_4}  } 
	|u(x,0)|^{q_c} dx  \right)^\frac{2}{q_c}
	\left( T_2^\frac{1}{2} R_1^{M_4} \right)^{n(1-\frac{2}{q_c})}, 
\end{aligned}
\]
and so 
\[
\begin{aligned}
	&\int_{T_2^\frac{1}{2} R_1 \leq |x| \leq T_2^\frac{1}{2} R_1^{M_4}  } 
	|u(x,0)|^{q_c} dx  
	\geq 
	\frac{1}{C} R_1^{-n(\frac{q_c}{2}-1)M_4} 
	e^{-\frac{q_c}{2}R_1^{ 4M_4 }} \\
	&\geq 
	\frac{1}{C} e^{-R_1^{M_4}} 
	e^{-\frac{q_c}{2}R_1^{ 4M_4 }} 
	\geq 
	e^{-R_1^{ 5M_4 }}. 
\end{aligned}
\]

\noindent\textbf{Step 12: final summation}

In analogy to \cite{Ta21} for the Navier-Stokes equations, 
we use the conclusion of Step 11 
and a summing of disjoint spatial scales to conclude. 
This then implies the lower bound \eqref{eq:lowerboundconctime}.

From \eqref{eq:recallRhatnotat} and \eqref{eq:laY3case}, it follows that 
\[
\begin{aligned}
	&\int_{T_2^\frac{1}{2} R_1 \leq |x| \leq T_2^\frac{1}{2} R_1^{M_4}  } 
	|u(x,0)|^{q_c} dx  
	\geq 
	e^{-R_1^{ 5M_4 }}
	\geq 
	e^{-(M_4^{{M_4}^{M^{q_c + 1+ \beta}}})^{5 M_4}}  \\
	&= e^{-M_4^{{5 M_4 \times M_4}^{M^{q_c + 1+ \beta}}}} 
	= e^{-M_4^{{5 M_4}^{1+M^{q_c + 1+ \beta}}}} 
	\geq e^{-M_4^{e^{M_4}}} \geq e^{-e^{e^{M_5}}}. 
\end{aligned}
\]
Similarly, we have $1\leq M_4 \leq R_1 \leq R_1^{M_4} \leq e^{ e^{M_5} }$. 
Thus, 
\begin{equation}\label{eq:concloweee}
	\int_{ T_2^\frac{1}{2} \leq |x| \leq e^{ {e}^{M_5}} T_2^\frac{1}{2}} 
	|u(x,0)|^{q_c}  dx 
	\geq 
	e^{-e^{e^{M_5}}}. 
\end{equation} 
Recall that $T_2=M^{-2} T_1 = M(-t'')$ and 
$t''$ is arbitrarily chosen as in \eqref{eq:t''def}, that is, 
\[
	M^{(p-1)(2p+4)}  A^{\frac{n(p-1)}{p+1}}(-t')
	\leq - t''\leq M^{-3}. 
\]
Thus, we can take any $T_2\in [T_{\mi},T_{\ma}]$, where 
\[
	T_\mi:= M^{(p-1)(2p+4)+1} A^\frac{n(p-1)}{p+1} (-t'), 
	\quad T_\ma:= M^{-2}. 
\]
This together with \eqref{eq:concloweee} implies that 
\begin{equation}\label{eq:lowkanqc}
	\int_{ (e^{{e}^{M_5}})^{k-1} T_\mi^\frac{1}{2} 
	\leq |x| \leq (e^{ {e}^{M_5}})^k T_\mi^\frac{1}{2} } 
	|u(x,0)|^{q_c}  dx 
	\geq 
	e^{- e^{ e^{M_5}}} 
\end{equation}
holds for $k=1,\ldots,k''-1$, where $k''$ is the integer satisfying 
\begin{equation}\label{eq:kddesti}
	(e^{ {e}^{M_5}})^{k''-1} T_\mi^\frac{1}{2} 
	\leq 
	e^{ {e}^{M_5}} T_\ma^\frac{1}{2} 
	\leq 
	(e^{ {e}^{M_5}})^{k''} T_\mi^\frac{1}{2}
	\left( = e^{ k'' {e}^{M_5}} T_\mi^\frac{1}{2}\right). 
\end{equation}

Since $A$ in the definition of $T_\mi$ 
is determined in \eqref{eq:epsfixed} 
(see also Proposition \ref{pro:epsreg} 
for $\eps=\eps_0/2$ in the definition of $A$), 
we have 
\[
\begin{aligned}
	k'' &\geq {e}^{-M_5} \log \left( e^{ {e}^{M_5}} 
	\frac{T_\ma^{1/2}}{T_\mi^{1/2}} \right) 
	= {e}^{-M_5} \log \left( e^{ {e}^{M_5}} 
	\frac{1}{M^{(p-1)(p+2)+\frac{3}{2}} A^\frac{n(p-1)}{2(p+1)} 
	(-t')^\frac{1}{2} } \right) \\
	&= 1+ {e}^{-M_5} \log \left( 
	\frac{1}{M^{(p-1)(p+2)+\frac{3}{2}} 
	( 48\log(\frac{M^{p+1}}{\eps_0/2})  )^\frac{n(p-1)}{2(p+1)} 
	(-t')^\frac{1}{2} } \right) \\
	&\geq 
	1+ {e}^{-M_5} \log \left( 
	\frac{1}{M^{(p-1)(p+2)+\frac{3}{2}} 
	(48\log(\frac{C M^{p+1}}{ M^{-2(p+1)^2}}))^\frac{n(p-1)}{2(p+1)} 
	(-t' )^\frac{1}{2} } \right).  
\end{aligned}
\]
For $M$ sufficiently large, we note that 
\begin{equation}\label{eq:48npnp12}
	\left(48 \log \frac{C M^{p+1}}{ M^{-2(p+1)^2}}  \right)^\frac{n(p-1)}{2(p+1)} 
	\leq 
	\left( M^\frac{(p+1)}{n(p-1)} \right)^\frac{n(p-1)}{2(p+1)} 
	\leq M^\frac{1}{2}. 
\end{equation}
Now choose $0<\sigma_p<1/2$ depending only on $p$ and satisfying 
\[
	(p-1)(p+2)+1-p(2p+3)\left( \frac{1}{2}-\sigma_p \right)
	= -1 -\frac{1}{2}p 
	+p(2p+3)\sigma_p <0. 
\]
Then, from this choice of $\sigma_p$ and \eqref{eq:48npnp12}, we get
\[
\begin{aligned}
	k'' &\geq 
	1+ {e}^{-M_5} \log \left( 
	\frac{1}{M^{(p-1)(p+2)+1} (-t' )^\frac{1}{2}} \right) \\
	&\geq 
	1+ {e}^{-M_5} \log \left( 
	\frac{1}{M^{(p-1)(p+2)+1} 
	(-t' )^{\frac{1}{2}-\sigma_p} (-t' )^{\sigma_p} } \right) \\
	&\geq 
	1+ {e}^{-M_5} \log \left( 
	\frac{1}{M^{(p-1)(p+2)+1} (-t')^{\frac{1}{2}-\sigma_p} } \right) 
	+ \sigma_p {e}^{-M_5} \log \frac{1}{-t'}. 
\end{aligned}
\]
Recall from our standing assumptions in Proposition \ref{pro:quantestheat} 
that $-t'<M^{-p(2p+3)}$.
Then we note that $M^{(p-1)(p+2)+1} (-t')^{(1/2)-\sigma_p} <1$ and obtain  
\[
	k'' \geq 
	1+ \sigma_p {e}^{-M_5} \log \frac{1}{-t'} 
	\geq 1+ {e}^{-2M_5} \log \frac{1}{-t'}. 
\]

Summing up the lower estimates \eqref{eq:lowkanqc} for $k=1,\ldots, k''-1$, 
we have 
\[
\begin{aligned}
	&\int_{ T_\mi^\frac{1}{2} 
	\leq |x| \leq (e^{ {e}^{M_5}})^{k''-1} T_\mi^\frac{1}{2} } 
	|u(x,0)|^{q_c}  dx 
	\geq (k''-1)e^{-e^{e^{M_5}}} \\
	&\geq {e}^{-2M_5} e^{-e^{e^{M_5}}} \log \frac{1}{-t'} 
	\geq e^{-e^{e^{M_6}}} \log \frac{1}{-t'}. 
\end{aligned}
\]
By \eqref{eq:kddesti} and $T_\ma=M^{-2}$, we also have 
\[
	\{ T_\mi^\frac{1}{2} 
	\leq |x| \leq (e^{ {e}^{M_5}})^{k''-1} T_\mi^\frac{1}{2} \} 
	\subset 
	\{ |x| \leq e^{ {e}^{M_5}} T_\ma^\frac{1}{2} \} 
	\subset 
	\{ |x| \leq e^{ {e}^{M_6}}  \}. 
\]
Thus, 
\[
\begin{aligned}
	&N:=\int_{ |x| \leq e^{ {e}^{M_6}} } 
	|u(x,0)|^{q_c}  dx 
	\geq e^{-e^{e^{M_6}}} \log \frac{1}{-t'}. 
\end{aligned}
\]
This is equivalent to 
$-t' \geq \exp( -N e^{e^{e^{ M_6}}})$. 
Then, \eqref{eq:lowerboundconctime} with $x_0=0$ follows. 
On the other hand, 
from \eqref{eq:concloweee} and $T_2\leq T_\ma \leq M^{-2}$, 
it follows that 
\[
	N\geq 
	\int_{ T_2^\frac{1}{2} \leq |x| \leq e^{ {e}^{M_5}} T_2^\frac{1}{2}} 
	|u(x,0)|^{q_c}  dx 
	\geq 
	e^{-e^{e^{M_5}}} \geq  M e^{-e^{e^{M_6}}}. 
\]
This shows \eqref{eq:lowerboundconcN}. 
The proof of Proposition \ref{pro:quantestheat} is complete. 
\end{proof}

By using Proposition \ref{pro:quantestheat}, 
we show our main quantitative estimate.

\begin{proof}[Proof of Proposition \ref{pro:quantregheat}] 
Assume that $u$ satisfies the Lorentz norm bound \eqref{eq:Type1bound2.0} and
define $N$ by \eqref{eq:Ndef}. 
We consider $2$ cases. 

\noindent\textbf{Case 1}: $N<Me^{-e^{e^{M^{C_0}}}}$

Since \eqref{eq:lowerboundconcN} does not hold, 
the contraposition of Proposition \ref{pro:quantestheat} shows that 
\begin{equation}\label{eq:tpepup}
	(-t')^{\frac{2}{p-1}-\frac{n}{2}}
	\int_{t'}^{t'/2} 
	\int_{B(x_0, A(-t')^{\frac{1}{2}})} |u(x,t)|^{p+1} dxdt
	\leq \eps
\end{equation}
for all $t'\in (-M^{-p(2p+3)},0)$. 
Here $\eps=\eps_{0}/2$, $\eps_{0}$ is defined in Proposition \ref{pro:epsreg} 
and $A$ is defined in \eqref{eq:epsfixed}.
Now $t_*$ defined by \eqref{tstardef} satisfies 
\begin{equation}\label{eq:itsMp2p3}
	- t_*= M^{-p(2p+3)-1} \exp\left(-Ne^{e^{e^{M^{C_{0}}}}}\right) 
	\leq M^{-p(2p+3)-1} < M^{-p(2p+3)}. 
\end{equation}
By \eqref{eq:tpepup} with $t'= t_*$, we get 
\begin{equation}\label{eq:noconctstar}
	(-t_*)^{\frac{2}{p-1}-\frac{n}{2}}
	\int_{t_*}^{t_*/2} 
	\int_{B(x_0, A(-t_*)^{\frac{1}{2}})} |u(x,t)|^{p+1} dxdt
	\leq \eps.
\end{equation}
Then by Proposition \ref{pro:epsreg} 
with $t_0=0$ and  $\delta=(-t_*)^{\frac{1}{2}}$ 
and by $\eps<M^{-2(p+1)^2} \eps_1^{2(p+1)}$ 
with the fact that $\eps_1$ depends only on $n$ and $p$, 
we get 
\[
\begin{aligned}
	\|u\|_{L^{\infty}(B(x_0, \frac{1}{4}(-t_*)^{\frac{1}{2}})
	\times (\frac{1}{16}t_{*},0))}
	&\leq 
	C M \eps^\frac{1}{2(p+1)^2} 
	(-t_*)^{-\frac{1}{p-1}} \\
	&\leq 
	C \eps_1^\frac{1}{p+1} 
	(-t_*)^{-\frac{1}{p-1}}
	\leq 
	C(n,p) (-t_*)^{-\frac{1}{p-1}}. 
\end{aligned}
\]
Hence we obtain \eqref{eq:quantest}. 

\noindent\textbf{Case 2}: $N\geq Me^{-e^{e^{M^{C_0}}}}$

By \eqref{eq:itsMp2p3}, we have $t_*\in (-M^{-p(2p+3)},0)$. 
Since 
\[
	- t_*= M^{-p(2p+3)-1} \exp\left(-Ne^{e^{e^{M^{C_{0}}}}}\right) 
	< \exp\left(-Ne^{e^{e^{M^{C_{0}}}}}\right), 
\]
we observe that \eqref{eq:lowerboundconctime} does not hold. 
Then by the contraposition of Proposition \ref{pro:quantestheat}, 
we must have \eqref{eq:noconctstar}. 
By the same reasoning as Case 1, we get \eqref{eq:quantest}. 
\end{proof}

\section{Proof of main theorems}\label{sec:proof}
We are now in a position to prove 
our main theorems 
(Theorems \ref{thm:critrateheat}, \ref{thm:critratetypeIheat} and 
\ref{thm:critratetypeIheatinfinite}).

\begin{proof}[Proof of Theorem \ref{thm:critrateheat}] 
By applying Proposition \ref{pro:quantregheat} 
(specifically \eqref{eq:quantest}) to 
\[
	u_{t}(x,s):=(t^{\frac{1}{2}})^{\frac{2}{p-1}}u(t^{\frac{1}{2}} x,ts+t), 
\]
the following statement holds true.  
If $M_0$ is sufficiently large 
and if $u$ is a classical solution on $\R^n\times (0,t]$
with
\begin{equation}\label{eq:limsup}
	\max\left(M_0,\sup_{s\in (0,t)}\|u(\cdot,s)\|_{L^{q_c}(\R^n)}\right)
	\leq M_t 
\end{equation}
for some $M_t>0$, 
then there exists $C_*>0$ depending only on $n$ and $p$ such that
\begin{equation}\label{eq:prop4supest}
t^{\frac{1}{p-1}}\|u(\cdot,t)\|_{L^{\infty}(\R^n)}\leq e^{e^{e^{e^{ M_t^{C_*}}}}}.
\end{equation}
With this in hand, we proceed with proving Theorem \ref{thm:critrateheat}.

We suppose for contradiction that 
\[
	\limsup_{t\to T} 
	\frac{\|u(\cdot,t)\|_{L^{q_{c}}(\R^n)}}
	{\left( \log\log\log\log \left( \frac{1}{(T-t)^{ 1/(2(p-1))}}
	\right)\right)^{ 1/(2C_*)}}<\infty.
\]
This implies that there exists $N_*>0$ such that for all $t\in (0,T)$, 
\begin{equation}\label{eq:thm1contra2}
	\sup_{s\in (0,t)} {\|u(\cdot,s)\|_{L^{q_{c}}(\R^n)}}
	\leq N_{*} \left( 
	\log\log\log\log \left( \frac{1}{(T-t)^{ 1/(2(p-1))}} 
	\right) \right)^{ 1/(2C_*) }.
\end{equation}
Since $\lim_{t\to T} \log\log\log\log ( 1/(T-t)^{1/(2(p-1))} )=\infty$, 
there exists $0<\alpha<T$ depending only on 
$p$, $C_*$, $T$, $M_0$ and $N_*$ such that for all $t\in (T-\alpha, T)$, 
\begin{align}
	&\label{eq:alphaN*M0}
	\max(N_*^2, M_{0})
	\leq \left( \log\log\log\log \left( 
	\frac{1}{(T-t)^{ 1/(2(p-1))}} \right) \right)^{1/C_*}, \\
	&\label{eq:blowupratecontrapoint}
	t^{-\frac{1}{p-1}} (T-t)^{\frac{1}{2(p-1)}}
	<\frac{1}{2}(p-1)^{-\frac{1}{p-1}}.
\end{align}

Using \eqref{eq:thm1contra2} 
and \eqref{eq:alphaN*M0} gives that for all $t\in (T-\alpha, T)$, 
\[
\begin{aligned}
	&\begin{aligned}
	\sup_{s\in (0,t)} \|u(\cdot,s)\|_{L^{q_{c}}(\R^n)}
	&\leq N_{*} \left( 
	\log\log\log\log \left( 
	\frac{1}{(T-t)^{1/(2(p-1))}} \right)\right)^{1/(2C_*)}\\
	&\leq \left( \log\log\log\log \left( 
	\frac{1}{(T-t)^{1/(2(p-1))}} \right)\right)^{1/C_*}, 
	\end{aligned} \\
	&M_0\leq 
	\left( \log\log\log\log \left( 
	\frac{1}{(T-t)^{1/(2(p-1))}} \right)\right)^{1/C_*}.
\end{aligned}
\]
This shows that \eqref{eq:limsup} is satisfied with 
\[
	M_t:=\left( \log\log\log\log \left( 
	\frac{1}{(T-t)^{1/(2(p-1))}} \right)\right)^{1/C_*}.
\]
By substituting this into \eqref{eq:prop4supest} 
and using \eqref{eq:blowupratecontrapoint}, 
we have for all $t\in (T-\alpha, T)$, 
\begin{equation}\label{eq:thm1contra4}
	\|u(\cdot,t)\|_{L^\infty(\R^n)}
	\leq 
	t^{-\frac{1}{p-1}}(T-t)^{-\frac{1}{2(p-1)}}
	<\frac{(p-1)^{-\frac{1}{p-1}}}{2(T-t)^{\frac{1}{p-1}}}.
\end{equation}
But from \cite[Proposition 23.1]{QSbook2}, 
we have that for all $t\in (0,T)$, 
\begin{equation}\label{eq:QSlowpoint}
	\|u(\cdot,t)\|_{L^{\infty}(\R^n)}
	\geq \frac{(p-1)^{-\frac{1}{p-1}}}{(T-t)^{\frac{1}{p-1}}}. 
\end{equation}
This contradicts \eqref{eq:thm1contra4}.
Hence, we must have
\[
	\limsup_{t\to T} 
	\frac{ \|u(\cdot,t)\|_{L^{q_{c}}(\R^n)} }{\left(
	\log\log\log\log\left(\frac{1}{(T-t)^{1/(2(p-1))}}
	\right)\right)^{ 1/(2C_*)}}=\infty. 
\]
This shows the desired conclusion with $c:=1/(2C_*)$. 
The proof is complete. 
\end{proof}

\begin{proof}[Proof of Theorem \ref{thm:critratetypeIheat}]
Let 
\begin{equation}\label{eq:theorem2Trestrict}
	( 1+ M^{-2(2p^2+3p+2)})^{-1} T<t<T.
\end{equation}
We see that the assumptions of Theorem \ref{thm:critratetypeIheat} 
allow us to apply Proposition \ref{pro:quantregheat} 
(specifically \eqref{eq:quantest}) to 
\[
	u_{t}(x,s):=(t^{\frac{1}{2}})^{\frac{2}{p-1}}u(t^{\frac{1}{2}} x,ts+t).
\]
Then we can check that 
\[
\begin{aligned}
	t^{\frac{1}{p-1}}\|u(\cdot,t)\|_{L^{\infty}(\R^n)}
	&\leq C(n,p) M^\frac{p(2p+3)+1}{p-1} 
	\exp\left( \frac{1}{p-1} 
	e^{e^{e^{M^{C_0}}}} \int_{\R^n} |u(x,t)|^{q_c} dx\right) \\
	&\leq 
	M^\frac{p(2p+3)+2}{p-1} 
	\exp\left( \frac{1}{p-1} 
	e^{e^{e^{M^{C_0}}}} \int_{\R^n} |u(x,t)|^{q_c} dx\right)
\end{aligned}
\]
for $M$ sufficiently large. 
This together with \eqref{eq:QSlowpoint} gives 
\[
	(p-1)^{-\frac{1}{p-1}}
	\left( \frac{t}{T-t} \right)^{\frac{1}{p-1}} M^{-\frac{2p^2+3p+2}{p-1}}
	\leq \exp\left( 
	\frac{1}{p-1}e^{e^{e^{M^{C_0}}}} \int_{\R^n} |u(x,t)|^{q_c} dx\right).
\]
From \eqref{eq:theorem2Trestrict}, we see that
$1\leq (t/(T-t))^{1/(2(p-1))}  M^{-(2p^2+3p+2)/(p-1)}$. 
Thus,
\[
	(p-1)^{-\frac{1}{p-1}}
	\left( \frac{t}{T-t} \right)^{\frac{1}{2(p-1)}}
	\leq \exp\left( 
	\frac{1}{p-1}e^{e^{e^{M^{C_0}}}} 
	\int_{\R^n} |u(x,t)|^{q_{c}} dx \right), 
\]
and so 
\[
	\int_{\R^n} |u(x,t)|^{q_{c}} dx \geq 
	\frac{1}{2} e^{ - e^{e^{M^{C_0}} } } \log \left( \frac{t}{T-t} \right) 
	- e^{ - e^{e^{M^{C_0}} } } \log(p-1). 
\]
This together with \eqref{eq:theorem2Trestrict} 
gives the desired conclusion with $C:=2C_0$, for $M$ sufficiently large. 
\end{proof}

\begin{proof}[Proof of Theorem \ref{thm:critratetypeIheatinfinite}]
We see that the assumptions of Theorem \ref{thm:critratetypeIheatinfinite} 
allow us to apply Proposition \ref{pro:quantregheat} 
(specifically \eqref{eq:quantest}) to 
\[
	u_{t}(x,s):=(t^{\frac{1}{2}})^{\frac{2}{p-1}}u(t^{\frac{1}{2}} x,ts+t).
\]
This gives 
\begin{equation}\label{eq:uest1thm3}
	M^{-\frac{2p^2+3p+2}{p-1}} |u(0,t)| t^\frac{1}{p-1}
	\leq 
	\exp\left( \frac{1}{p-1} e^{e^{e^{M^{C_0}}}} 
	\int_{|x|\leq t^{\frac{1}{2}}{e^{e^{M^{C_0}}}}} |u(x,t)|^{q_{c}} dx \right).
\end{equation}
As (by assumption) $\limsup_{t\to \infty}|u(0,t)|=\infty$, 
there exists $t_{n}\to \infty$ such that
\[
	M^{-\frac{2p^2+3p+2}{p-1}} |u(0,t_{n})|\geq 1
	\quad\mbox{ for all }n\in\N. 
\]
Substituting this into \eqref{eq:uest1thm3} then readily gives the desired conclusion.
\end{proof}

\section*{Acknowledgments}
The first author thanks Philippe Souplet for references 
and the Institute of Science Tokyo for its hospitality. 
The second author was supported in part by JSPS KAKENHI 
Grant Numbers 17K05312, 21H00991 and 21H04433. 
The third author was supported in part 
by JSPS KAKENHI Grant Numbers 22H01131, 22KK0035 and 23K12998.

\end{document}